\documentclass[a4paper,12pt]{article}

\usepackage{amsmath}
\usepackage{ntheorem}
\usepackage{amssymb}
\usepackage{latexsym}
\usepackage{url,mathrsfs}
\usepackage[all]{xy}\UseComputerModernTips
\usepackage{graphicx}

\theoremstyle{plain}
\newtheorem{proposition}{Proposition}[section]
\newtheorem{theorem}[proposition]{Theorem}
\newtheorem{lemma}[proposition]{Lemma}
\newtheorem{corollary}[proposition]{Corollary}

\theoremstyle{plain}
\theorembodyfont{\normalfont\rm}
\newtheorem{definition}[proposition]{Definition}
\newtheorem{remark}[proposition]{Remark}
\newtheorem{example}[proposition]{Example}

\theoremstyle{nonumberplain}
\theoremheaderfont{\normalfont\itshape}
\theoremseparator{.}
\theorembodyfont{\normalfont}
\newtheorem{proof}{Proof}

\newcommand{\qed}{\hfill $\Box$}

\makeatletter
\renewcommand{\theequation}{%
\thesection.\arabic{equation}}
\@addtoreset{equation}{section}
\makeatother

\newcommand{\ZZ}{{\mathbb Z}}
\newcommand{\NN}{{\mathbb N}}
\newcommand{\RR}{{\mathbb R}}
\newcommand{\CC}{{\mathbb C}}
\newcommand{\PP}{{\mathbb P}}
\newcommand{\QQ}{{\mathbb Q}}
\newcommand{\LL}{{\mathbb L}}

\renewcommand{\d}{{\rm dim}}

\newcommand{\Vol}{{\rm Vol}}

\newcommand{\e}{\varepsilon}
\renewcommand{\SS}{{\mathcal S}}
\newcommand{\Spec}{{\rm Spec}}
\newcommand{\Con}{{\rm Con}}

\newcommand{\id}{{\rm id}}
\newcommand{\supp}{{\rm supp}}

\newcommand{\Int}{{\rm Int}}
\newcommand{\Ker}{{\rm Ker}}
\newcommand{\Var}{{\rm Var}}
\newcommand{\HSm}{{\rm HS}^{\rm mon}}
\renewcommand{\Im}{{\rm Im}}
\newcommand{\relint}{{\rm rel.int}}

\newcommand{\Cone}{{\rm Cone}}

\newcommand{\dist}{{\rm dist}}
\newcommand{\height}{{\rm ht}}

\newcommand{\dd}{b}

\newcommand{\Db}{{\bf D}^{b}}
\newcommand{\Dbc}{{\bf D}_{c}^{b}}

\newcommand{\KK}{{\rm K}}

\newcommand{\F}{{\cal F}}
\newcommand{\G}{{\cal G}}

\newcommand{\M}{{\mathcal M}}
\renewcommand{\O}{{\cal O}}

\newcommand{\CF}{{\rm CF}}
\newcommand{\J}{{S}}

\renewcommand{\div}{{\rm div}}

\renewcommand{\sp}{{\rm sp}}

\newcommand{\tl}[1]{\widetilde{#1}}

\newcommand{\simto}{\overset{\sim}{\longrightarrow}}

\newcommand{\dsum}{\displaystyle \sum}

\renewcommand{\(}{\left(}
\renewcommand{\)}{\right)}
\newcommand{\inun}{\text{\rotatebox[origin=c]{90}{$\in$}}}
\newcommand{\longhookrightarrow}{\DOTSB\lhook\joinrel\longrightarrow}
\newcommand{\longtwoheadrightarrow}{\relbar\joinrel\twoheadrightarrow}

\newcommand{\coker}{{\rm coker}}
\newcommand{\DR}{{\rm DR}}

\newcommand{\IC}{{\rm IC}}
\newcommand{\MHM}{{\rm MHM}}
\newcommand{\var}{{\rm var}}

\newcommand{\Clt}{\CC[t]\langle\partial_t\rangle}
\newcommand{\Cltau}{\CC[\tau]\langle\partial_\tau\rangle}
\newcommand{\Fou}[1]{{}^F\!#1}

\def\Afu{\mathbb{A}^{\!1}}
\def\DD{\mathbf{D}}
\def\bH{\boldsymbol{H}}

\def\cO{\mathscr{O}}
\def\cD{\mathscr{D}}
\def\cH{\mathscr{H}}
\def\cE{\mathscr{E}}
\def\cM{\mathscr{M}}

\def\bbA{{\mathbb A}}

\def\rH{\mathrm{H}}
\def\rM{\mathrm{M}}
\def\rN{\mathrm{N}}
\def\rP{\mathrm{P}}
\def\gr{\mathrm{gr}}

\let\wt\widetilde
\def\isom{\stackrel{\sim}{\longrightarrow}}
\def\to{\mathchoice{\longrightarrow}{\rightarrow}{\rightarrow}{\rightarrow}}
\def\mto{\mathchoice{\longmapsto}{\mapsto}{\mapsto}{\mapsto}}
\def\hto{\mathrel{\lhook\joinrel\to}}

\let\moins\smallsetminus
\def\bbullet{{\scriptscriptstyle\bullet}}

\def\cf{cf.\kern.3em}
\def\eg{e.g.\kern.3em}
\def\ie{i.e.,\ }

\newenvironment{enumeratea}
{\bgroup\begin{enumerate}}
{\end{enumerate}\egroup}

\renewcommand{\labelenumi}{{\rm (\roman{enumi})}}

\begin{document}

\title{Monodromy at infinity of polynomial maps and Newton polyhedra (with Appendix by C. Sabbah)}

\author{Yutaka \textsc{Matsui}\footnote{Department of Mathematics, Kinki University, 3-4-1, Kowakae, Higashi-Osaka, Osaka, 577-8502, Japan. E-mail: matsui@math.kindai.ac.jp} \and Kiyoshi \textsc{Takeuchi}\footnote{Institute of Mathematics, University  of Tsukuba, 1-1-1, Tennodai, Tsukuba, Ibaraki, 305-8571, Japan. E-mail: takemicro@nifty.com}}

\date{}

\sloppy

\maketitle

\begin{abstract}
By introducing motivic Milnor fibers at infinity of polynomial maps, we propose some methods for the study of nilpotent parts of monodromies at infinity. The numbers of Jordan blocks in the monodromy at infinity will be described by the Newton polyhedron at infinity of the polynomial.
\end{abstract}

\section{Introduction}\label{sec:1}

The aim of this paper is to study the nilpotent parts of monodromies at infinity of polynomial maps. More precisely, following the construction of motivic Milnor fibers in Denef-Loeser \cite{D-L-1} and \cite{D-L-2}, we introduce motivic reincarnations of global (Milnor) fibers of polynomial maps and give some methods for the calculations of their mixed Hodge numbers. Since by construction these mixed Hodge numbers contain the information on the monodromy at infinity of the map, we thereby determine its Jordan normal form. In particular, we will describe the numbers of Jordan blocks in the monodromy at infinity in terms of its Newton polyhedron at infinity.

In order to explain our results more precisely, we recall the definition and the basic properties of monodromies at infinity. After two fundamental papers \cite{Broughton} and \cite{S-T-1}, many authors studied the global behavior of polynomial maps $f \colon \CC^n \longrightarrow \CC$. For a polynomial map $f \colon \CC^n \longrightarrow \CC$, it is well-known that there exists a finite subset $B \subset \CC$ such that the restriction
\begin{equation}
\CC^n \setminus f^{-1}(B) \longrightarrow \CC \setminus B
\end{equation}
of $f$ is a locally trivial fibration. We denote by $B_f$ the smallest subset $B \subset \CC$ satisfying this condition. Let $C_R=\{x\in \CC\ |\ |x|=R\}$ ($R\gg 0$) be a sufficiently large circle in $\CC$ such that $B_f\subset \{x \in \CC\ |\ |x|<R\}$. Then by restricting the locally trivial fibration $\CC^n \setminus f^{-1}(B_f) \longrightarrow \CC \setminus B_f$ to $C_R$ we obtain a geometric monodromy automorphism $\Phi_f^{\infty} \colon f^{-1}(R) \simto f^{-1}(R)$ and the linear maps
\begin{equation}
\Phi_j^{\infty} \colon H^j(f^{-1}(R) ;\CC) \overset{\sim}{\longrightarrow} H^j(f^{-1}(R) ;\CC) \ \ (j=0,1,\ldots)
\end{equation}
associated to it, where the orientation of $C_R$ is taken to be counter-clockwise as usual. We call $\Phi_j^{\infty}$'s the (cohomological) monodromies at infinity of $f$. Various formulas for their eigenvalues (i.e. the semisimple parts) were obtained by Libgober-Sperber \cite{L-S} etc. Also, some important results on their nilpotent parts were obtained by Garc{\'i}a-L{\'o}pez-N{\'e}methi \cite{L-N-2} and Dimca-Saito \cite{D-S-1} etc. For example, Dimca-Saito \cite{D-S-1} obtained an upper bound of the sizes of Jordan blocks for the eigenvalue $1$ in $\Phi_j^{\infty}$. For the special case $n=2$, see \cite{Dimca2} etc. However, to the best of our knowledge, the nilpotent parts have not yet been fully understood. The monodromies at infinity $\Phi_j^{\infty}$ are important, because after a basic result \cite{N-N} of Neumann-Norbury, Dimca-N{\'e}methi \cite{D-N} proved that the monodromy representations
\begin{equation}
\pi_1(\CC \setminus B_f, c) \longrightarrow {\rm Aut} ( H^j(f^{-1}(c); \CC) ) \qquad (c \in \CC \setminus B_f)
\end{equation}
are completely determined by $\Phi_j^{\infty}$'s. In this paper, assuming that $f$ is convenient and non-degenerate at infinity (see Definition \ref{dfn:3-3}) we describe the nilpotent parts (i.e. the Jordan normal forms) explicitly. Note that the second condition is satisfied by generic polynomials $f(x)\in \CC [x_1, x_2, \ldots , x_n]$. By the results of Broughton \cite{Broughton}, $f$ is tame at infinity (see Definition \ref{dfn:tame}) and there exists a strong concentration $H^j(f^{-1}(R) ;\CC) \simeq 0$ ($j\neq 0, n-1$) of the cohomology groups of the generic fiber $f^{-1}(R)$ ($R \gg 0$) of $f$. Since $\Phi_0^{\infty}=\id_{\CC}$, $\Phi_{n-1}^{\infty}$ is the only non-trivial monodromy. Following \cite{L-S}, we call the convex hull of $\{0\}$ and the Newton polytope $NP(f)$ of $f$ in $\RR^n$ the Newton polyhedron of $f$ at infinity and denote it by $\Gamma_{\infty}(f)$. Let $q_1,\ldots,q_l$ (resp. $\gamma_1,\ldots, \gamma_{l^{\prime}}$) be the $0$-dimensional (resp. $1$-dimensional) faces of $\Gamma_{\infty}(f)$ such that $q_i\in \Int (\RR_+^n)$ (resp. the relative interior $\relint(\gamma_i)$ of $\gamma_i$ is contained in $\Int(\RR_+^n)$). For each $q_i$ (resp. $\gamma_i$), denote by $d_i >0$ (resp. $e_i>0$) its lattice distance $\dist(q_i, 0)$ (resp. $\dist(\gamma_i,0)$) from the origin $0\in \RR^n$. For $1\leq i \leq l^{\prime}$, let $\Delta_i$ be the convex hull of $\{0\}\sqcup \gamma_i$ in $\RR^n$. Then for $\lambda \in \CC \setminus \{1\}$ and $1 \leq i \leq l^{\prime}$ such that $\lambda^{e_i}=1$ we set
\begin{equation}
n(\lambda)_i
= \sharp\{ v\in \ZZ^n \cap \relint(\Delta_i) \ |\ \height (v, \gamma_i)=k\} +\sharp \{ v\in \ZZ^n \cap \relint(\Delta_i) \ |\ \height (v, \gamma_i)=e_i-k\},
\end{equation}
where $k$ is the minimal positive integer satisfying $\lambda=\zeta_{e_i}^{k}$ ($\zeta_{e_i}=\exp (2\pi\sqrt{-1}/e_i)$) and for $v\in \ZZ^n \cap \relint(\Delta_i)$ we denote by $\height (v, \gamma_i)$ the lattice height of $v$ from the base $\gamma_i$ of $\Delta_i$. Then in Section \ref{sec:7-2} we prove the following result which describes the number of Jordan blocks for each fixed eigenvalue $\lambda \neq 1$ in $\Phi_{n-1}^{\infty}$. Recall that by the monodromy theorem the sizes of such Jordan blocks are bounded by $n$.

\begin{theorem}\label{thm:1-1}
Assume that $f$ is convenient and non-degenerate at infinity. Then for any $\lambda \in\CC^* \setminus \{1\}$ we have
\begin{enumerate}
\item The number of the Jordan blocks for the eigenvalue $\lambda$ with the maximal possible size $n$ in $\Phi_{n-1}^{\infty} \colon H^{n-1}(f^{-1}(R);\CC) \simto H^{n-1}(f^{-1}(R);\CC)$ ($R \gg 0$) is equal to $\sharp \{q_i \ |\ \lambda^{d_i}=1\}$.
\item The number of the Jordan blocks for the eigenvalue $\lambda$ with size $n-1$ in $\Phi_{n-1}^{\infty}$ is equal to $\sum_{i \colon \lambda^{e_i}=1} n(\lambda)_i$.
\end{enumerate}
\end{theorem}

\noindent Namely the nilpotent parts for the eigenvalues $\lambda \neq 1$ in the monodromy at infinity $\Phi_{n-1}^{\infty}$ are determined by the lattice distances of the faces of $\Gamma_{\infty}(f)$ from the origin $0 \in \RR^n$. The monodromy theorem asserts also that the sizes of the Jordan blocks for the eigenvalue $1$ in $\Phi_{n-1}^{\infty}$ are bounded by $n-1$. In this case, we have the following result. Denote by $\Pi_f$ the number of the lattice points on the $1$-skeleton of $\partial \Gamma_{\infty}(f) \cap \Int (\RR^n_+)$. We say also that $\gamma \prec \Gamma_{\infty}(f)$ is a face at infinity of $\Gamma_{\infty}(f)$ if $0 \notin \gamma$. For a face at infinity $\gamma \prec \Gamma_{\infty}(f)$, denote by $l^*(\gamma)$ the number of the lattice points on the relative interior $\relint(\gamma)$ of $\gamma$.

\begin{theorem}\label{thm:1-2}
In the situation of Theorem \ref{thm:1-1} we have 
\begin{enumerate}
\item The number of the Jordan blocks for the eigenvalue $1$ with the maximal possible size $n-1$ in $\Phi_{n-1}^{\infty}$ is $\Pi_f$.
\item The number of the Jordan blocks for the eigenvalue $1$ with size $n-2$ in $\Phi_{n-1}^{\infty}$ is equal to $2 \sum_{\gamma} l^*(\gamma)$, where $\gamma$ ranges through the faces of $\Gamma_{\infty}(f)$ at infinity such that $\d \gamma =2$ and $\relint(\gamma) \subset \Int (\RR^n_+)$. In particular, this number is even.
\end{enumerate}
\end{theorem}

\noindent Roughly speaking, the nilpotent part for the eigenvalue $1$ in the monodromy at infinity $\Phi_{n-1}^{\infty}$ is determined by the convexity of the hypersurface $\partial \Gamma_{\infty}(f) \cap \Int (\RR^n_+)$. Thus Theorems \ref{thm:1-1} and \ref{thm:1-2} generalize the well-known fact that the monodromies of quasi-homogeneous polynomials are semisimple. Moreover we will also give a general algorithm for computing the numbers of Jordan blocks with smaller sizes. See Section \ref{sec:7-2} for the detail.

This paper is organized as follows. In Section \ref{sec:2}, after recalling some basic notions we give some generalizations of the results in Danilov-Khovanskii \cite{D-K} which will be effectively used later. In Section \ref{sec:3}, we recall some basic definitions on monodromies at infinity and review our new proof in \cite{M-T-new3} of Libgober-Sperber's theorem \cite{L-S} on the semisimple parts of monodromies at infinity. In Section \ref{sec:7}, we prove global analogues of the results in Denef-Loeser \cite{D-L-1} and \cite{D-L-2}. Namely by mimicking their construction, we introduce motivic Milnor fibers at infinity and prove basic results. Note that after our introduction of motivic Milnor fibers at infinity in the preliminary version arXiv:0809.3149v7 of \cite{M-T-new3} Raibaut \cite{Raibaut} introduced the same notion. Some deep results in Sabbah \cite{Sabbah}, \cite{Sabbah-2} and \cite{Sabbah-3} will be used to justify our arguments. Finally in Section \ref{sec:7-2}, by rewriting these results in terms of the Newton polyhedron at infinity $\Gamma_{\infty}(f)$ with the help of the results in Section \ref{sec:2}, we prove some combinatorial formulas for the Jordan normal form of the monodromy at infinity $\Phi_{n-1}^{\infty}$. We obtain also a global analogue of the Steenbrink conjecture proved by Varchenko-Khovanskii \cite{K-V} and Saito \cite{Saito-3}. In our another paper \cite{M-T-1} we apply our methods also to local Milnor monodromies and obtain results completely parallel to Theorems \ref{thm:1-1} and \ref{thm:1-2} etc. We thus find a striking symmetry between local and global. Note also that in the paper \cite{E-T} the results in these papers were generalized to the monodromies over complete intersection subvarieties in $\CC^n$. Moreover in \cite{M-T-2}, without assuming that $f$ is tame at infinity, we prove some general results on the upper bounds for the sizes and the numbers of Jordan blocks in the monodromies at infinity $\Phi_j^{\infty}$. 

\bigskip
\noindent{\bf Acknowledgement}: We thank Prof. Sch{\"u}rmann and Dr. Raibaut for pointing to us the fact that our motivic Milnor fiber at infinity $\SS_f^{\infty}$ of $f$ does not depend on the compactification of $\CC^n$ by \cite[Theorem 3.9]{G-L-M}. We are also grateful to Prof. Sabbah for several discussions and kindly permitting us to use his unpublished important results in this paper.

\section{Preliminary notions and results}\label{sec:2}

In this section, we introduce basic notions and results which will be used in this paper. In this paper, we essentially follow the terminology of \cite{Dimca}, \cite{H-T-T} and \cite{K-S}. For example, for a topological space $X$ we denote by $\Db(X)$ the derived category whose objects are bounded complexes of sheaves of $\CC_X$-modules on $X$. For an algebraic variety $X$ over $\CC$, let $\Dbc(X)$ be the full subcategory of $\Db(X)$ consisting of constructible complexes of sheaves. In this case, for an abelian group $G$ we denote by $\CF_G(X)$ the abelian group of $G$-valued constructible functions on $X$. Let $\CC(t)^*=\CC(t) \setminus \{0\}$ be the multiplicative group of the function field $\CC(t)$ of the scheme $\CC$. In this paper, we consider $\CF_G(X)$ only for $G=\ZZ$ or $\CC(t)^*$. For a $G$-valued constructible function $\rho \colon X \longrightarrow G$, by taking a stratification $X=\bigsqcup_{\alpha}X_{\alpha}$ of $X$ such that $\rho|_{X_{\alpha}}$ is constant for any $\alpha$, we set $\int_X \rho :=\sum_{\alpha}\chi(X_{\alpha}) \cdot \rho(x_{\alpha}) \in G$, where $x_{\alpha}$ is a reference point in $X_{\alpha}$. Then we can easily show that $\int_X\rho \in G$ does not depend on the choice of the stratification $X=\bigsqcup_{\alpha} X_{\alpha}$ of $X$. More generally, for any morphism $f \colon X \longrightarrow Y$ of algebraic varieties over $\CC$ and $\rho \in \CF_G(X)$, we define the push-forward $\int_f \rho \in \CF_G(Y)$ of $\rho$ by $(\int_f \rho)(y):=\int_{f^{-1}(y)} \rho$ for $y \in Y$. Now recall that for a non-constant regular function $f \colon X \longrightarrow \CC$ on a variety $X$ over $\CC$ and the hypersurface $X_0:= \{x\in X\ |\ f(x)=0\} \subset X$ there exists a nearby cycle functor
\begin{equation}
\psi_f \colon \Dbc(X) \longrightarrow \Dbc(X_0)
\end{equation}
defined by Deligne (see \cite[Section 4.2]{Dimca} for an excellent survey of this subject). As we see in the next proposition, the nearby cycle functor $\psi_f$ generalizes the classical notion of Milnor fibers. In the above situation, for $x \in X_0$ denote by $F_x$ the Milnor fiber of $f\colon X \longrightarrow \CC$ at $x$ (see for example \cite{Takeuchi} for a review on this subject). 

\begin{proposition}{\rm \bf(\cite[Proposition 4.2.2]{Dimca})}\label{prp:2-7-2} 
For any $\F\in \Dbc(X)$, $x \in X_0$ and $j \in \ZZ$, there exists a natural isomorphism
\begin{equation}\label{eq:2-24}
H^j(F_x ;\F) \simeq H^j(\psi_f(\F))_x.
\end{equation}
\end{proposition}

By this proposition, we can study the cohomology groups $H^j(F_x;\CC)$ of the Milnor fiber $F_x$ by using sheaf theory. Recall also that in the above situation we can define the Milnor monodromy operators
\begin{equation}
\Phi_{j,x} \colon H^j(F_x;\CC) \overset{\sim}{\longrightarrow} H^j(F_x;\CC) \ (j=0,1,\ldots)
\end{equation}
and the zeta function
\begin{equation}
\zeta_{f,x}(t):=\prod_{j=0}^{\infty} \det(\id -t\Phi_{j,x})^{(-1)^j} \in \CC(t)^*
\end{equation}
associated with it. This classical notion of Milnor monodromy zeta functions can be also generalized as follows.

\begin{definition}\label{dfn:2-8}
Let $f \colon X \longrightarrow \CC$ be a non-constant regular function on $X$ and $X_0 :=\{x\in X\ |\ f(x)=0\}$ the hypersurface defined by it. Then for $\F \in \Dbc(X)$ there exists a monodromy automorphism
\begin{equation}
\Phi(\F) \colon \psi_f(\F) \simto \psi_f(\F)
\end{equation}
of $\psi_f(\F)$ in $\Dbc(X_0)$ (see e.g. \cite[Section 4.2]{Dimca}). We define a $\CC(t)^*$-valued constructible function $\zeta_f(\F) \in \CF_{\CC(t)^*}(X_0)$ on $X_0$ by
\begin{equation}
\zeta_{f,x}(\F)(t):=\prod_{j \in \ZZ} \det\left(\id -t\Phi(\F)_{j,x}\right)^{(-1)^j}\in \CC(t)^*
\end{equation}
for $x \in X_0$, where $\Phi(\F)_{j,x} \colon (H^j(\psi_f(\F)))_x \simto (H^j(\psi_f(\F)))_x$ are induced by $\Phi(\F)$.
\end{definition}

For the proof of the following proposition, see for example, \cite[p.170-173]{Dimca}.

\begin{proposition}\label{prp:2-9}
Let $\pi \colon Y \longrightarrow X$ be a proper morphism of algebraic varieties over $\CC$ and $f \colon X \longrightarrow \CC$ a non-constant regular function on $X$. Set $g:=f \circ \pi \colon Y \longrightarrow \CC$, $X_0:=\{x\in X\ |\ f(x)=0\}$ and $Y_0:=\{y\in Y\ |\ g(y)=0\}=\pi^{-1}(X_0)$. Then for any $\G\in \Dbc(Y)$ we have
\begin{equation}
\int_{\pi|_{Y_0}} \zeta_g(\G) =\zeta_f(R\pi_*\G)
\end{equation}
in $\CF_{\CC(t)^*}(X_0)$, where $\int_{\pi|_{Y_0}}\colon \CF_{\CC(t)^*}(Y_0) \longrightarrow \CF_{\CC(t)^*}(X_0)$ is the push-forward of $\CC(t)^*$-valued constructible functions by $\pi|_{Y_0} \colon Y_0 \longrightarrow X_0$.
\end{proposition}

From now on, let us introduce our slight generalizations of the results in Danilov-Khovanskii \cite{D-K}.

\begin{definition}\label{dfn:2-6-2}
Let $g(x)=\sum_{v \in \ZZ^n} a_vx^v$ ($a_v\in \CC$) be a Laurent polynomial on $(\CC^*)^n$. 
\begin{enumerate}
\item We call the convex hull of $\supp(g):=\{v\in \ZZ^n \ |\ a_v\neq 0\} \subset \ZZ^n $ in $\RR^n$ the Newton polytope of $g$ and denote it by $NP(g)$.
\item For $u\in (\RR^n)^*$, we set $\Gamma(g;u):=\left\{ v\in NP(g)\ \left| \ \langle u,v\rangle =\min_{w\in NP(g)} \langle u,w\rangle \right.\right\}$.
\item For $u \in (\RR^n)^*$, we define the $u$-part of $g$ by $g^u(x):=\sum_{v \in \Gamma(g;u)} a_vx^v$.
\end{enumerate}
\end{definition}

\begin{definition}[\cite{Kushnirenko}]
Let $g$ be a Laurent polynomial on $(\CC^*)^n$. Then we say that the hypersurface $Z^*=\{ x\in (\CC^*)^n \ |\ g(x)=0 \}$ of $(\CC^*)^n$ is non-degenerate if for any $u \in (\RR^n)^*$ the hypersurface $\{ x\in (\CC^*)^n \ |\ g^u(x)=0 \}$ is smooth and reduced.
\end{definition}

In the sequel, let us fix an element $\tau =(\tau_1,\ldots, \tau_n) \in T:=(\CC^*)^n$ and let $g$ be a Laurent polynomial on $(\CC^*)^n$ such that $Z^*=\{ x\in (\CC^*)^n \ |\ g(x)=0\}$ is non-degenerate and invariant by the automorphism $l_{\tau} \colon (\CC^*)^n \underset{\tau \times}{\simto}(\CC^*)^n$ induced by the multiplication by $\tau$. Set $\Delta =NP(g)$ and for simplicity assume that $\d \Delta=n$. Then there exists $\beta \in \CC$ such that $l_{\tau}^*g= g \circ l_{\tau}=\beta g$. This implies that for any vertex $v$ of $\Delta =NP(g)$ we have ${\tau}^v={\tau}_1^{v_1} \cdots {\tau}_n^{v_n}=\beta$. Moreover by the condition $\d \Delta=n$ we see that $\tau_1, \tau_2, \ldots , \tau_n$ are roots of unity. For $p,q \geq 0$ and $k \geq 0$, let $h^{p,q}(H_c^k(Z^*;\CC))$ be the mixed Hodge number of $H_c^k(Z^*;\CC)$ and set
\begin{equation}
e^{p,q}(Z^*)=\dsum_k (-1)^k h^{p,q}(H_c^k(Z^*;\CC))
\end{equation}
as in \cite{D-K}. The above automorphism of $(\CC^*)^n$ induces a morphism of mixed Hodge structures $l_{\tau}^* \colon H_c^k(Z^*;\CC) \simto H_c^k(Z^*;\CC)$ and hence $\CC$-linear automorphisms of the $(p,q)$-parts $H_c^k(Z^*;\CC)^{p,q}$ of $H_c^k(Z^*;\CC)$. For $\alpha \in \CC$, let $h^{p,q}(H_c^k(Z^*;\CC))_{\alpha}$ be the dimension of the $\alpha$-eigenspace $H_c^k(Z^*;\CC)_{\alpha}^{p,q}$ of this automorphism of $H_c^k(Z^*;\CC)^{p,q}$ and set
\begin{equation}
e^{p,q}(Z^*)_{\alpha}=\dsum_k (-1)^k h^{p,q}(H_c^k(Z^*;\CC))_{\alpha}.
\end{equation}
Since we have $l_{\tau}^r =\id_{Z^*}$ for some $r \gg 0$, these numbers are zero unless $\alpha$ is a root of unity. Obviously we have
\begin{equation}
e^{p,q}(Z^*)=\dsum_{\alpha \in \CC} e^{p,q}(Z^*)_{\alpha}, \qquad e^{p,q}(Z^*)_{\alpha}=e^{q,p}(Z^*)_{\overline{\alpha}}.
\end{equation}
In this setting, along the lines of Danilov-Khovanskii \cite{D-K} we can give an algorithm for computing these numbers $e^{p,q}(Z^*)_{\alpha}$ as follows. First of all, as in \cite[Section 3]{D-K} we can easily obtain the following result.

\begin{proposition}\label{prp:2-15}
For $p,q \geq 0$ such that $p+q >n-1$, we have
\begin{equation}
e^{p,q}(Z^*)_{\alpha}=
\begin{cases}
(-1)^{n+p+1}\binom{n}{p+1} & (\text{$\alpha=1$ and $p=q$}),\\
\hspace*{10mm}0 & (\text{otherwise}),
\end{cases}
\end{equation}
(we used the convention $\binom{a}{b}=0$ ($0 \leq a <b$) for binomial coefficients).
\end{proposition}

\begin{proof}
If $p+q >n-1$, we have $H_c^k(Z^*;\CC)^{p,q}=0$ for $k \leq n-1$. Moreover for $p,q \geq 0$ such that $p+q>n-1$ and $k>n-1$ the Gysin homomorphism
\begin{equation}
\Theta_{p,q} \colon H_c^k(Z^*;\CC)^{p,q} \longrightarrow H_c^{k+2}((\CC^*)^n;\CC)^{p+1,q+1}
\end{equation}
is an isomorphism by \cite[Proposition 3.2]{D-K}. Since for such $p,q$ and $k$ there exists a commutative diagram
\begin{equation}
\xymatrix@C=30mm{
H_c^k(Z^*;\CC)^{p,q} \ar[r]_{\Theta_{p,q}}^{\sim} \ar[d]^{l_{\tau}^*}_{\wr}&H_c^{k+2}((\CC^*)^n;\CC)^{p+1,q+1} \ar[d]^{l_{\tau}^*}_{\wr}\\
H_c^k(Z^*;\CC)^{p,q} \ar[r]_{\Theta_{p,q}}^{\sim} & H_c^{k+2}((\CC^*)^n;\CC)^{p+1,q+1}
}
\end{equation}
and $l_{\tau} \colon (\CC^*)^n \simto (\CC^*)^n$ is homotopic to the identity of $(\CC^*)^n$, we obtain isomorphisms
\begin{equation}
H_c^k(Z^*;\CC)_{\alpha}^{p,q} \simeq H_c^{k+2}((\CC^*)^n;\CC)_{\alpha}^{p+1,q+1}=
\begin{cases}
\CC^{\binom{n}{p+1}} &(\text{$k=n+p-1$, $\alpha=1$ and $p=q$}), \\
0 & (\text{otherwise}).
\end{cases}
\end{equation}
Then the result follows from the definition of $e^{p,q}(Z^*)_{\alpha}$. This completes the proof. \qed
\end{proof}

For a vertex $w$ of $\Delta$, consider the translated polytope $\Delta^w:=\Delta -w$ such that $0 \prec \Delta^w$ and ${\tau}^v=1$ for any vertex $v$ of $\Delta^w$. Then for $\alpha \in \CC$ and $k \geq 0$ set
\begin{equation}
l^*(k\Delta)_{\alpha}=\sharp \{ v \in \Int (k\Delta^w) \cap \ZZ^n \ |\ {\tau}^v =\alpha\} \in \ZZ_+:=\ZZ_{\geq 0}
\end{equation}
and
\begin{equation}
l(k\Delta)_{\alpha}=\sharp \{ v \in (k\Delta^w) \cap \ZZ^n \ |\ {\tau}^v =\alpha\} \in \ZZ_+.
\end{equation}
We can easily see that these numbers $l^*(k\Delta)_{\alpha}$ and $l(k\Delta)_{\alpha}$ do not depend on the choice of the vertex $w$ of $\Delta$. Next, define two formal power series $P_{\alpha}(\Delta;t)=\sum_{i \geq 0}\varphi_{\alpha, i}(\Delta)t^i$ and $Q_{\alpha}(\Delta;t)=\sum_{i \geq 0}\psi_{\alpha,i}(\Delta)t^i$ by
\begin{equation}
P_{\alpha}(\Delta;t)=(1-t)^{n+1} \left\{ \dsum_{k \geq 0} l^*(k\Delta)_{\alpha}t^k\right\}
\end{equation}
and
\begin{equation}
Q_{\alpha}(\Delta;t)=(1-t)^{n+1} \left\{ \dsum_{k \geq 0} l(k\Delta)_{\alpha}t^k\right\}
\end{equation}
respectively. Then we can easily show that $P_{\alpha}(\Delta;t)$ is actually a polynomial as in \cite[Section 4.4]{D-K}. Moreover as in Macdonald \cite{Macdonald}, we can easily prove that for any $\alpha \in \CC^*$ the function $h_{\Delta,\alpha}(k):=l(k\Delta)_{\alpha^{-1}}$ of $k \geq 0$ is a polynomial of degree $n$ with coefficients in $\QQ$. By a straightforward generalization of the Ehrhart reciprocity proved by \cite{Macdonald}, we obtain also an equality
\begin{equation}
h_{\Delta,\alpha}(-k)=(-1)^n l^*(k\Delta)_{\alpha}
\end{equation}
for $k> 0$. By an elementary computation (see \cite[Remark 4.6]{D-K}), this implies that we have
\begin{equation}\label{E:sym}
\varphi_{\alpha , i}(\Delta)= \psi_{\alpha^{-1}, n+1-i}(\Delta ) \qquad (i\in\ZZ ).
\end{equation}
In particular, $Q_{\alpha}(\Delta;t)=\sum_{i \geq 0}\psi_{\alpha, i}(\Delta)t^i$ is a polynomial for any $\alpha \in \CC^*$.

\begin{theorem}\label{thm:2-14}
In the situation as above, we have
\begin{equation}
\dsum_q e^{p,q}(Z^*)_{\alpha}
=\begin{cases}
(-1)^{p+n+1}\binom{n}{p+1} +(-1)^{n+1} \varphi_{\alpha, n-p}(\Delta) & (\alpha=1), \\
(-1)^{n+1} \varphi_{\alpha, n-p}(\Delta) & (\alpha \neq 1).
\end{cases}
\end{equation}
\end{theorem}

\begin{proof}
Let $\Sigma_1$ be the dual fan of $\Delta$ in $\RR^n$. Then we can construct a subdivision $\Sigma$ of $\Sigma_1$ such that the toric variety $X_{\Sigma}$ associated with it is smooth and projective. Moreover, there exists a $T$-equivariant line bundle $\O_{X_{\Sigma}}(\Delta)$ on $X_{\Sigma}$ whose global section $\varGamma(X_{\Sigma};\O_{X_{\Sigma}}(\Delta))$ is naturally isomorphic to the space $\{\sum_{v \in \Delta \cap \ZZ^n}a_vx^v \ |\ a_v\in \CC\}$ of Laurent polynomials with support in $\Delta \cap \ZZ^n$ (see \cite[Section 2]{D-K} and \cite[Section 2.1]{Oda} etc.). Since the Laurent polynomial $g$ is a section of $\O_{X_{\Sigma}}(\Delta)$, we obtain an isomorphism $\O_{X_{\Sigma}}(\Delta)\simeq \O_{X_{\Sigma}}(\overline{Z^*})$. Then by using the isomorphism
\begin{equation}
\varGamma(X_{\Sigma};\O_{X_{\Sigma}}(\Delta)) \underset{{\bf A}}{\simto} \varGamma(X_{\Sigma}; \O_{X_{\Sigma}}(\overline{Z^*}))
\end{equation}
and the pull-back of the meromorphic functions in $\varGamma(X_{\Sigma}; \O_{X_{\Sigma}}(\overline{Z^*}))$ by $l_{\tau}$, we can define an action $l_{\tau}^*$ of $\tau \in (\CC^*)^n$ on $\varGamma(X_{\Sigma};\O_{X_{\Sigma}}(\Delta))\simeq \{\sum_{v \in \Delta \cap \ZZ^n} a_vx^v\ |\ a_v\in \CC\}$. Note that this action $l_{\tau}^*$ is different from the one constructed in \cite[Section 2.1 and 2.2]{Oda} by using the $T$-equivariance of the line bundle $\O_{X_{\Sigma}}(\Delta)$. From now on, we shall describe the action $l_{\tau}^*$ explicitly. For an $n$-dimensional cone $\sigma \in \Sigma$, let $v_{\sigma}\prec \Delta$ be the $0$-dimensional supporting face of $\sigma$ in $\Delta$ and $U_{\sigma} \simeq \CC_y^n$ the affine open subset of $X_{\Sigma}$ which corresponds to $\sigma$. More precisely we set $U_{\sigma}=\Spec(\CC[\sigma^{\vee} \cap \ZZ^n])$. Then on $U_{\sigma} \simeq \CC_y^n$ we have
\begin{equation}
g(y)=y_1^{a_1}\cdots y_n^{a_n}\times g_{\sigma}(y) \hspace{10mm}(a_i\in \ZZ),
\end{equation}
where $g_{\sigma}$ is a polynomial such that $NP(g_{\sigma}) =\Delta^{v_{\sigma}}=\Delta -v_{\sigma}$. Namely, in $U_{\sigma}$ the hypersurface $\overline{Z^*}\subset X_{\Sigma}$ is defined by $\overline{Z^*}=\{g_{\sigma}=0\}$. Hence there exists an isomorphism
\begin{equation}
\varGamma(U_{\sigma};\O_{X_{\Sigma}}) \underset{{\bf B}}{\simto} \varGamma(U_{\sigma};\O_{X_{\Sigma}}(\overline{Z^*}))
\end{equation}
given by
\begin{equation}
\sum_{v \in \sigma^{\vee} \cap \ZZ^n} a_vx^v \longmapsto \dfrac{1}{g_{\sigma}} \sum_{v \in \sigma^{\vee} \cap \ZZ^n} a_vx^v.
\end{equation}
Since we have $l_{\tau}^*g_{\sigma}=g_{\sigma}\circ l_{\tau}=g_{\sigma}$ by the construction of $g_{\sigma}$, via the isomorphism ${\bf B}$, the action $l_{\tau}^*$ of $\tau \in (\CC^*)^n$ on $\varGamma(U_{\sigma};\O_{X_{\Sigma}}(\overline{Z^*}))$ corresponds to the automorphism of $\varGamma(U_{\sigma};\O_{X_{\Sigma}})\simeq \CC[\sigma^{\vee} \cap \ZZ^n]$ defined by
\begin{equation}
\sum_{v\in \sigma^{\vee} \cap \ZZ^n} a_vx^v \longmapsto \sum_{v\in \sigma^{\vee} \cap \ZZ^n}a_v{\tau}^vx^v.
\end{equation}
On the other hand, there exists also a natural injection
\begin{equation}
\varGamma(X_{\Sigma};\O_{X_{\Sigma}}(\Delta)) \underset{{\bf C}}{\longhookrightarrow} \varGamma(U_{\sigma};\O_{X_{\Sigma}})
\end{equation}
given by
\begin{equation}
\sum_{v\in \Delta \cap \ZZ^n}a_vx^v \longmapsto \sum_{v \in \sigma^{\vee} \cap \ZZ^n} a_vx^{v-v_{\sigma}}.
\end{equation}
Then by the commutative diagram
\begin{equation}
\xymatrix@C=15mm{
\varGamma(X_{\Sigma};\O_{X_{\Sigma}}(\Delta)) \ar@{^{(}->}[r]^{\bf C} \ar[d]^{\bf A}_{\wr}& \varGamma(U_{\sigma};\O_{X_{\Sigma}}) \ar[d]^{\bf B}_{\wr}\\
\varGamma(X_{\Sigma};\O_{X_{\Sigma}}(\overline{Z^*})) \ar@{^{(}->}[r]& 
\varGamma(U_{\sigma};\O_{X_{\Sigma}}(\overline{Z^*})),}
\end{equation}
we see that the action $l_{\tau}^*$ on $\varGamma(X_{\Sigma};\O_{X_{\Sigma}}(\Delta))$ is given by
\begin{equation}
\sum_{v \in \Delta \cap \ZZ^n} a_vx^v \longmapsto \sum_{v\in \Delta \cap \ZZ^n} a_v{\tau}^{v-v_{\sigma}}x^v.
\end{equation}
Note that this morphism does not depend on the choice of the $n$-dimensional cone $\sigma \in \Sigma$. From now on, we will describe also a natural action of $\tau \in T=(\CC^*)^n$ on $\varGamma(X_{\Sigma};\O_{X_{\Sigma}}(k\Delta))$ for $k \geq 1$. For $k \geq 1$, let $g_k$ be a Laurent polynomial on $(\CC^*)^n$ such that $NP(g_k)=k\Delta$ and $Z^*_k=\{x\in (\CC^*)^n\ |\ g_k(x)=0\}$ is non-degenerate and stable by the automorphism $l_{\tau} \colon (\CC^*)^n \simto (\CC^*)^n$. Such a Laurent polynomial $g_k$ always exists (see Lemma \ref{lem:key} below). Since we have $\O_{X_{\Sigma}}(\Delta)^{\otimes k} \simeq \O_{X_{\Sigma}}(k\Delta)$, the $k$-th power $g^k$ of the Laurent polynomial $g$ is a section of $\O_{X_{\Sigma}}(k\Delta)$ satisfying the condition $\div g^k=k\overline{Z^*}$. Therefore we obtain isomorphisms
\begin{equation}
\O_{X_{\Sigma}}(k\Delta) \simeq \O_{X_{\Sigma}}(\overline{Z^*_k})\simeq \O_{X_{\Sigma}}(k\overline{Z^*})
\end{equation}
and the Weil divisors $\overline{Z^*_k}$ and $k\overline{Z^*}$ are naturally equivalent. Now let $\sigma \in \Sigma$ be an $n$-dimensional cone. Then on $U_{\sigma}\simeq \CC_y^n$ we have
\begin{eqnarray}
g^k(y)&=& y_1^{ka_1}\cdots y_n^{ka_n}\times g_{\sigma}^k(y),\\
g_k(y)&=& y_1^{ka_1}\cdots y_n^{ka_n}\times (g_k)_{\sigma}(y),
\end{eqnarray}
where $(g_k)_{\sigma}$ is an $l_{\tau}^*$-invariant polynomial on $U_{\sigma}\simeq \CC_y^n$. From this, we see that the rational function $\frac{g^k}{g_k}$ on $X_{\Sigma}$ is $l_{\tau}^*$-invariant and
\begin{equation}
\div \( \dfrac{g^k}{g_k}\)=k\overline{Z^*}-\overline{Z^*_k}
\end{equation}
on the whole $X_{\Sigma}$. Then there exists an isomorphism
\begin{equation}
\begin{array}{ccc}
\varGamma(X_{\Sigma};\O_{X_{\Sigma}}(k\overline{Z^*})) & \underset{{\bf D}}{\simto}& \varGamma(X_{\Sigma}; \O_{X_{\Sigma}}(\overline{Z^*_k}))\\
\inun & & \inun \\
\varphi & \longmapsto & \dfrac{g^k}{g_k}\times \varphi
\end{array}
\end{equation}
and a commutative diagram
\begin{equation}
\xymatrix{
\varGamma(X_{\Sigma};\O_{X_{\Sigma}}(k\overline{Z^*})) \ar[r]^{\sim}_{l_{\tau}^*} \ar[d]^{\wr}_{{\bf D}}& \varGamma(X_{\Sigma};\O_{X_{\Sigma}}(k\overline{Z^*}))\ar[d]^{\wr}_{{\bf D}}\\
\varGamma(X_{\Sigma};\O_{X_{\Sigma}}(\overline{Z^*_k})) \ar[r]^{\sim}_{l_{\tau}^*} \ar[d]^{\wr}& \varGamma(X_{\Sigma};\O_{X_{\Sigma}}(\overline{Z^*_k}))\ar[d]^{\wr}\\
\varGamma(X_{\Sigma};\O_{X_{\Sigma}}(k\Delta)) \ar[r]^{\sim} & \varGamma(X_{\Sigma};\O_{X_{\Sigma}}(k\Delta)),
}
\end{equation}
where the upper and middle horizontal arrows are the pull-backs of meromorphic functions on $X_{\Sigma}$ by $l_{\tau}$ and by taking a vertex $w$ of $k\Delta$ the bottom horizontal arrow is defined by
\begin{equation}
\sum_{v\in k\Delta \cap \ZZ^n} a_vx^v \longmapsto \sum_{v \in k\Delta \cap \ZZ^n} a_v \tau^{v-w}x^v.
\end{equation}
Moreover, let $D_1,\ldots, D_m$ be the (smooth) toric divisors on $X_{\Sigma}$ such that $X_{\Sigma} \setminus \( \bigcup_{i=1}^m D_i\)=(\CC^*)^n$ and for $0 \leq p \leq n$ set $D=\bigcup_{i=1}^m D_i$ and
\begin{equation}
\Omega_{(X_{\Sigma},D)}^p=\Ker \left[ \Omega_{X_{\Sigma}}^p \longrightarrow \bigoplus_{i=1}^m \Omega_{D_i}^p\right]
\end{equation}
as in \cite[Section 1.11]{D-K}. Then there exists a well-known isomorphism
\begin{equation}
\bigwedge^p\ZZ^n \otimes_{\ZZ} \O_{X_{\Sigma}}(-\sum_{i=1}^mD_i) \simto \Omega_{(X_{\Sigma},D)}^p
\end{equation}
given by
\begin{equation}
(v_1\wedge \cdots \wedge v_p) \otimes \varphi \longmapsto \varphi \times \dfrac{dx^{v_1}}{x^{v_1}} \wedge \cdots \wedge \dfrac{dx^{v_p}}{x^{v_p}},
\end{equation}
where $v_i \in \ZZ^n$ and $\varphi \in \O_{X_{\Sigma}}(-\sum_{i=1}^m 
D_i)\subset \O_{X_{\Sigma}}$. Since we have
\begin{equation}
l_{\tau}^*\( \dfrac{dx^{v_1}}{x^{v_1}} \wedge \cdots \wedge \dfrac{dx^{v_p}}{x^{v_p}}\) =\dfrac{dx^{v_1}}{x^{v_1}} \wedge \cdots \wedge \dfrac{dx^{v_p}}{x^{v_p}},
\end{equation}
for any $k\geq 1$ and $0 \leq p \leq n$ we obtain a commutative diagram
\begin{equation}
\xymatrix@R=5mm{
\varGamma(X_{\Sigma};\Omega_{(X_{\Sigma},D)}^p(k \overline{Z^*})) 
\ar[r]^{\sim}_{l_{\tau}^*} \ar@{-}[d]^{\wr}& \varGamma(X_{\Sigma};\Omega_{(X_{\Sigma},D)}^p (k \overline{Z^*}))\ar@{-}[d]^{\wr}\\
\bigwedge^p\ZZ^n \otimes_{\ZZ} \left\{ \sum_{v\in \Int(k\Delta) \cap \ZZ^n}a_vx^v \right\} \ar[r]^{\sim} & \bigwedge^p\ZZ^n \otimes_{\ZZ} \left\{ \sum_{v\in \Int(k\Delta) \cap \ZZ^n}a_vx^v \right\},
}
\end{equation}
where we set $\Omega_{(X_{\Sigma},D)}^p(k \overline{Z^*})=\Omega^p_{(X_{\Sigma},D)}\otimes_{\O_{X_{\Sigma}}}\O_{X_{\Sigma}}(k \overline{Z^*})$ and by taking a vertex $w$ of $k\Delta$ the bottom horizontal arrow is induced by
\begin{equation}
\sum_{v \in \Int(k\Delta) \cap \ZZ^n}a_vx^v \longmapsto \sum_{v \in \Int(k\Delta) \cap \ZZ^n} a_v\tau^{v-w}x^v.
\end{equation}
By using this explicit description of
\begin{equation}
l_{\tau}^* \colon \varGamma(X_{\Sigma}; \Omega_{(X_{\Sigma},D)}^p(k \overline{Z^*})) \simto \varGamma(X_{\Sigma};\Omega_{(X_{\Sigma},D)}^p(k \overline{Z^*})),
\end{equation}
the assertion can be proved just by following the proof for the formula in \cite[Section 4.4]{D-K}. This completes the proof. \qed
\end{proof}

With Proposition \ref{prp:2-15} and Theorem \ref{thm:2-14} at hands, we can now easily calculate the numbers $e^{p,q}(Z^*)_{\alpha}$ on the non-degenerate hypersurface $Z^* \subset (\CC^*)^n$ for any $\alpha \in \CC$ as in \cite[Section 5.2]{D-K}. Indeed for a projective toric compactification $X$ of $(\CC^*)^n$ such that the closure $\overline{Z^*}$ of $Z^*$ in $X$ is smooth, the variety $\overline{Z^*}$ is smooth projective and hence there exists a perfect pairing
\begin{equation}
H^{p,q}(\overline{Z^*};\CC)_{\alpha} \times H^{n-1-p, n-1-q}(\overline{Z^*};\CC)_{\alpha^{-1}} \longrightarrow \CC
\end{equation}
for any $p,q \geq 0$ and $\alpha \in \CC^*$ (see for example \cite[Section 5.3.2]{Voisin}). Therefore, we obtain equalities 
$e^{p,q}(\overline{Z^*})_{\alpha}=e^{n-1-p,n-1-q}(\overline{Z^*})_{\alpha^{-1}}$ which are necessary to proceed the algorithm in \cite[Section 5.2]{D-K}. We have also the following analogue of \cite[Proposition 5.8]{D-K}.

\begin{proposition}\label{prp:new}
For any $\alpha \in \CC$ and $p> 0$ we have
\begin{equation}
e^{p,0}(Z^*)_{\alpha}=e^{0,p}(Z^*)_{\overline{\alpha}}= (-1)^{n-1} \sum_{\begin{subarray}{c} \Gamma \prec \Delta\\ \d \Gamma =p+1\end{subarray}}l^*(\Gamma)_{\alpha}.
\end{equation}
\end{proposition}

The following result is an analogue of \cite[Corollary 5.10]{D-K}. For $\alpha \in \CC$, denote by $\Pi(\Delta)_{\alpha}$ the number of the lattice points $v=(v_1,\ldots, v_n)$ on the $1$-skeleton of $\Delta^w=\Delta-w$ such that ${\tau}^v=\alpha$, where $w$ is a vertex of $\Delta$.

\begin{proposition}\label{prp:2-19}
In the situation as above, for any $\alpha \in \CC^*$ we have
\begin{equation}
e^{0,0}(Z^*)_{\alpha}=
\begin{cases}
(-1)^{n-1} \left(\Pi(\Delta)_{1}-1\right) & (\alpha=1), \\
(-1)^{n-1}  \Pi(\Delta)_{\alpha^{-1}} & (\alpha \neq 1).
\end{cases}
\end{equation}
\end{proposition}

\begin{proof}
By Theorem \ref{thm:2-14}, Proposition \ref{prp:new} and the equality \eqref{E:sym}, the assertion can be proved as in the proof \cite[Corollary 5.10]{D-K}. \qed
\end{proof}

For a vertex $w$ of $\Delta$, we define a closed convex cone $\Con(\Delta, w)$ by $\Con(\Delta,w)=\{ r \cdot (v -w) \ |\ r \in \RR_+, \ v \in \Delta\} \subset \RR^n$.

\begin{definition}\label{dfn:2-16}
Let $\Delta$ be an $n$-dimensional integral polytope in $(\RR^n, \ZZ^n)$.
\begin{enumerate}
\item (see \cite[Section2.3]{D-K}) We say that $\Delta$ is prime if for any vertex $w$ of $\Delta$ the cone $\Con(\Delta,w)$ is generated by a basis of $\RR^n$.
\item We say that $\Delta$ is pseudo-prime if for any $1$-dimensional face $\gamma \prec \Delta$ the number of the $2$-dimensional faces $\gamma^{\prime} \prec \Delta$ such that $\gamma \prec \gamma^{\prime}$ is $n-1$.
\end{enumerate}
\end{definition}

By definition, prime polytopes are pseudo-prime. Moreover any face of a pseudo-prime polytope is again pseudo-prime.

\begin{definition}\label{dfn:7-6-4}(\cite{D-K}) Let $\Delta$ and $\Delta^{\prime}$ be two $n$-dimensional integral polytopes in $(\RR^n, \ZZ^n)$. We denote by ${\rm som}(\Delta )$ (resp. ${\rm som}(\Delta^{\prime})$) the set of vertices of $\Delta$ (resp. $\Delta^{\prime}$). Then we say that $\Delta^{\prime}$ majorizes $\Delta$ if there exists a map $\Psi \colon {\rm som}(\Delta^{\prime}) \longrightarrow {\rm som}(\Delta)$ such that $\Con(\Delta, \Psi(w)) \subset \Con(\Delta^{\prime}, w)$ for any vertex $w$ of $\Delta^{\prime}$.
\end{definition}

For an integral polytope $\Delta$ in $(\RR^n, \ZZ^n)$, we denote by $X_{\Delta}$ the toric variety associated with the dual fan of $\Delta$. Recall that if $\Delta^{\prime}$ majorizes $\Delta$ there exists a natural morphism $X_{\Delta^{\prime}} \longrightarrow X_{\Delta}$. 

\begin{proposition}\label{prp:7-6-5}
Let $\Delta$ and $Z^*_{\Delta}=Z^*$ with an action of $l_{\tau}$ be as above. Assume that an $n$-dimensional integral polytope $\Delta^{\prime}$ in $(\RR^n, \ZZ^n)$ majorizes $\Delta$ by the map $\Psi \colon {\rm som}(\Delta^{\prime}) \longrightarrow {\rm som}(\Delta)$. Then for the closure $\overline{Z^*}$ of $Z^*$ in $X_{\Delta^{\prime}}$ we have
\begin{eqnarray}
\sum_q e^{p,q}(\overline{Z^*})_1
&=&\sum_{\Gamma\prec \Delta^{\prime}} (-1)^{\d \Gamma+p+1} \left\{\binom{\d \Gamma}{p+1}-\binom{\dd_{\Gamma}}{p+1}\right\}\nonumber \\
& & +\sum_{\Gamma \prec \Delta^{\prime}}(-1)^{\d \Gamma +1}\sum_{i=0}^{\min\{\dd_{\Gamma},p\}}\binom{\dd_{\Gamma}}{i}(-1)^i \varphi_{1,\d \Psi(\Gamma)-p+i}(\Psi(\Gamma)),\label{eq:7-6-5}
\end{eqnarray}
where for $\Gamma \prec \Delta^{\prime}$ we set $\dd_{\Gamma}=\d \Gamma -\d \Psi(\Gamma)$.
\end{proposition}

\begin{proof}
We do not have to assume here that $\Delta^{\prime}$ is prime. Although the 
toric compactification $X_{\Delta^{\prime}}$ of $T=(\CC^*)^n$ might be very 
singular, it always admits the standard action of $T$. Hence $l_{\tau}: Z^* 
\simto Z^*$ naturally extends to an automorphism of the closure $\overline{Z^*}$ in $X_{\Delta^{\prime}}$. Let $X_{\Delta^{\prime}}=\bigsqcup_{\Gamma \prec \Delta^{\prime}} T_{\Gamma}$ be the decomposition of the toric variety $X_{\Delta^{\prime}}$ into $T$-orbits $T_{\Gamma} \simeq (\CC^*)^{\d \Gamma}$. Then by the non-degeneracy of $Z^*$ the hypersurfaces $\overline{Z^*} \cap T_{\Gamma} \subset T_{\Gamma}$ are also non-degenerate. Moreover we have the product decomposition  $\overline{Z^*} \cap T_{\Gamma} \simeq (\CC^*)^{\dd_{\Gamma}} \times 
Z_{\Psi(\Gamma)}^*$. Note that the action on its first component $(\CC^*)^{\dd_{\Gamma}}$ is homotopic to the identity. Therefore by Theorem \ref{thm:2-14} we have
\begin{eqnarray}
\sum_q e^{p,q}(\overline{Z^*})_1&=& \sum_q\sum_{\Gamma \prec \Delta^{\prime}} e^{p,q}((\CC^*)^{\dd_{\Gamma}} \times Z_{\Psi(\Gamma)}^*)_1\\
&=& \sum_{\Gamma \prec \Delta^{\prime}}\sum_{i=0}^{\min\{\dd_{\Gamma}, p\}} \binom{\dd_{\Gamma}}{i}(-1)^{i+\dd_{\Gamma}} \sum_qe^{p-i,q-i}(Z_{\Psi(\Gamma)}^*)_1\\
&=& \sum_{\Gamma \prec \Delta^{\prime}}\sum_{i=0}^{\min\{\dd_{\Gamma}, p\}} \binom{\dd_{\Gamma}}{i}(-1)^{i+\dd_{\Gamma}}\times (-1)^{\d \Psi(\Gamma)+1} \nonumber \\
& & \times \left\{ (-1)^{p-i}\binom{\d \Psi(\Gamma)}{p+1-i}+\varphi_{1,\d \Psi(\Gamma)-p+i}(\Psi(\Gamma))\right\}.
\end{eqnarray}
Then the result follows from the simple calculations 
\begin{equation}
\binom{\d \Gamma}{p+1}=\binom{\d \Psi(\Gamma)+\dd_{\Gamma}}{p+1}
= \sum_{i=0}^{\min\{ \dd_{\Gamma},p\}} \binom{\d \Psi(\Gamma)}{p+1-i}\binom{\dd_{\Gamma}}{i}+\binom{\dd_{\Gamma}}{p+1}.
\end{equation}
\qed
\end{proof}

From now on, we assume that $\Delta=NP(g)$ is pseudo-prime. Let $\Sigma$ be the dual fan of $\Delta$ and $X_{\Sigma}$ the toric variety associated to it. Then except finite points $X_{\Sigma}$ is an orbifold and the closure $\overline{Z^*}$ of $Z^*$ in $X_{\Sigma}$ does not intersect such points by the non-degeneracy of $g$. Hence $\overline{Z^*}$ is an orbifold i.e. quasi-smooth in the sense of \cite[Proposition 2.4]{D-K}. In particular, there exists a Poincar{\'e} duality isomorphism
\begin{equation}
[H^{p,q}(\overline{Z^*};\CC)_{\alpha}]^* \simeq H^{n-1-p,n-1-q}(\overline{Z^*};\CC)_{\alpha^{-1}}
\end{equation}
for any $\alpha \in \CC^*$ (see for example \cite{Danilov-2} and \cite[Corollary 8.2.22]{H-T-T}). Then by slightly generalizing the arguments in \cite{D-K} we obtain the following analogue of \cite[Section 5.5 and Theorem 5.6]{D-K}.

\begin{proposition}\label{prp:2-17}
In the situation as above, for any $\alpha \in \CC \setminus \{1\}$ and $p,q \geq 0$, we have
\begin{eqnarray}
e^{p,q}(\overline{Z^*})_{\alpha}
&=& \begin{cases}
\hspace*{5mm}\displaystyle -\sum_{\Gamma \prec \Delta}(-1)^{\d \Gamma} \varphi_{\alpha, \d \Gamma -p}(\Gamma) & (p+q=n-1),\\
\hspace*{20mm}0 & (\text{otherwise}),
\end{cases}\\
e^{p,q}(Z^*)_{\alpha}
&=& (-1)^{n+p+q} \sum_{\begin{subarray}{c} \Gamma \prec \Delta\\ \d \Gamma =p+q+1\end{subarray}} \left\{ \sum_{\Gamma^{\prime} \prec \Gamma} (-1)^{\d \Gamma^{\prime}} \varphi_{\alpha, \d \Gamma^{\prime}-p}(\Gamma^{\prime})\right\}.
\end{eqnarray}
\end{proposition}

\begin{proposition}\label{prp:2-17-2}
In the situation as above, we have
\begin{enumerate}
\item For $p,q \geq 0$ such that $p \neq q$, we have
\begin{equation}
e^{p,q}(\overline{Z^*})_1=
\begin{cases}
\hspace*{5mm}\displaystyle-\sum_{\Gamma \prec \Delta} (-1)^{\d \Gamma} \varphi_{1, \d \Gamma -\max\{p,q\}}(\Gamma) & (\text{$p+q=n-1$}),\\
\hspace*{25mm}0 & (\text{otherwise}).
\end{cases}
\end{equation}
\item For $p\geq 0$, we have
\begin{equation}
e^{p,p}(\overline{Z^*})_1=
\begin{cases}
\hspace*{5mm}\displaystyle(-1)^{p+1} \sum_{\begin{subarray}{c} \Gamma\prec \Delta \\ \d \Gamma \geq p+1\end{subarray}} (-1)^{\d \Gamma} \binom{\d \Gamma}{p+1} & (2p>n-1),\\
\hspace*{5mm}\displaystyle(-1)^{n-p} \sum_{\begin{subarray}{c} \Gamma\prec \Delta \\ \d \Gamma \geq n-p\end{subarray}} (-1)^{\d \Gamma} \binom{\d \Gamma}{n-p} & (2p<n-1),\\
\hspace*{5mm}\displaystyle \sum_{\Gamma\prec \Delta} (-1)^{\d \Gamma} \left\{(-1)^{p+1}\binom{\d \Gamma}{p+1} -\varphi_{1, \d \Gamma -p}(\Gamma) \right\} & (2p=n-1).
\end{cases}
\end{equation}
\end{enumerate}
\end{proposition}

From this proposition and the proof of \cite[Theorem 5.6]{D-K}, we obtain also the formula for $e^{p,q}(Z^*)_1$. For $\alpha \in \CC \setminus \{1\}$ and a face $\Gamma \prec \Delta$, set $\tl{\varphi}_{\alpha}(\Gamma)=\sum_{i=0}^{\d \Gamma} \varphi_{\alpha, i}(\Gamma)$. Then Proposition \ref{prp:2-17} can be rewritten as follows.

\begin{corollary}\label{cor:2-18}
For any $\alpha \in \CC \setminus \{1\}$ and $r \geq 0$, we have
\begin{equation}
\sum_{p+q=r}e^{p,q}(Z^*)_{\alpha}=(-1)^{n+r} \sum_{\begin{subarray}{c} \Gamma \prec \Delta\\ \d \Gamma =r+1\end{subarray}} \left\{ \sum_{\Gamma^{\prime} \prec \Gamma} (-1)^{\d \Gamma^{\prime}}\tl{\varphi}_{\alpha}(\Gamma^{\prime})\right\}.
\end{equation}
\end{corollary}

\section{Semisimple part of monodromies at infinity}\label{sec:3}

In this section, we recall some basic definitions on monodromies at infinity and review our new proof in \cite{M-T-new3} of Libgober-Sperber's theorem \cite{L-S}. Let $f(x)$ be a polynomial on $\CC^n$. Then as we explained in Introduction, there exist a locally trivial fibration $\CC^n \setminus f^{-1}(B_f) \longrightarrow \CC \setminus B_f$ and the linear maps 
\begin{equation}
\Phi_j^{\infty} \colon H^j(f^{-1}(R) ;\CC) \overset{\sim}{\longrightarrow} H^j(f^{-1}(R) ;\CC) \ \ (j=0,1,\ldots)
\end{equation}
($R \gg 0$) associated to it. To study the monodromies at infinity $\Phi_j^{\infty}$, we often impose the following natural condition.

\begin{definition}[\cite{Kushnirenko}]\label{dfn:tame}
Let $\partial f\colon \CC^n \longrightarrow \CC^n$ be the map defined by $\partial f(x)=(\partial_1f(x), \ldots, \partial_n f(x))$. Then we say that $f$ is tame at infinity if the restriction $(\partial f)^{-1}(B(0;\e )) \longrightarrow B(0;\e )$ of $\partial f$ to a sufficiently small ball $B(0;\e )$ centered at the origin $0 \in \CC^n$ is proper. 
\end{definition}

The following result is fundamental in the study of monodromies at infinity.

\begin{theorem}[Broughton \cite{Broughton} and Siersma-Tib{\u a}r \cite{S-T-1}]\label{tame}
Assume that $f$ is tame at infinity. Then the generic fiber $f^{-1}(c)$ ($c \in \CC \setminus B_f$) has the homotopy type of the bouquet of $(n-1)$-spheres. In particular, we have
\begin{equation}
H^j(f^{-1}(c);\CC)=0 \quad (j \neq 0, n-1).
\end{equation}
\end{theorem}

By this theorem if $f$ is tame at infinity, then $\Phi_{n-1}^{\infty}$ is the only non-trivial monodromy at infinity and its characteristic polynomial is calculated by the following zeta function $\zeta_f^{\infty}(t) \in \CC(t)^*$. 

\begin{definition}\label{dfn:3-4}
We define the monodromy zeta function at infinity $\zeta_f^{\infty}(t)$ of $f$ by
\begin{equation}
\zeta_{f}^{\infty}(t):=\prod_{j=0}^{\infty} \det(\id -t\Phi_j^{\infty})^{(-1)^j} \in \CC(t)^*.
\end{equation}
\end{definition}

\begin{definition}[\cite{L-S}]\label{dfn:3-1}
We call the convex hull of $\{0\}\cup NP(f)$ in $\RR^n$ the Newton polyhedron at infinity of $f$ and denote it by $\Gamma_{\infty}(f)$.
\end{definition}

For a subset $\J\subset \{1,2,\ldots,n\}$, let us set
\begin{equation}
\RR^{\J}:=\{ v=(v_1,v_2,\ldots,v_n)\in \RR^n\ |\ \text{$v_i=0$ for $i \notin \J$}\} \simeq \RR^{\sharp S}.
\end{equation}
We set also $\Gamma_{\infty}^{\J}(f)=\Gamma_{\infty}(f) \cap \RR^{\J}$. Recall that $f$ is convenient if we have $\dim \Gamma_{\infty}^{\J}(f)=\sharp {\J}$ for any ${\J} \subset \{1,2,\ldots,n\}$.

\begin{definition}[\cite{Kushnirenko}]\label{dfn:3-3}
We say that $f(x)=\sum_{v\in \ZZ_+^n} a_vx^v$ ($a_v\in \CC$) is non-degenerate at infinity if for any face $\gamma$ of $\Gamma_{\infty}(f)$ such that $0 \notin \gamma$ the complex hypersurface $\{x \in (\CC^*)^n\ |\ f_{\gamma}(x)=0\}$ in $(\CC^*)^n$ is smooth and reduced, where we set $f_{\gamma}(x)=\sum_{v \in \gamma \cap \ZZ_+^n} a_vx^v$.
\end{definition}

If $f$ is convenient and non-degenerate at infinity, then by a result of Broughton \cite{Broughton} it is tame at infinity. In this case, the monodromy zeta function $\zeta_f^{\infty}(t)$ has the following beautiful expression. For each non-empty subset ${\J} \subset \{1,2,\ldots, n\}$, let $\{\gamma_1^{\J},\gamma_2^{\J}, \ldots,\gamma_{n({\J})}^{\J}\}$ be the $(\sharp {\J}-1)$-dimensional faces of $\Gamma_{\infty}^{\J}(f)$ such that $0 \notin \gamma_i^{\J}$. For $1 \leq i \leq n(\J)$, let $u_i^{\J} \in (\RR^{\J})^* \cap \ZZ^{\J}$ be the unique non-zero primitive vector which takes its maximum in $\Gamma_{\infty}^{\J}(f)$ exactly on $\gamma_i^{\J}$ and set
\begin{equation}
d_i^{\J}: = \max_{v\in \Gamma_{\infty}^{\J}(f)} \langle u_i^{\J} ,v \rangle \in \ZZ_{>0}.
\end{equation}
We call $d_i^{\J}$ the lattice distance of $\gamma_i^{\J}$ from the origin $0 \in \RR^{\J}$. For each face $\gamma_i^{\J} \prec \Gamma_{\infty}^{\J}(f)$, let $\LL(\gamma_i^{\J})$ be the smallest affine linear subspace of $\RR^n$ containing $\gamma_i^{\J}$ and $\Vol_{\ZZ}(\gamma_i^{\J}) \in \ZZ_{>0}$ the normalized $(\sharp \J -1)$-dimensional volume (i.e. the $(\sharp S-1)!$ times the usual volume) of $\gamma_i^{\J}$ with respect to the lattice $\ZZ^n \cap \LL(\gamma_i^{\J})$.

\begin{theorem}\label{thm:3-5}{\rm (\cite{L-S}, see also \cite{M-T-new3} for a slight generalization )} 
Assume that $f$ is convenient and non-degenerate at infinity. Then we have 
\begin{equation}
\zeta_f^{\infty}(t)=\prod_{{\J} \neq \emptyset }\zeta^{\infty}_{f, {\J}}(t),
\end{equation}
where for each non-empty subset $\J \subset \{1,2,\ldots,n\}$ we set
\begin{equation}
\zeta^{\infty}_{f,{\J}}(t):=\prod_{i=1}^{n({\J})}(1-t^{d_i^{\J}})^{(-1)^{\sharp {\J}-1}\Vol_{\ZZ}(\gamma_i^{\J})}.
\end{equation}
\end{theorem}

This theorem was first proved by Libgober-Sperber \cite{L-S}. Here for the reader's convenience, we briefly recall our new proof in \cite{M-T-new3} which will be frequently used in this paper.

\begin{proof}
Let $j \colon \CC \longhookrightarrow \PP^1=\CC\sqcup \{\infty\}$ be the 
compactification and set $\F:=j_!(Rf_!\CC_{\CC^n})\in \Dbc(\PP^1)$. Take a local coordinate $h$ of $\PP^1$ in a neighborhood of $\infty \in \PP^1$ such that $\infty=\{h=0\}$. Then by the isomorphism $H_j(f^{-1}(R);\CC) \simeq H_c^{2n-2-j}(f^{-1}(R);\CC)$ we see that
\begin{equation}
\zeta_f^{\infty}(t)=\zeta_{h, \infty}(\F)(t) \in \CC(t)^*.
\end{equation}
Now let us consider $\CC^n$ as a toric variety associated with the fan $\Sigma_0$ in $\RR^n$ formed by the all faces of the first quadrant $\RR_+^n:=(\RR_{\geq 0})^n\subset \RR^n$. Let $T\simeq (\CC^*)^n$ be the open dense torus in it. Then by the convenience of $f$, $\Sigma_0$ is a subfan of the dual fan $\Sigma_1$ of $\Gamma_{\infty}(f)$ and we can construct a smooth subdivision $\Sigma$ of $\Sigma_1$ without subdividing the cones in $\Sigma_0$ (see e.g. \cite[Lemma (2.6), Chapter II, page 99]{Oka}). This implies that the toric variety $X_{\Sigma}$ associated with $\Sigma$ is a smooth compactification of $\CC^n$. Recall that $T$ acts on $X_{\Sigma}$ and the $T$-orbits are parametrized by the cones in $\Sigma$. For a cone $\sigma \in \Sigma$ denote by $T_{\sigma} \simeq (\CC^*)^{n-\dim \sigma}$ the corresponding $T$-orbit. We have also natural affine open subsets $\CC^n(\sigma) \simeq \CC^n$ of $X_{\Sigma}$ associated to $n$-dimensional cones $\sigma$ in $\Sigma$. Let $\sigma$ be an $n$-dimensional cone in $\Sigma$ and $\{w_1,\ldots, w_n\} \subset \ZZ^n$ the set of the primitive vectors on the edges of $\sigma$. Then there exists an affine open subset $\CC^n(\sigma)$ of $X_{\Sigma}$ such that $\CC^n(\sigma) \simeq \CC^n_y$ and $f$ has the following form on it:
\begin{equation} 
f(y)=\sum_{v \in \ZZ_+^n}a_{v}y_1^{\langle w_1,v \rangle}\cdots y_n^{\langle w_n, v \rangle} =y_1^{b_1} \cdots y_n^{b_n} \times f_{\sigma}(y), 
\end{equation}
where we set $f=\sum_{v \in \ZZ_+^n}a_{v}x^{v}$,
\begin{equation}
b_i=\min_{v\in \Gamma_{\infty}(f)} \langle w_i,v \rangle \leq 0 \qquad (i=1,2,\ldots,n)
\end{equation}
and $f_{\sigma}(y)$ is a polynomial on $\CC^n(\sigma) \simeq \CC^n_y$. In $\CC^n(\sigma) \simeq \CC^n_y$ the hypersurface $Z:= \overline{f^{-1}(0)} \subset X_{\Sigma}$ is explicitly written as $\{ y \in \CC^n(\sigma) \ | \ f_{\sigma}(y)=0 \}$. The variety $X_{\Sigma}$ is covered by such affine open subsets. Let $\tau$ be a $d$-dimensional face of the $n$-dimensional cone $\sigma \in \Sigma$. For simplicity, assume that $w_1,\ldots, w_d$ generate $\tau$. Then in the affine chart $\CC^n(\sigma) \simeq \CC^n_y$ the $T$-orbit $T_{\tau}$ associated to $\tau$ is explicitly defined by
\begin{equation*}
T_{\tau}=\{(y_1,\ldots,y_n)\in \CC^n(\sigma)\  |\ y_1=\cdots =y_d=0,\ y_{d+1},\ldots, y_n\neq 0\}\simeq (\CC^*)^{n-d}.
\end{equation*}
Hence we have
\begin{equation}
 X_{\Sigma}=\bigcup_{\dim \sigma =n}\CC^n(\sigma)=\bigsqcup_{\tau \in\Sigma}T_{\tau}.
\end{equation}
Now $f$ was extended to a meromorphic function $\tl{f}$ on $X_{\Sigma}$, but $\tl{f}$ has still points of indeterminacy. From now on, we will eliminate such points by blowing up $X_{\Sigma}$. For a cone $\sigma$ in $\Sigma$ by taking a non-zero vector $u$ in the relative interior $\relint(\sigma)$ of $\sigma$ we define a face $\gamma(\sigma)$ of $\Gamma_{\infty}(f)$ by
\begin{equation}
\gamma(\sigma) =\left\{ v \in \Gamma_{\infty}(f) \ \Big| \ \langle u ,v \rangle = \min_{w \in \Gamma_{\infty}(f)} \langle u,w \rangle \right\}.
\end{equation}
This face $\gamma(\sigma)$ does not depend on the choice of $u \in \relint(\sigma)$ and is called the supporting face of $\sigma$ in $\Gamma_{\infty}(f)$. Following \cite{L-S}, we say that a $T$-orbit $T_{\sigma}$ in $X_{\Sigma}$ is at infinity if $0 \notin \gamma(\sigma)$. In our situation (i.e. $f$ is convenient), this is equivalent to the condition $\sigma \not\subset \RR^n_+$. We can easily see that $\tl{f}$ has poles on the union of $T$-orbits at infinity. Let $\rho_1, \rho_2, \ldots, \rho_m$ be the $1$-dimensional cones in $\Sigma$ such that $\rho_i \not\subset \RR^n_+$ and set $T_i=T_{\rho_i}$. Then $T_1,T_2,\ldots, T_m$ are the $(n-1)$-dimensional $T$-orbits at infinity in $X_{\Sigma}$. For any $i=1,2,\ldots,m$ the toric divisor $D_i:=\overline{T_i}$ is a smooth hypersurface in $X_{\Sigma}$ and the poles of $\tl{f}$ are contained in $D_1 \cup \cdots \cup D_m$. Moreover by the non-degeneracy at infinity of $f$, the hypersurface $Z= \overline{f^{-1}(0)}$ in $X_{\Sigma}$ intersects $D_I:= \bigcap_{i \in I}D_i$ transversally for any non-empty subset $I \subset \{ 1,2, \ldots, m \}$. At such intersection points, $\tl{f}$ has indeterminacy. Furthermore we denote the (unique non-zero) primitive vector in $\rho_i \cap \ZZ^n$ by $u_i$. Then the order $a_i>0$ of the pole of $\tl{f}$ along $D_i$ is given by
\begin{equation}
a_i=-\min_{v\in \Gamma_{\infty}(f)} \langle u_i,v \rangle.
\end{equation} 
Now, in order to eliminate the indeterminacy of the meromorphic function $\tl{f}$ on $X_{\Sigma}$, we first consider the blow-up $\pi_1 \colon X_{\Sigma}^{(1)} \longrightarrow X_{\Sigma}$ of $X_{\Sigma}$ along the $(n-2)$-dimensional smooth subvariety $D_1\cap Z$. Then the indeterminacy of the pull-back $\tl{f}\circ \pi_1$ of $\tl{f}$ to $X_{\Sigma}^{(1)}$ is improved. If $\tl{f}\circ \pi_1$ still has points of indeterminacy on the intersection of the exceptional divisor $E_1$ of $\pi_1$ and the proper transform $Z^{(1)}$ of $Z$, we construct the blow-up $\pi_2 \colon X_{\Sigma}^{(2)} \longrightarrow X_{\Sigma}^{(1)}$ of $X_{\Sigma}^{(1)}$ along $E_1 \cap Z^{(1)}$. By repeating this procedure $a_1$ times, we obtain a tower of blow-ups
\begin{equation}
X_{\Sigma}^{(a_1)} \underset{\pi_{a_1}}{\longrightarrow} 
\cdots \cdots
\underset{\pi_2}{\longrightarrow} X_{\Sigma}^{(1)} 
\underset{\pi_1}{\longrightarrow} X_{\Sigma}.
\end{equation}
Then the pull-back of $\tl{f}$ to $X_{\Sigma}^{(a_1)}$ has no indeterminacy over $T_1$ (see the figures below).

\noindent\begin{minipage}{0.17\textwidth}
\begin{center}
\includegraphics[scale=0.85]{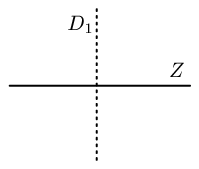}

Figure 1
\end{center}
\end{minipage}
\begin{minipage}{0.27\textwidth}
\begin{center}

\includegraphics[scale=.68]{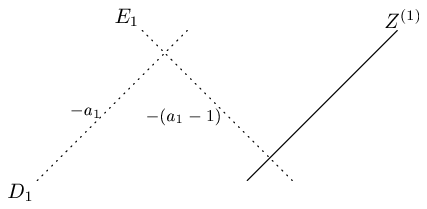}

Figure 2
\end{center}
\end{minipage}
\begin{minipage}{0.53\textwidth}
\begin{center}
\includegraphics[scale=.68]{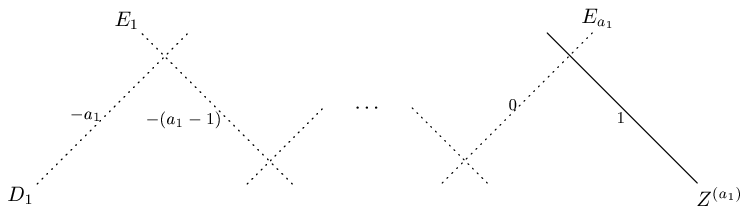}

Figure 3
\end{center}
\end{minipage}

\vspace{1mm}
Next we apply this construction to the proper transforms of $D_2$ and $Z$ in $X_{\Sigma}^{(a_1)}$. Then we obtain also a tower of blow-ups
\begin{equation}
X_{\Sigma}^{(a_1)(a_2)} \longrightarrow \cdots \cdots \longrightarrow X_{\Sigma}^{(a_1)(1)} \longrightarrow X_{\Sigma}^{(a_1)}
\end{equation}
and the indeterminacy of the pull-back of $\tl{f}$ to $X_{\Sigma}^{(a_1)(a_2)}$ is eliminated over $T_1 \sqcup T_2$. By applying the same construction to (the proper transforms of) $D_3, D_4,\ldots, D_m$, we finally obtain a birational morphism $\pi \colon \tl{X_{\Sigma}} \longrightarrow X_{\Sigma}$ such that $g:=\tl{f} \circ \pi$ has no point of indeterminacy on the whole $\tl{X_{\Sigma}}$. Note that the smooth compactification $\tl{X_{\Sigma}}$ of $\CC^n$ thus obtained is not a toric variety any more. We shall explain the geometry of $\tl{X_{\Sigma}}$ more precisely in the proof of Lemma \ref{geometry}. In particular, we will see that the union of the exceptional divisors of $\pi \colon \tl{X_{\Sigma}} \longrightarrow X_{\Sigma}$ and the proper transforms of $D_1, \ldots, D_m$ in $\tl{X_{\Sigma}}$ is normal crossing. Finally we get a commutative diagram of holomorphic maps
\begin{equation}
\xymatrix{
\CC^n \ar@{^{(}->}[r]^{\iota} \ar[d]_f & \tl{X_{\Sigma}} \ar[d]^g\\
\CC \ar@{^{(}->}[r]^j & \PP^1,}
\end{equation}
where $g$ is proper. Therefore we obtain an isomorphism $\F = j_!(Rf_!\CC_{\CC^n}) \simeq Rg_* (\iota_! \CC_{\CC^n})$ in $\Dbc(\PP^1)$. Let us apply Proposition \ref{prp:2-9} to the proper morphism $g \colon \tl{X_{\Sigma}} \longrightarrow \PP^1$. Then by calculating the monodromy zeta function of $\psi_{h \circ g}(\iota_! \CC_{\CC^n})$ at each point of $(h\circ g)^{-1}(0) =g^{-1}(\infty ) \subset \tl{X_{\Sigma}}$, we can calculate $\zeta_{h, \infty}(\F)(t)$ with the help of Bernstein-Khovanskii-Kushnirenko's theorem (see \cite{Khovanskii} etc.). This completes the proof.\qed
\end{proof}

\section{Motivic Milnor fibers at infinity}\label{sec:7}

In this section, following Denef-Loeser \cite{D-L-1} and \cite{D-L-2} we introduce motivic reincarnations of global (Milnor) fibers of polynomial maps and give a general formula for the nilpotent parts (i.e. the numbers of Jordan blocks of arbitrary sizes) in their monodromies at infinity. Namely, we formulate a global analogue of the results in \cite{D-L-1} and \cite{D-L-2}. Let $f \colon \CC^n \longrightarrow \CC$ be a polynomial map. We take a smooth compactification $X$ of $\CC^n$. Then by eliminating the points of indeterminacy of the meromorphic extension of $f$ to $X$ we obtain a commutative diagram
\begin{equation}
\xymatrix{
\CC^n \ar@{^{(}->}[r]^{\iota} \ar[d]_f & \tl{X} \ar[d]^g
\\
\CC \ar@{^{(}->}[r]^j & \PP^1}
\end{equation}
such that $g$ is a proper holomorphic map and $\tl{X} \setminus \CC^n$, $Y:=g^{-1}( \infty )$ are normal crossing divisors in $\tl{X}$. Take a local coordinate $h$ of $\PP^1$ in a neighborhood of $\infty\in \PP^1$ such that $\infty=\{h=0\}$ and set $\tl{g}=h\circ g$. Note that $\tl{g}$ is a holomorphic function defined on a neighborhood of the closed subvariety $Y=\tl{g}^{-1}(0)=g^{-1}(\infty) \subset \tl{X} \setminus \CC^n$ of $\tl{X}$. Then for $R \gg 0$ we have
\begin{equation}\label{eq:4-2}
H_c^j(f^{-1}(R);\CC) \simeq H^j(Y; \psi_{\tl{g}}(\iota_! \CC_{\CC^n})).
\end{equation}
Let us define an open subset $\Omega$ of $\tl{X}$ by
\begin{equation}
\Omega=\Int (\iota(\CC^n) \sqcup Y)
\end{equation}
and set $U=\Omega \cap Y$. Then $U$ (resp. the complement of $\Omega$ in $\tl{X}$) is a normal crossing divisor in $\Omega$ (resp. $\tl{X}$). Hence we can easily prove the isomorphisms
\begin{equation}
H^j(Y;\psi_{\tl{g}}(\iota_! \CC_{\CC^n}))\simeq H^j(Y; \psi_{\tl{g}}(\iota_!^{\prime}\CC_{\Omega}))\simeq H_c^j(U; \psi_{\tl{g}}(\CC_{\tl{X}})),
\label{eq:4-5}
\end{equation}
where $\iota^{\prime} \colon \Omega \longhookrightarrow \tl{X}$ is the inclusion. Now let $E_1, E_2, \ldots, E_k$ be the irreducible components of the normal crossing divisor $U=\Omega \cap Y$ in $\Omega \subset \tl{X}$. For each $1 \leq i \leq k$, let $b_i>0$ be the order of the zero of $\tl{g}$ along $E_i$. For a non-empty subset $I \subset \{1,2,\ldots, k\}$, let us set
\begin{equation}
E_I=\bigcap_{i \in I} E_i,\hspace{10mm}E_I^{\circ}=E_I \setminus \bigcup_{i \not\in I}E_i
\end{equation}
and $d_I=\gcd (b_i)_{i \in I}>0$. Then, as in \cite[Section 3.3]{D-L-2}, we can construct an unramified Galois covering $\tl{E_I^{\circ}} \longrightarrow E_I^{\circ}$ of $E_I^{\circ}$ as follows. First, for a point $p \in E_I^{\circ}$ we take an affine open neighborhood $W \subset \Omega \setminus ( \cup_{i \notin I} E_i)$ of $p$ on which there exist regular functions $\xi_i$ ($i\in I$) such that $E_i \cap W=\{ \xi_i=0 \}$ for any $i \in I$. Then on $W$ we have $\tl{g}=\tl{g_{1,W}} (\tl{g_{2,W}})^{d_I}$, where we set $\tl{g_{1,W}}=\tl{g} \prod_{i \in I}\xi_i^{-b_i}$ and $\tl{g_{2,W}}=\prod_{i \in I} \xi_i^{\frac{b_i}{d_I}}$. Note that $\tl{g_{1,W}}$ is a unit on $W$ and $\tl{g_{2,W}} \colon W \longrightarrow \CC$ is a regular function. It is easy to see that $E_I^{\circ}$ is covered by such affine open subsets $W$ of $\Omega \setminus ( \cup_{i \notin I} E_i)$. Then as in \cite[Section 3.3]{D-L-2} by gluing the varieties
\begin{equation}\label{eq:6-26}
\tl{E_{I,W}^{\circ}}=\{(t,z) \in \CC^* \times (E_I^{\circ} \cap W) \ |\ t^{d_I} =(\tl{g_{1,W}})^{-1}(z)\}
\end{equation}
together in the following way, we obtain the variety $\tl{E_I^{\circ}}$ over $E_I^{\circ}$. If $W^{\prime}$ is another such open subset and $\tl{g}=\tl{g_{1,W^{\prime}}} (\tl{g_{2,W^{\prime}}})^{d_I}$ is the decomposition of $\tl{g}$ on it, we patch $\tl{E_{I,W}^{\circ}}$ and $\tl{E_{I,W^{\prime}}^{\circ}}$ by the morphism $(t,z) \longmapsto (\tl{g_{2,W^{\prime}}}(z)( \tl{g_{2,W}})^{-1}(z) \cdot t, z)$ defined over $W \cap W^{\prime}$.

\begin{remark}
Let $N>0$ be the least common multiple of $b_1, \ldots, b_k$. As in Steenbrink \cite{Steenbrink}, by taking the normalization of the base change of $\tl{g}\colon \Omega \longrightarrow \CC$ by the $N$-th power map $\CC \longrightarrow \CC$ we obtain a morphism $\Omega^{\prime} \longrightarrow \Omega$. Then it is well-known that the variety $\tl{E_I^{\circ}}$ is obtained as a connected component of the inverse image of $E_I^{\circ}$ by $\Omega^{\prime} \longrightarrow \Omega$ (see Looijenga \cite{Looijenga}). Moreover Looijenga \cite[Lemma 5.3]{Looijenga} proved that $\tl{E_I^{\circ}} \longrightarrow E_I^{\circ}$ is the Stein factorization of a fiber bundle over $E_I^{\circ}$, which shows why the term $(1-\LL)^{\sharp I -1}$ appears in \eqref{MMF}. 
\end{remark}

Now for $d \in \ZZ_{>0}$, let $\mu_d \simeq \ZZ/\ZZ d$ be the multiplicative group consisting of the $d$-roots in $\CC$. We denote by $\hat{\mu}$ the projective limit $\underset{d}{\varprojlim} \mu_d$ of the projective system $\{ \mu_i \}_{i \geq 1}$ with morphisms $\mu_{id} \longrightarrow \mu_i$ given by $t \longmapsto t^d$. Then the unramified Galois covering $\tl{E_I^{\circ}}$ of $E_I^{\circ}$ admits a natural $\mu_{d_I}$-action defined by assigning the automorphism $(t,z) \longmapsto (\zeta_{d_I} t, z)$ of $\tl{E_I^{\circ}}$ to the generator $\zeta_{d_I}:=\exp (2\pi\sqrt{-1}/d_I) \in \mu_{d_I}$. Namely the variety $\tl{E_I^{\circ}}$ is equipped with a good $\hat{\mu}$-action in the sense of \cite[Section 2.4]{D-L-2}. Following the notations in \cite{D-L-2}, denote by $\M_{\CC}^{\hat{\mu}}$ the ring obtained from the Grothendieck ring $\KK_0^{\hat{\mu}}(\Var_{\CC})$ of varieties over $\CC$ with good $\hat{\mu}$-actions by inverting the Lefschetz motive $\LL\simeq \CC \in \KK_0^{\hat{\mu}}(\Var_{\CC})$. Recall that $\LL \in \KK_0^{\hat{\mu}}(\Var_{\CC})$ is endowed with the trivial action of $\hat{\mu}$.

\begin{definition}
We define the motivic Milnor fiber at infinity $\SS_f^{\infty}$ of the polynomial map $f \colon \CC^n \longrightarrow \CC$ by
\begin{equation}\label{MMF}
\SS_f^{\infty} =\sum_{I \neq \emptyset} (1-\LL)^{\sharp I -1} [\tl{E_I^{\circ}}] \in \M_{\CC}^{\hat{\mu}}.
\end{equation}
\end{definition}

\begin{remark}
By Guibert-Loeser-Merle \cite[Theorem 3.9]{G-L-M}, the motivic Milnor fiber at infinity $\SS_f^{\infty}$ of $f$ does not depend on the compactification $X$ of $\CC^n$. This fact was informed to us by Sch{\"u}rmann (a private communication) and Raibaut \cite{Raibaut}. 
\end{remark}

As in \cite[Section 3.1.2 and 3.1.3]{D-L-2}, we denote by $\HSm$ the abelian category of Hodge structures with a quasi-unipotent endomorphism. Then, to the object $\psi_h(j_!Rf_!\CC_{\CC^n})\in \Dbc(\{\infty\})$ and the semisimple part of the monodromy automorphism acting on it, we can associate an element
\begin{equation}
[H_f^{\infty}] \in \KK_0(\HSm)
\end{equation}
in an obvious way. Similarly, to $\psi_h(Rj_*Rf_*\CC_{\CC^n})\in \Dbc(\{\infty\})$ we associate an element
\begin{equation}
[G_f^{\infty}] \in \KK_0(\HSm).
\end{equation}
According to a deep result \cite[Theorem 13.1]{Sabbah-2} of Sabbah, if $f$ is tame at infinity then the weights of the element $[G_f^{\infty}]$ are defined by the monodromy filtration up to some Tate twists (see also \cite{Saito-1} and \cite{Saito-2}). This implies that for the calculation of the monodromy at infinity $\Phi_{n-1}^{\infty}\colon H^{n-1}(f^{-1}(R);\CC) \simto H^{n-1}(f^{-1}(R);\CC)$ ($R \gg 0$) of $f$ it suffices to calculate $[H_f^{\infty}]\in \KK_0(\HSm)$ which is the dual of $[G_f^{\infty}]$.

To describe the element $[H_f^{\infty}]\in \KK_0(\HSm)$ in terms of $\SS_f^{\infty}\in \M_{\CC}^{\hat{\mu}}$, let
\begin{equation}
\chi_h \colon \M_{\CC}^{\hat{\mu}} \longrightarrow \KK_0(\HSm)
\end{equation}
be the Hodge characteristic morphism defined in \cite{D-L-2} which associates to a variety $Z$ with a good $\mu_d$-action the Hodge structure
\begin{equation}
\chi_h ([Z])=\sum_{j \in \ZZ} (-1)^j [H_c^j(Z;\QQ)] \in \KK_0(\HSm)
\end{equation}
with the actions induced by the one $z \longmapsto \exp (2\pi\sqrt{-1}/d)z$ ($z\in Z$) on $Z$. Then by applying the proof of \cite[Theorem 4.2.1]{D-L-1} to our situation \eqref{eq:4-2} and \eqref{eq:4-5}, we obtain the following result.

\begin{theorem}\label{thm:7-6}
In the Grothendieck group $\KK_0(\HSm)$, we have
\begin{equation}
[H_f^{\infty}]=\chi_h(\SS_f^{\infty}).
\end{equation}
\end{theorem}

On the other hands, the results in \cite{Sabbah} and \cite{Sabbah-2} imply the following symmetry of the weights of the element $[H_f^{\infty}] \in \KK_0(\HSm)$ when $f$ is tame at infinity. See Appendix for the details. Another proof was found later also in \cite{D-S-new}. Recall that if $f$ is tame at infinity we have $H_c^j(f^{-1}(R); \CC )=0$ ($R \gg 0$) for $j \not= n-1, 2n-2$ and $H_c^{2n-2}(f^{-1}(R); \CC ) \simeq [H^0(f^{-1}(R); \CC )]^* \simeq \CC$. For an element $[V] \in \KK_0(\HSm)$, $V \in \HSm $ with a quasi-unipotent endomorphism $\Theta \colon V \simto V$, $p, q \geq 0$ and $\lambda \in \CC$ denote by $e^{p,q}([V])_{\lambda}$ the dimension of the $\lambda$-eigenspace of the morphism $V^{p,q} \simto V^{p,q}$ induced by $\Theta$ on the $(p,q)$-part $V^{p,q}$ of $V$.

\begin{theorem}[Sabbah \cite{Sabbah} and \cite{Sabbah-2}]\label{cor:7-6-2}
Assume that $f$ is tame at infinity. Then
\begin{enumerate}
\item Let $\lambda \in \CC^* \setminus \{1\}$. Then we have $e^{p,q}( [H_f^{\infty}])_{\lambda}=0$ for $(p,q) \notin [0,n-1] \times [0,n-1]$. Moreover for $(p,q) \in [0,n-1] \times [0,n-1]$ we have
\begin{equation}
e^{p,q}( [H_f^{\infty}])_{\lambda}=e^{n-1-q,n-1-p}( [H_f^{\infty}])_{\lambda}.
\end{equation}
\item We have $e^{p,q}( [H_f^{\infty}])_{1}=0$ for $(p,q) \notin \{(n-1, n-1)\} \sqcup ([0,n-2] \times [0,n-2])$ and $e^{n-1,n-1}( [H_f^{\infty}])_{1}=1$. Moreover for $(p,q) \in [0,n-2] \times [0,n-2]$ we have
\begin{equation}
e^{p,q}( [H_f^{\infty}])_{1}=e^{n-2-q,n-2-p}( [H_f^{\infty}])_{1}.
\end{equation}
\end{enumerate}
\end{theorem}

Using our results in Section \ref{sec:7-2}, we can check the above symmetry by explicitly calculating $\chi_h(\SS_f^{\infty})$ for small $n$'s. Since the weights of $[G_f^{\infty}] \in \KK_0(\HSm)$ are defined by the monodromy filtration and $[G_f^{\infty}]$ is the dual of $[H_f^{\infty}]$ up to some Tate twist, we obtain the following result.

\begin{theorem}\label{thm:7-6-2}
Assume that $f$ is tame at infinity. Then
\begin{enumerate}
\item Let $\lambda \in \CC^* \setminus \{1\}$ and $k \geq 1$. Then the number of the Jordan blocks for the eigenvalue $\lambda$ with sizes $\geq k$ in $\Phi_{n-1}^{\infty} \colon H^{n-1}(f^{-1}(R);\CC) \simto H^{n-1}(f^{-1}(R);\CC)$ ($R \gg 0$) is equal to
\begin{equation}
(-1)^{n-1} \sum_{p+q=n-2+k, n-1+k} e^{p,q}( \chi_h(\SS_f^{\infty}))_{\lambda}.
\end{equation}
\item For $k \geq 1$, the number of the Jordan blocks for the eigenvalue $1$ with sizes $\geq k$ in $\Phi_{n-1}^{\infty}$ is equal to
\begin{equation}
(-1)^{n-1} \sum_{p+q=n-2-k, n-1-k} e^{p,q}( \chi_h(\SS_f^{\infty}))_{1}.
\end{equation}
\end{enumerate}
\end{theorem}

By using Newton polyhedrons at infinity, we can rewrite the result of Theorem \ref{thm:7-6} more neatly as follows. Let $f\in \CC[x_1,\ldots,x_n]$ be a convenient polynomial. Assume that $f$ is non-degenerate at infinity. Then $f$ is tame at infinity and it suffices to calculate $\Phi_j^{\infty}$ only for $j=n-1$. From now on, we will freely use the notations in the proof of Theorem \ref{thm:3-5}. For example, $\rho_1, \ldots, \rho_m$ are the $1$-dimensional cones in the smooth fan $\Sigma$ such that $\rho_i \not\subset \RR_+^n$. We call these cones the rays at infinity. Each ray $\rho_i$ at infinity corresponds to the toric divisor $D_i$ in $X_{\Sigma}$ and the divisor $D:=D_1\cup \cdots \cup D_m=X_{\Sigma} \setminus \CC^n$ in $X_{\Sigma}$ is normal crossing. We denote by $a_i>0$ the order of the poles of $f$ along $D_i$. By eliminating the points of indeterminacy of the meromorphic extension of $f$ to $X_{\Sigma}$ we constructed the commutative diagram:
\begin{equation}
\xymatrix{
\CC^n \ar@{^{(}->}[r]^{\iota} \ar[d]_f & \tl{X_{\Sigma}} \ar[d]^g\\\CC \ar@{^{(}->}[r]^j & \PP^1.}
\end{equation}
Recall that in the construction of $\tl{X_{\Sigma}}$ we first construct a tower of blow-ups over $D_1 \cap \overline{f^{-1}(0)}$ and next apply the same operation to the remaining divisors $D_2, \ldots, D_m$ (in this order). See the proof of Theorem \ref{thm:3-5} and Lemma \ref{geometry} for the details. Take a local coordinate $h$ of $\PP^1$ in a neighborhood of $\infty \in \PP^1$ such that $\infty=\{h=0\}$ and set $\tl{g}=h\circ g$, $Y=\tl{g}^{-1}(0)=g^{-1}(\infty) \subset \tl{X_{\Sigma}}$ and $\Omega=\Int(\iota(\CC^n) \sqcup Y)$. For simplicity, let us set $\tl{g}=\frac{1}{f}$. Then the divisor $U=Y \cap \Omega$ in $\Omega$ contains not only the proper transforms $D_1^{\prime}, \ldots, D_m^{\prime}$ of $D_1, \ldots, D_m$ in $\tl{X_{\Sigma}}$ but also the exceptional divisors of the blow-up: $\tl{X_{\Sigma}} \longrightarrow X_{\Sigma}$. From now on, we will show that these exceptional divisors are not necessary to compute the monodromy at infinity of $f \colon \CC^n \longrightarrow \CC$ by Theorem \ref{thm:7-6}. For each non-empty subset $I \subset \{1,2,\ldots, m\}$, set $D_I= \bigcap_{i \in I} D_i$,
\begin{equation}
D_I^{\circ}=D_I \setminus \left\{ \( \bigcup_{i \notin I}D_i\) \cup \overline{f^{-1}(0)}\right\} \subset X_{\Sigma}
\end{equation}
and $d_I =\gcd(a_i)_{i \in I} >0$. Then the function $\tl{g}=\frac{1}{f}$ is regular on $D_I^{\circ}$ and we can decompose it as $\frac{1}{f}=\tl{g_{1}}(\tl{g_{2}})^{d_I}$ globally on a Zariski open neighborhood $W$ of $D_I^{\circ}$ in $X_{\Sigma}$, where $\tl{g_1}$ is a unit on $W$ and $\tl{g_2} \colon W \longrightarrow \CC$ is regular. Therefore we can construct an unramified Galois covering $\tl{D_I^{\circ}}$ of $D_I^{\circ}$ with a natural $\mu_{d_I}$-action as in \eqref{eq:6-26}. Let $[\tl{D_I^{\circ}}]$ be the element of the ring $\M_{\CC}^{\hat{\mu}}$ which corresponds to $\tl{D_I^{\circ}}$.

\begin{theorem}\label{thm:7-7}
Assume that $f$ is convenient and non-degenerate at infinity. Then we have the equality
\begin{equation}
\chi_h\(\SS_f^{\infty}\)= \dsum_{I \neq \emptyset}\chi_h\( (1-\LL)^{\sharp I -1} [\tl{D_I^{\circ}}]\)
\end{equation}
in the Grothendieck group $\KK_0(\HSm)$.
\end{theorem}

\begin{proof}
First, we prove the assertion for $n=2$. In this case, we number the rays at infinity in $\Sigma$ in the clockwise direction as in the figure below.

\vspace{1mm}
\begin{minipage}[t]{0.4\textwidth}
\begin{center}
\includegraphics[scale=0.8]{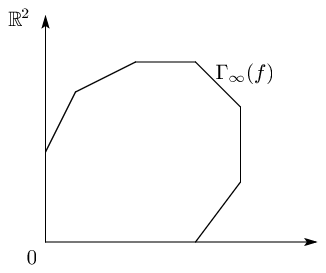}

Figure 4
\end{center}
\end{minipage}
\hspace{20mm}\begin{minipage}[t]{0.4\textwidth}

\begin{center}
\includegraphics[scale=0.75]{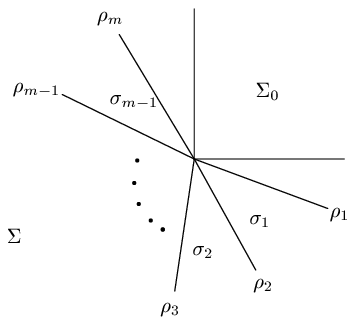}

Figure 5
\end{center}
\end{minipage}

\vspace{2mm}
Let $\sigma_i=\RR_{+} \rho_i +\RR_{+} \rho_{i+1}$ ($1 \leq i \leq m-1$) be the $2$-dimensional cone in $\Sigma$ between $\rho_i$ and $\rho_{i+1}$. Then the cone $\sigma_i$ corresponds to an affine open subset $\CC^2(\sigma_i)\simeq \CC^2_{\xi, \eta}$ of $X_{\Sigma}$ on which the meromorphic extension of $\frac{1}{f}$ to $X_{\Sigma}$ has the form
\begin{equation}
\(\dfrac{1}{f}\) (\xi, \eta)=\dfrac{\xi^{a_i}\eta^{a_{i+1}}}{f_{\sigma_i}(\xi,\eta)},
\end{equation}
where $f_{\sigma_i}$ is a polynomial of $\xi$ and $\eta$. In this situation, we have $D_i \cap \CC^2(\sigma_i)=\{\xi=0\}$ and $D_{i+1}\cap \CC^2(\sigma_i)=\{\eta=0\}$. Moreover by the non-degeneracy at infinity of $f(x,y) \in \CC[x,y]$ we have $f_{\sigma_i}(0,0)\neq 0$ and the algebraic curve $f_{\sigma_i}^{-1}(0)=\{(\xi,\eta) \ |\ f_{\sigma_i}(\xi, \eta)=0\}$ intersects $D_i \cap \CC^2(\sigma_i)$ and $D_{i+1} \cap \CC^2(\sigma_i)$ transversally. 

\vspace{1mm}
\begin{center}
\includegraphics[scale=.8]{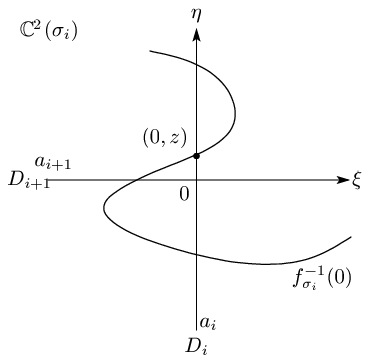}

Figure 6
\end{center}

\vspace{1mm}
Let $(0,z)$, $z\neq 0$ be a point of $D_i \cap f_{\sigma_i}^{-1}(0)$. In constructing the variety $\tl{X_{\Sigma}}$, we constructed a tower of blow-ups over this point $(0,z)$ as in the figure below.

\vspace{1mm}
\begin{center}
\includegraphics[scale=.9]{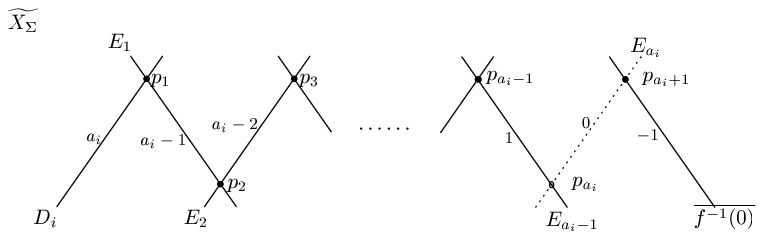}

Figure 7
\end{center}

\vspace{1mm}
Here we have $E_j \simeq \PP^1$ and the function $\frac{1}{f} =\tl{g}$ has the zero of order $a_i-j$ along the exceptional divisor $E_j$. The open set $\Omega$ is the complement of $E_{a_i} \simeq \PP^1$ in this figure. For $1 \leq j \leq a_i -1$, set $E_j^{\circ}:=E_j \setminus \{p_j, p_{j+1}\}$ and let $\tl{E_j^{\circ}}$ be the unramified Galois covering of $E_j^{\circ}$ with a $\mu_{a_i-j}$-action (in the construction of $\SS_f^{\infty}$). The motivic Milnor fiber at infinity $\SS_f^{\infty}$ also contains $(1-\LL) \cdot [p_j]\in \M_{\CC}^{\hat{\mu}}$ with the trivial $\hat{\mu}$-action for $1\leq j \leq a_i-1$. 

\begin{lemma}\label{lem:7-7-2}
For $1 \leq j \leq a_i-1$, we have
\begin{equation}
\chi_h\((1-\LL)\cdot [p_j]+[\tl{E_j^{\circ}}]\)=0
\end{equation}
in the Grothendieck group $\KK_0(\HSm)$.
\end{lemma}

\begin{proof}
First set $X=\CC^2(\sigma_i) \setminus \{\eta=0\}$ and $Y=\{(0,z)\}\subset X$ and consider the regular functions
\begin{equation}
f_1(\xi,\eta)=\xi, \hspace{10mm}f_2(\xi,\eta)=\dfrac{f_{\sigma_i}(\xi,\eta)}{\eta^{a_{i+1}}}
\end{equation}
on $X$. Then (in a neighborhood of $Y$) the blow-up $\tl{X_Y}$ of $X$ along $Y$ is isomorphic to the closure of the image of the morphism
\begin{equation}
X\setminus Y \longrightarrow X \times \PP^1
\end{equation}
given by
\begin{equation}
(\xi,\eta) \longmapsto (\xi,\eta, (f_1(\xi,\eta):f_2(\xi,\eta))).
\end{equation}
Let $\pi \colon \tl{X_Y} \longtwoheadrightarrow X$ be the natural morphism and define an open subset $W$ of $\tl{X_Y}$ by
\begin{equation}
W=\{(\xi,\eta, (1:\alpha)) \in \tl{X_Y} \ |\ \alpha \in \CC\}.
\end{equation}
Then considering $\alpha$ as a regular function on $W$, on $W \subset \tl{X_Y}$ we have
\begin{equation}
\(\frac{1}{f}\) \circ \pi =\dfrac{(f_1\circ \pi)^{a_i}}{(f_2\circ \pi)}=\dfrac{(f_1\circ \pi)^{a_i}}{\alpha(f_1\circ\pi)}=\dfrac{(f_1\circ \pi)^{a_i-1}}{\alpha}.
\end{equation}
Moreover in $W \subset \tl{X_Y}$ the exceptional divisor $E=\pi^{-1}(Y)$ ($\simeq E_1$) is defined by $E=\{f_1\circ \pi =0\}$. Therefore the unramified Galois covering $\tl{E_1^{\circ}}$ of $E_1^{\circ}\simeq \{\alpha \in \CC\ |\ \alpha \neq 0\}\simeq \CC^*$ is $\{(t ,\alpha) \in (\CC^*)^2\ |\ t^{a_i-1}\alpha^{-1}=1\}$, which is isomorphic to $\CC^*$ with an automorphism homotopic to the identity. Hence its Hodge characteristic $\chi_h \( [\tl{E_1^{\circ}}] \) \in \KK_0(\HSm)$ is isomorphic to $\chi_h\( \LL -1 \)$. We thus obtain the equality
\begin{equation}
\chi_h\((1-\LL)\cdot [p_1] +[\tl{E_1^{\circ}}]\)=0
\end{equation}
in $\KK_0(\HSm)$. By repeating this argument, we can similarly prove the remaining assertions. \qed
\end{proof}

Recall that the motivic Milnor fiber at infinity $\SS_f^{\infty}$ is a sum of the unramified Galois coverings of some Zariski locally closed subvarieties of $\Omega$. Then Lemma \ref{lem:7-7-2} above implies that the Hodge characteristics of the base changes of $\SS_f^{\infty}$ to the exceptional divisors of $\Omega \longrightarrow X_{\Sigma}$ are zero in $\KK_0(\HSm)$. In other words, for the calculation of $\chi_h\(\SS_f^{\infty}\)\in \KK_0(\HSm)$, we can forget the parts of $\chi_h\(\SS_f^{\infty}\) \in \KK_0(\HSm)$ coming from the exceptional divisors of $\Omega \longrightarrow X_{\Sigma}$. So the theorem was proved in the case $n=2$.

\medskip
From now on, we shall prove Theorem \ref{thm:7-7} in the case $n>2$. Let $\pi_{\Omega} \colon \Omega \longtwoheadrightarrow X_{\Sigma}$ be the restriction of the morphism $\pi \colon \tl{X_{\Sigma}} \longtwoheadrightarrow X_{\Sigma}$ to $\Omega$. For each non-empty subset $I \subset \{1,2,\ldots,m\}$ set $D_I^*:=D_I\setminus \( \bigcup_{i \notin I}D_i\)$. Then, to prove the theorem, it suffices to show that for any non-empty subset $I \subset \{1,2,\ldots,m\}$ the Hodge characteristic of the base change of $\SS_f^{\infty}$ to $\pi_{\Omega}^{-1}(D_I^* \cap \overline{f^{-1}(0)}) \subset \Omega$ is zero in $\KK_0(\HSm)$. First, let us consider the case where $I=\{ i\}$ for some $1 \leq i \leq m$. Then by the construction of $\tl{X_{\Sigma}}$ and $\Omega$, the morphism $\pi_{\Omega}$ induces a fiber bundle
\begin{equation}
\pi_{\Omega}^{-1}\( D_{\{i\}}^* \cap \overline{f^{-1}(0)}\) \longtwoheadrightarrow D_{\{i\}}^* \cap \overline{f^{-1}(0)}
\end{equation}
over $D_{\{i\}}^* \cap \overline{f^{-1}(0)}$ whose fiber is isomorphic to the (locally closed) curve $(E_1 \cup \cdots \cup E_{a_i-1})\setminus \{p_{a_i}\}$ as in Figure 7. By the proof of Lemma \ref{lem:7-7-2}, this fiber bundle is locally trivial with respect to the Zariski topology of $D_{\{i\}}^* \cap \overline{f^{-1}(0)}$, and the Hodge characteristic of the base change of $\SS_f^{\infty}$ to $\pi_{\Omega}^{-1}(D_{\{i\}}^* \cap \overline{f^{-1}(0)}) \subset \Omega$ is zero in $\KK_0(\HSm)$. Next, consider a general non-empty subset $I=\{ i_1<\cdots <i_k\}$ of $\{1,2,\ldots, m\}$. Then we have the following lemma.

\begin{lemma}\label{geometry}
In the situation as above, the restriction of $\pi_{\Omega}$ to $\pi_{\Omega}^{-1}(D_I^* \cap \overline{f^{-1}(0)})$:
\begin{equation}
\pi_{\Omega}^{-1}(D_I^* \cap \overline{f^{-1}(0)}) \longrightarrow D_I^* \cap \overline{f^{-1}(0)}
\end{equation}
is a Zariski locally trivial bundle over $D_I^* \cap \overline{f^{-1}(0)}$ whose fiber is isomorphic to the (locally closed) curve $(E_1 \cup \cdots\cup E_{a_{i_1-1}}) \setminus \{p_{a_{i_1}}\}$ in Figure 7 for $i=i_1$.
\end{lemma}

\begin{proof}
First, let us consider the case where $\sharp I=2$. Without loss of generality, we may assume that $I=\{ 1,2 \}$ as in the figure:

\begin{center}
\includegraphics[scale=.6]{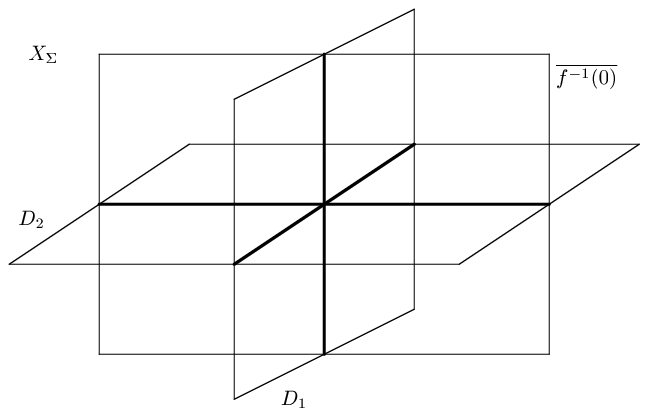}

Figure 8
\end{center}

\vspace{1mm}
Although we consider the problem in a neighborhood of $D_I^* \subset X_{\Sigma}$, the total space of Figure 8 is denoted simply by $X_{\Sigma}$. In the construction of $\tl{X_{\Sigma}}$, we first construct a tower of blow-ups over $D_1 \cap \overline{f^{-1}(0)}$ as

\vspace{1mm}
\begin{center}
\includegraphics[scale=.75]{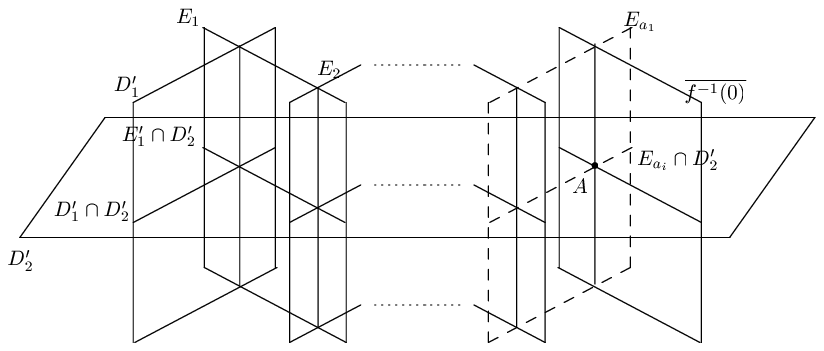}

Figure 9
\end{center}

\vspace{1mm}
Here $D_1^{\prime}$, $D_2^{\prime}$ are the proper transforms of $D_1$, $D_2$ respectively. In Figure 9, the meromorphic function $\tl{g}$ still has points of indeterminacy on $D_2^{\prime} \cap \overline{f^{-1}(0)}$. Then we construct a tower of blow-ups over $D_2^{\prime} \cap \overline{f^{-1}(0)}$ as

\vspace{1mm}
\begin{center}
\includegraphics[scale=.75]{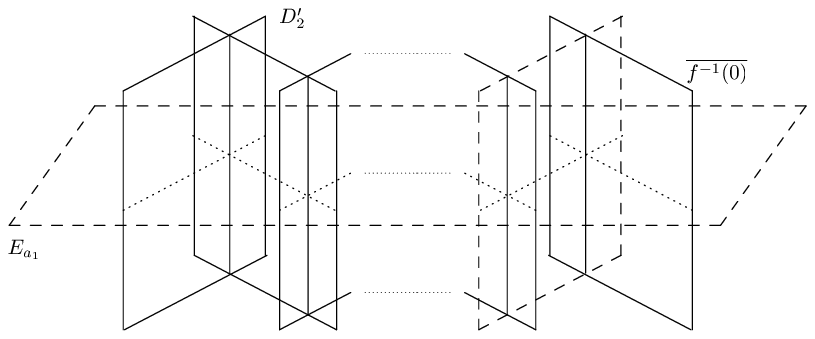}

Figure 10
\end{center}

\vspace{1mm}
Since we are considering the problem only in a neighborhood of $D_I^*$ and already finished the necessary blow-ups over $D_I^*$, the total space of Figure 10 is denoted simply by $\tl{X_{\Sigma}}$. In Figure 10, the open set $\Omega \subset \tl{X_{\Sigma}}$ is the complement of the union of two dotted divisors ($E_{a_1}^{\prime}$ is one of them). Moreover we see that the inverse image of the set $A$ (in Figure 9) in $\tl{X_{\Sigma}}$ is contained in $\tl{X_{\Sigma}}\setminus \Omega$. This implies that $\pi_{\Omega}^{-1}(D_I^* \cap \overline{f^{-1}(0)})$ is the (locally closed) variety of the form:

\vspace{1mm}
\begin{center}
\includegraphics[scale=.75]{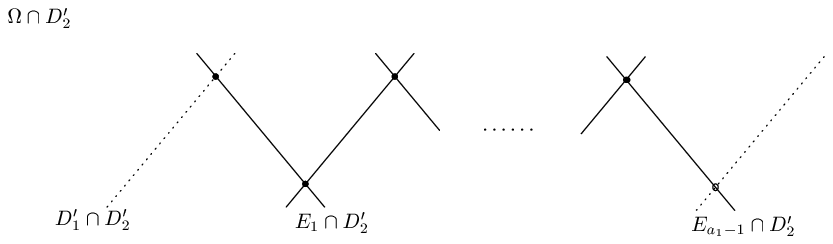}

Figure 11
\end{center}

\vspace{1mm}
This completes the proof for the case $\sharp I=2$. The general case can be proved similarly. \qed
\end{proof}

By the proof of Lemma \ref{geometry} we see that $\pi_{\Omega}^{-1}(D_I^* \cap \overline{f^{-1}(0)})$ has a geometric structure as the figure below in $\Omega \setminus \pi_{\Omega}^{-1}\( \bigcup_{i \notin I}D_i\)$:

\vspace{1mm}
\begin{center}
\includegraphics[scale=.75]{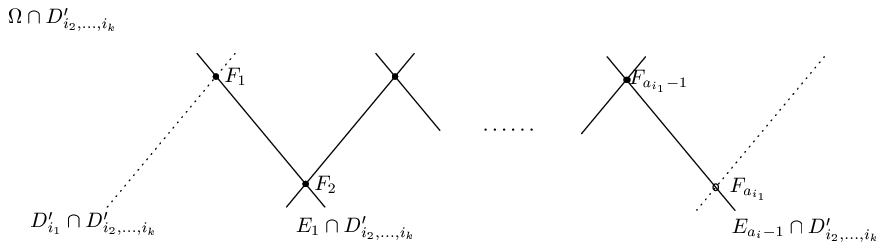}

Figure 12
\end{center}

\vspace{1mm}
Here $E_1,\ldots, E_{a_{i_1}}$ are the exceptional divisors in $\tl{X_{\Sigma}}$ constructed when we made a tower of blow-ups over $D_{i_1} \cap \overline{f^{-1}(0)}$ (We used essentially the condition $i_1=\min \{i_1,\ldots, i_k\}$. To simplify the notations, we denote $E_j \setminus \pi_{\Omega}^{-1}\( \bigcup_{i \notin I}D_i\)$ simply by $E_j$ etc.). Moreover we set $D_{i_2, \ldots, i_k}^{\prime} := D_{i_2}^{\prime} \cap \cdots \cap D_{i_k}^{\prime}$ and $F_j:=E_{j-1} \cap E_j \cap D_{i_2, \ldots, i_k}^{\prime}$. Note that $E_j \cap D_{i_2, \ldots, i_k}^{\prime}$ is a $\PP^1$-bundle over $D_I^* \cap \overline{f^{-1}(0)}$. Let us set $(E_j \cap D_{i_2, \ldots, i_k}^{\prime})^{\circ}:= (E_j \cap D_{i_2,\ldots, i_k}^{\prime}) \setminus (F_j \sqcup F_{j+1})$. Then for each point $p$ of $(E_j \cap D_{i_2,\ldots, i_k}^{\prime})^{\circ}$ there exists a Zariski open neighborhood $W$ of $p$ in $\Omega$ and a local coordinate system $\xi_1,\xi_2,\ldots, \xi_n$ on it such that
\begin{gather}
D_{i_2}^{\prime}=\{\xi_2=0\}, \ldots, D_{i_k}^{\prime}=\{\xi_k=0\},\\
E_j \cap D_{i_2,\ldots, i_k}^{\prime}=\{\xi_1=\xi_2=\cdots =\xi_k=0\}
\end{gather}
and the function $\frac{1}{f}=\tl{g}$ can be written in the form
\begin{equation}
\tl{g}(\xi_1,\ldots, \xi_n)=\xi_1^{a_{i_1}-j}\xi_2^{a_{i_2}}\cdots \xi_k^{a_{i_k}} \times (\text{ a unit on $W$})
\end{equation}
on $W$. Set $d_{I,j}=\gcd(a_{i_1}-j, a_{i_2}, \ldots, a_{i_k})>0$ ($1 \leq j \leq a_{i_1}-1$). Then the base change of $\SS_f^{\infty}$ to $(E_j \cap D_{i_2,\ldots, i_k}^{\prime})^{\circ} \subset \Omega$ is an unramified Galois covering $\tl{(E_j \cap D_{i_2,\ldots, i_k}^{\prime})^{\circ}}$ of $(E_j \cap D_{i_2,\ldots, i_k}^{\prime})^{\circ}$ with a natural $\mu_{d_{I,j}}$-action. Moreover by the proof of Lemma \ref{lem:7-7-2}, we observe that $\tl{(E_j \cap D_{i_2,\ldots, i_k}^{\prime})^{\circ}}$ is a (Zariski) locally trivial family over $D_I^* \cap \overline{f^{-1}(0)}$. By using this fact (and an analogue of \cite[Proposition 1.6]{D-K}), in the same way as the final part of the proof of Lemma \ref{lem:7-7-2} we obtain the equality 
\begin{equation}
\chi_h\((1-\LL)\cdot [F_j] +[\tl{(E_j \cap D_{i_2,\ldots, i_k}^{\prime})^{\circ}}]\)=0
\end{equation}
in $\KK_0(\HSm)$ for $1 \leq j \leq a_{i_1}-1$, where $[F_j] \in \M_{\CC}^{\hat{\mu}}$ is endowed with the trivial action of $\hat{\mu}$. This implies that the Hodge characteristic of the base change of $\SS_f^{\infty}$ to $\pi_{\Omega}^{-1}(D_I^* \cap \overline{f^{-1}(0)}) \subset \Omega$ is zero in $\KK_0(\HSm)$. In other words, the contribution to $\chi_h\(\SS_f^{\infty}\) \in \KK_0(\HSm)$ from the exceptional divisors of $\pi_{\Omega} \colon \Omega \longrightarrow X_{\Sigma}$ is zero. This completes the proof of Theorem \ref{thm:7-7}. \qed
\end{proof}

\begin{remark}\label{EQH} 
It seems that the equality $\SS_f^{\infty}= \sum_{I \neq \emptyset} (1-\LL)^{\sharp I -1} [\tl{D_I^{\circ}}]$ does not hold in $\M_{\CC}^{\hat{\mu}}$. Indeed, we used a homotopy in the last part of the proof of Lemma \ref{lem:7-7-2} (and that of Theorem \ref{thm:7-7}). 
\end{remark}

\section{Combinatorial descriptions of monodromies at infinity}\label{sec:7-2}

In this section, by rewriting Theorem \ref{thm:7-7} in terms of the Newton polyhedron at infinity $\Gamma_{\infty}(f)$ of $f$ we prove some combinatorial formulas for the Jordan normal form of its monodromy at infinity $\Phi_{n-1}^{\infty}$. We inherit the situation and the notations in the last half of Section \ref{sec:7}. Namely we assume that $f$ is convenient and non-degenerate at infinity. Recall that $\rho_1, \rho_2, \ldots, \rho_m$ are $1$-dimensional cones in the smooth fan $\Sigma$ such that $\rho_i \not\subset \RR^n_+$. 

\begin{definition}
We say that $\gamma \prec \Gamma_{\infty}(f)$ is a face at infinity of $\Gamma_{\infty}(f)$ if $0 \notin \gamma$.
\end{definition}

For a cone $\sigma \in \Sigma$ whose supporting face $\gamma (\sigma) \prec \Gamma_{\infty}(f)$ is at infinity ($\Longleftrightarrow$ $\sigma \not\subset \RR^n_+$) we set $I_{\sigma}= \{ 1 \leq  i \leq m \ | \ \rho_i \prec \sigma \}$, $T_{\sigma}^{\circ}=T_{\sigma} \setminus \overline{f^{-1}(0)}$ and
\begin{equation}
\tl{T_{\sigma}^{\circ}} = \tl{D_{I_{\sigma}}^{\circ}} \cap (\CC_t^* \times T_{\sigma}) \subset \CC^*_t \times T_{\sigma}.
\end{equation}
Then $\tl{T_{\sigma}^{\circ}}$ is a hypersurface in the algebraic torus $\CC_t^* \times T_{\sigma} \simeq (\CC^*)^{n- \d \sigma +1}$ and a finite covering of $T_{\sigma}^{\circ}$. Moreover for any non-empty subset $I \subset \{ 1,2, \ldots, m \}$ we have the decomposition:
\begin{equation}
\tl{D_I^{\circ}}= \bigsqcup_{\sigma \colon  I_{\sigma}=I} \tl{T_{\sigma}^{\circ}}.
\end{equation}
Therefore, for the calculation of $\chi_h([\tl{D_I^{\circ}}])\in \KK_0(\HSm)$ by the results in Section \ref{sec:2} we have to show that the hypersurfaces $\tl{T_{\sigma}^{\circ}} \subset \CC_t^* \times T_{\sigma} \simeq (\CC^*)^{n- \d \sigma +1}$ are non-degenerate. Indeed, for such a cone $\sigma \in \Sigma$ let $\sigma_0 \in \Sigma$ be an $n$-dimensional cone such that $\sigma \prec \sigma_0$ and $\{ w_1, w_2, \ldots, w_n \} \subset \ZZ^n$ the set of the primitive vectors on the edges of $\sigma_0$. Set $\d \sigma =k >0$. Then we may assume that $w_1, \ldots, w_k$ generate $\sigma$ so that in the affine open subset $\CC^n (\sigma_0) \simeq \CC^n_y$ of $X_{\Sigma}$ associated to $\sigma_0$ we have
\begin{equation}
T_{\sigma}=\{(y_1,\ldots,y_n) \in \CC^n(\sigma_0)\ |\ y_1=\cdots=y_k=0, \ y_{k+1},\ldots, y_n \neq 0 \}.
\end{equation}
On $\CC^n (\sigma_0) \simeq \CC^n_y$ the function $\tl{g}= \frac{1}{f}$ has the form:
\begin{equation}\label{eq:new5-4}
\tl{g}(y)= y_1^{c_1} \cdots y_n^{c_n} \times \frac{1}{f_{\sigma_0}(y)},
\end{equation}
where we set
\begin{equation}
c_j=-\min_{v\in \Gamma_{\infty}(f)} \langle w_j,v \rangle \geq 0 \qquad (j=1,2, \ldots, n)
\end{equation}
and $f_{\sigma_0}(y)$ is a polynomial on $\CC^n (\sigma_0)$. Set $d= \gcd(c_1,\ldots,c_k):=\gcd (\{c_j\ |\ 1\leq j \leq k, \ c_j\neq 0\})>0$. Then in $\CC_t^* \times T_{\sigma} \simeq (\CC^*)^{n- \d \sigma +1}_{t,y_{k+1}, \ldots, y_n}$ we have
\begin{equation}\label{DEQ} 
\tl{T_{\sigma}^{\circ}} = \{ (t, y_{k+1}, \ldots, y_n) \ | \ t^{-d}y_{k+1}^{-c_{k+1}} \cdots y_n^{-c_n} \times (f_{\sigma_0}|_{T_{\sigma}})(y_{k+1}, \ldots, y_n)=1 \}
\end{equation}
and the action $\Psi_{\sigma}$ of the generator of the cyclic group $\mu_d$ on it is given by the multiplication of $( \zeta_d, 1, \ldots, 1) \in \CC_t^* \times T_{\sigma}$. To show that the hypersurfaces $\tl{T_{\sigma}^{\circ}} \subset \CC_t^* \times T_{\sigma}$ are non-degenerate, we use the following elementary lemma.

\begin{lemma}\label{lem:key}
Let $g_0$ be a Laurent polynomial on $(\CC^*)^n$ such that the hypersurface $Z^* =\{ x\in (\CC^*)^n \ |\ g_0(x)=0\}$ is non-degenerate and $x^v$ be a monomial. Then the set of complex numbers $\lambda \in \CC$ such that the hypersurface $Z_{\lambda}^*=\{x\in (\CC^*)^n \ | \ g_0(x)-\lambda x^v=0\}$ is non-degenerate is open dense in $\CC$.
\end{lemma}

\begin{proof}
It is easy to see that for $c \in \CC$ and $x \in (\CC^*)^n$ the following two conditions are equivalent.
\begin{gather}
\text{$\dfrac{\partial}{\partial x_j} \{ g_0(x)-c x^v \} =0$ \ ($1 \leq j \leq n$), \quad $g_0(x)-c x^v=0$.}\\
\text{$\dfrac{\partial}{\partial x_j} \left\{ \dfrac{g_0(x)}{x^v} \right\} =0$ \ ($1 \leq j \leq n$), \quad $\dfrac{g_0(x)}{x^v}=c$.}
\end{gather}
Then by applying the Bertini-Sard theorem to the map $\dfrac{g_0}{x^v}\colon (\CC^*)^n \longrightarrow \CC$ we find that the hypersurface $\{ x \in (\CC^*)^n \ | \ g_0(x)-c x^v =0 \}$ in $(\CC^*)^n$ is smooth and reduced for generic $c \in \CC$. Note that for generic $c \in \CC$ the Newton polytope $NP(g_0-c x^v)$ of $g_0-c x^v$ is the convex hull of $NP(g_0) \cup \{ v \}$ and has only finitely many faces. For its face $\Gamma \prec NP(g_0-c x^v)$ let $(g_0-c x^v)^{\Gamma}$ be the $\Gamma$-part of the Laurent polynomial $g_0-c x^v$. Then the assertion follows by applying the above argument for $g_0-c x^v$ to $(g_0-c x^v)^{\Gamma}$. \qed
\end{proof}

\begin{proposition}\label{PND}
In the situation as above, the hypersurfaces $\tl{T_{\sigma}^{\circ}} \subset \CC_t^* \times T_{\sigma}$ are non-degenerate.
\end{proposition}

\begin{proof}
By the non-degeneracy at infinity of $f$ and Lemma \ref{lem:key} there exists $\lambda \in \CC^*$ such that the hypersurface
\begin{equation}
\{ (t, y_{k+1}, \ldots, y_n) \ | \ t^{-d}y_{k+1}^{-c_{k+1}} \cdots y_n^{-c_n} \times (f_{\sigma_0}|_{T_{\sigma}})(y_{k+1}, \ldots, y_n)= \lambda \}
\end{equation}
in $\CC_t^* \times T_{\sigma}$ is non-degenerate. Since it is isomorphic to $\tl{T_{\sigma}^{\circ}}$ by the multiplication by $( \lambda^{\prime}, 1, \ldots, 1) \in \CC_t^* \times T_{\sigma}$ for $\lambda^{\prime}$ satisfying $(\lambda^{\prime})^d= \lambda$, $\tl{T_{\sigma}^{\circ}}$ is also non-degenerate. \qed
\end{proof}

Note that if $\d \gamma (\sigma )=n- \d \sigma$ the integer $d>0$ above is equal to the lattice distance of $\gamma (\sigma )$ from the origin $0 \in \RR^n$. Moreover in this case, by \eqref{DEQ} the Newton polytope of the defining equation of $\tl{T_{\sigma}^{\circ}}$ in $\CC_t^* \times T_{\sigma} \simeq (\CC^*)^{n- \d \sigma +1}$ is the convex hull of $\{ 0 \} \sqcup \gamma (\sigma )$. For each face at infinity $\gamma \prec \Gamma_{\infty}(f)$ of $\Gamma_{\infty}(f)$, let $d_{\gamma}>0$ be the lattice distance of $\gamma$ from the origin $0 \in \RR^n$ and $\Delta_{\gamma}$ the convex hull of $\{0\} \sqcup \gamma$ in $\RR^n$. Let $\LL(\Delta_{\gamma})$ be the $(\dim \gamma +1)$-dimensional linear subspace of $\RR^n$ spanned by $\Delta_{\gamma}$ and consider the lattice $M_{\gamma}=\ZZ^n \cap \LL(\Delta_{\gamma}) \simeq \ZZ^{\dim \gamma+1}$ in it. Then we set $T_{\Delta_{\gamma}}:=\Spec (\CC[M_{\gamma}]) \simeq (\CC^*)^{\dim \gamma +1}$. Moreover let $\LL(\gamma)$ be the smallest affine linear subspace of $\RR^n$ containing $\gamma$ and for $v \in M_{\gamma}$ define their lattice heights $\height (v, \gamma) \in \ZZ$ from $\LL(\gamma)$ in $\LL(\Delta_{\gamma})$ so that we have $\height (0, \gamma)=d_{\gamma}>0$. Then to the group homomorphism $M_{\gamma} \longrightarrow \CC^*$ defined by $v \longmapsto \zeta_{d_{\gamma}}^{\height (v, \gamma)}$ we can naturally associate an element $\tau_{\gamma} \in T_{\Delta_{\gamma}}$. We define a Laurent polynomial $g_{\gamma}=\sum_{v \in M_{\gamma}}b_v x^v$ on $T_{\Delta_{\gamma}}$ by
\begin{equation}
b_v=\begin{cases}
a_v & (v \in \gamma),\\
-1 & (v=0),\\
\ 0 & (\text{otherwise}),
\end{cases}
\end{equation}
where $f=\sum_{v \in \ZZ^n_+} a_v x^v$. Then we have $NP(g_{\gamma}) =\Delta_{\gamma}$, $\supp g_{\gamma} \subset \{ 0\} \sqcup \gamma$ and the hypersurface $Z_{\Delta_{\gamma}}^*=\{ x \in T_{\Delta_{\gamma}}\ |\ g_{\gamma}(x)=0\}$ is non-degenerate by the proof of Proposition \ref{PND}. Since $Z_{\Delta_{\gamma}}^* \subset T_{\Delta_{\gamma}}$ is invariant by the multiplication $l_{\tau_{\gamma}} \colon  T_{\Delta_{\gamma}} \simto T_{\Delta_{\gamma}}$ by $\tau_{\gamma}$, $Z_{\Delta_{\gamma}}^*$ admits an action of $\mu_{d_{\gamma}}$. We thus obtain an element $[Z_{\Delta_{\gamma}}^*]$ of $\M_{\CC}^{\hat{\mu}}$. By the construction of $[Z_{\Delta_{\gamma}}^*]$ the following lemma is obvious. 

\begin{lemma}\label{MMFI}
Let $\gamma \prec \Gamma_{\infty}(f)$ be a face at infinity of $\Gamma_{\infty}(f)$ and $\sigma \in \Sigma$ a cone whose supporting face $\gamma (\sigma )$ in $\Gamma_{\infty}(f)$ is $\gamma$. Assume that $\d \gamma =n- \d \sigma$. Then in the Grothendieck ring $\M_{\CC}^{\hat{\mu}}$ we have the equality 
\begin{equation}
[ \tl{T_{\sigma}^{\circ}} ]=[Z_{\Delta_{\gamma}}^*].
\end{equation}
\end{lemma}

To rewrite Theorem \ref{thm:7-7} in terms of $\Gamma_{\infty}(f)$ we need the following result.

\begin{proposition}\label{prp:new5-6}
Let $\gamma\prec \Gamma_{\infty}(f)$ be a face at infinity of $\Gamma_{\infty}(f)$ and $\sigma_1, \sigma_2 \in \Sigma$ cones whose supporting faces $\gamma(\sigma_1)$ and $\gamma(\sigma_2)$ in $\Gamma_{\infty}(f)$ are the same and equal to $\gamma$. Then in the Grothendieck group $\KK_0(\HSm)$ we have
\begin{equation}\label{eq:new5-6}
\chi_h( \(\LL-1\)^{\d \sigma_1 -1} \cdot [\tl{T_{\sigma_1}^{\circ}}])= \chi_h( \(\LL-1\)^{\d \sigma_2 -1}\cdot [ \tl{T_{\sigma_2}^{\circ}}]).
\end{equation}
\end{proposition}

\begin{proof}
Without loss of generality we may assume that $\d \sigma_1 \leq \d \sigma_2$. Set $\d \sigma_i =k_i$ ($i=1, 2$). 

\noindent ({\bf Step 1}): First we prove \eqref{eq:new5-6} in the case $\sigma_1\prec \sigma_2$. Let $\sigma_0 \in \Sigma$ be an $n$-dimensional cone such that $\sigma_1 \prec \sigma_2 \prec \sigma_0$ and $\{ w_1, w_2, \ldots, w_n \} \subset \ZZ^n$ the set of the primitive vectors on the edges of $\sigma_0$. We may assume that $w_1, w_2, \ldots, w_{k_i}$ generate $\sigma_i$ for $i=1, 2$. Then in the affine open subset $\CC^n(\sigma_0) \simeq \CC^n_y$ of $X_{\Sigma}$ associated to $\sigma_0$ we have 
\begin{equation}
T_{\sigma_i}=\{ y \in \CC^n(\sigma_0) \ | \ y_1=\cdots =y_{k_i}=0, \ y_{k_i+1},\ldots, y_n \neq 0 \} 
\end{equation}
for $i=1, 2$. Moreover on $\CC^n(\sigma_0) \simeq \CC^n_y$ the function $\tl{g}=\frac{1}{f}$ has the form: 
\begin{equation} 
\tl{g}(y)= y_1^{c_1} \cdots y_n^{c_n} \times g_{\sigma_0}(y),
\end{equation}
where we set
\begin{equation}
c_j=-\min_{v\in \Gamma_{\infty}(f)} \langle w_j,v \rangle \geq 0 \qquad (j=1,2, \ldots, n)
\end{equation}
and $g_{\sigma_0}(y)$ is a meromorphic function on $\CC^n (\sigma_0)$. By the assumption $\gamma(\sigma_1)=\gamma(\sigma_2)=\gamma$, the restriction $g_{\sigma_0}|_{T_{\sigma_1}}$ of $g_{\sigma_0}$ to the larger torus $T_{\sigma_1}$ depends only on the variables $y_{k_2+1},\ldots, y_n$. Set $d_i=\gcd(c_1,\ldots,c_{k_i}):=\gcd(\{c_j\ |\ 1\leq j\leq k_i,\ c_j\neq 0\})>0$ ($i=1,2$). Then we have
\begin{equation}
\tl{T_{\sigma_i}^{\circ}}= \{(t_i,y_{k_i+1},\ldots,y_n) \ |\ t_i^{d_i} \cdot y_{k_i+1}^{c_{k_i+1}} \cdots y_n^{c_n} \times (g_{\sigma_0}|_{T_{\sigma_i}})(y_{k_2+1}, \ldots, y_n)=1\}
\end{equation}
in $\CC^*_{t_i} \times T_{\sigma_i} \simeq ( \CC^*)^{n-k_i+1}_{t_i,y_{k_i+1},\ldots,y_n}$ for $i=1, 2$. By the relation $d_2=\gcd(d_1,c_{k_1+1},\ldots,c_{k_2})$, for the integer $d=\frac{d_1}{d_2}\in \ZZ$ we have $\gcd(d,\frac{c_{k_1+1}}{d_2},\ldots, \frac{c_{k_2}}{d_2})=1$. Now let $A\in \{B\in M_{k_2-k_1+1}(\ZZ)\,|\,\det B=1\}$ be a unimodular matrix whose first row is the primitive vector $(d,\frac{c_{k_1+1}}{d_2},\ldots, \frac{c_{k_2}}{d_2})\in \ZZ^{k_2-k_1+1}$. Consider the automorphism $\Lambda_A$ of the algebraic torus $\CC^*\times (\CC^*)^{k_2-k_1}\simeq (\CC^*)^{k_2-k_1+1}$ defined by $A$:
\begin{equation}
(t_1, y_{k_1+1}, \ldots, y_{k_2})\longmapsto (t_2,z_{k_1+1},\ldots, z_{k_2}). 
\end{equation}
By this construction of $\Lambda_A$ obviously we have $t_2=t_1^dy_{k_1+1}^{\frac{c_{k_1+1}}{d_2}}\cdots y_{k_2}^{\frac{c_{k_2}}{d_2}}$. Therefore the automorphism $\Lambda_A \times \id_{T_{\sigma_2}}$ of $\CC^*\times (\CC^*)^{k_2-k_1}\times T_{\sigma_2}$ induces an isomorphism
\begin{equation}\label{AUTO}
\tl{T_{\sigma_1}^{\circ}} \simto \(\CC^*\)^{k_2-k_1}_{z_{k_1+1},\ldots,z_{k_2}}\times \tl{T_{\sigma_2}^{\circ}}.
\end{equation}
Moreover we have $\Lambda_A(\zeta_{d_1},1, \ldots,1)= (\zeta_{d_2},\beta_{k_1+1},\ldots, \beta_{k_2})$ for some $\beta_i \in \CC^*$. Since the action $\Psi_{\sigma_1}$ of the generator of $\mu_{d_1}$ on $\tl{T_{\sigma_1}^{\circ}}$ is the multiplication by the element $(\zeta_{d_1},1, \ldots,1)\in \CC^*\times T_{\sigma_1}$, the automorphism of $\(\CC^*\)^{k_2-k_1} \times \tl{T_{\sigma_2}^{\circ}}$ induced by $\Psi_{\sigma_1}$ via \eqref{AUTO} is given by 
\begin{equation} 
(z_{k_1+1},\ldots,z_{k_2},t_2, y_{k_2+1}, \ldots,y_n)
\longmapsto 
(\beta_{k_1+1}z_{k_1+1},\ldots, \beta_{k_2}z_{k_2}, \zeta_{d_2}t_2, y_{k_2+1}, \ldots, y_n). 
\end{equation}
This is obviously homotopic to $\id_{(\CC^*)^{k_2-k_1}}\times \Psi_{\sigma_2}$. Hence in the Grothendieck group $\KK_0(\HSm)$ we obtain the equality
\begin{equation} 
\chi_h([\tl{T_{\sigma_1}^{\circ}}])= \chi_h( \(\LL-1\)^{k_2 -k_1}\cdot [ \tl{T_{\sigma_2}^{\circ}}]), 
\end{equation}
from which \eqref{eq:new5-6} follows immediately. 

\noindent ({\bf Step 2}): Finally let us prove \eqref{eq:new5-6} in the general case. Let $\sigma$ be the unique cone in the dual fan $\Sigma_1$ of $\Gamma_{\infty}(f)$ whose supporting face in $\Gamma_{\infty}(f)$ is $\gamma$. Then our assumption $\gamma(\sigma_1)=\gamma(\sigma_2)=\gamma$ implies that $\relint(\sigma_{i}) \subset \relint (\sigma)$ for $i=1,2$. So there exists a continuous curve in $\relint (\sigma)$ which starts from a point in $\relint (\sigma_1)$ and ends at the one in $\relint (\sigma_2)$. Then applying ({\bf Step 1}) to each pair of two adjacent cones on it, we obtain \eqref{eq:new5-6}. This completes the proof. \qed
\end{proof}

\begin{remark}
In Proposition \ref{prp:new5-6} if $\d \sigma_1 \leq \d \sigma_2$ we can prove also a slightly stronger equality 
\begin{equation} 
\chi_h([\tl{T_{\sigma_1}^{\circ}}])= \chi_h( \(\LL-1\)^{\d \sigma_2 - \d \sigma_1}\cdot [ \tl{T_{\sigma_2}^{\circ}}]). 
\end{equation} 
Since we do not use it in this paper, we omit the proof. 
\end{remark}

For a face at infinity $\gamma \prec \Gamma_{\infty}(f)$ let $S_{\gamma} \subset \{1,2,\ldots, n\}$ be the minimal subset of $\{1,2,\ldots,n\}$ such that $\gamma \subset \RR^{S_{\gamma}}$ and set $m_{\gamma}=\sharp S_{\gamma}-\dim \gamma -1\geq 0$.

\begin{theorem}\label{thm:7-12}
Assume that $f$ is convenient and non-degenerate at infinity. Then we have the following results, where in the sums $\sum_{\gamma}$ below the face $\gamma$ of $\Gamma_{\infty}(f)$ ranges through those at infinity.
\begin{enumerate}
\item In the Grothendieck group $\KK_0(\HSm)$, we have
\begin{equation}
[H_f^{\infty}]=\chi_h(\SS_f^{\infty})=\sum_{\gamma} \chi_h((1-\LL)^{m_{\gamma}} \cdot [Z_{\Delta_{\gamma}}^*]).
\end{equation}
\item Let $\lambda \in \CC^* \setminus \{1\}$ and $k \geq 1$. Then the number of the Jordan blocks for the eigenvalue $\lambda$ with sizes $\geq k$ in $\Phi_{n-1}^{\infty} \colon H^{n-1}(f^{-1}(R);\CC) \simto H^{n-1}(f^{-1}(R);\CC)$ ($R \gg 0$) is equal to
\begin{equation}
(-1)^{n-1} \sum_{p+q=n-2+k, n-1+k} \left\{ \sum_{\gamma} e^{p,q} \left( \chi_h((1-\LL)^{m_{\gamma}}\cdot [Z_{\Delta_{\gamma}}^*])\right)_{\lambda} \right\}.
\end{equation}
\item For $k \geq 1$, the number of the Jordan blocks for the eigenvalue $1$ with sizes $\geq k$ in $\Phi_{n-1}^{\infty}$ is equal to
\begin{equation}
(-1)^{n-1} \sum_{p+q= n-2-k, n-1-k} \left\{\sum_{\gamma} e^{p,q} \left( \chi_h((1-\LL)^{m_{\gamma}}\cdot [Z_{\Delta_{\gamma}}^*])\right)_{1} 
\right\}.
\end{equation}
\end{enumerate}
\end{theorem}

\begin{proof}
\noindent (i) It suffices to rewrite Theorem \ref{thm:7-7}. Let $\gamma$ be a face at infinity of $\Gamma_{\infty}(f)$ such that $\sharp S_{\gamma}=n$. Denote by $\sigma$ the unique $(n-\d \gamma)$-dimensional cone in the dual fan $\Sigma_1$ of $\Gamma_{\infty}(f)$ whose supporting face in $\Gamma_{\infty}(f)$ is $\gamma$. Let $\Sigma$ be the smooth subdivision of $\Sigma_1$ in the proof of Theorem \ref{thm:3-5} and $\sigma_{j}$ ($1 \leq j\leq l$) the cones in $\Sigma$ such that $\relint(\sigma_{j}) \subset \relint (\sigma)$. Recall that $T_{\sigma_{j}}$ is the $(n- \d \sigma_j)$-dimensional $T$-orbit in $X_{\Sigma}$ which corresponds to $\sigma_j \in \Sigma$ and we set $T_{\sigma_{j}}^{\circ}=T_{\sigma_{j}} \setminus \overline{f^{-1}(0)}$. Then in the motivic Milnor fiber at infinity $\SS_f^{\infty}\in \M_{\CC}^{\hat{\mu}}$ of $f$ constructed by using the toric compactification $X_{\Sigma}$ of $\CC^n$, the following elements of $\M_{\CC}^{\hat{\mu}}$
\begin{equation}
(1-\LL)^{\dim \sigma_{j}-1} \cdot [ \tl{T_{\sigma_{j}}^{\circ}} ] \in \M_{\CC}^{\hat{\mu}} \hspace{10mm}(1 \leq j \leq l)
\end{equation}
are contained, where $\tl{T_{\sigma_{j}}^{\circ}}$ is the unramified Galois covering of $T_{\sigma_{j}}^{\circ}$. Let us fix $1 \leq j_0 \leq l$ such that $\d \sigma_{j_0}=\d \sigma =n- \d \gamma$. Then by Proposition \ref{prp:new5-6} for any $1 \leq j \leq l$ we have the equality
\begin{equation}\label{eq:add:5-5}
\chi_h((1-\LL)^{\dim \sigma_{j}-1} \cdot [\tl{T_{\sigma_{j}}^{\circ}}]) =(-1)^{\d \sigma_j - \d \sigma} \cdot \chi_h((1-\LL)^{n- \d \gamma -1}\cdot [ \tl{T_{\sigma_{j_0}}^{\circ}} ])
\end{equation}
in the Grothendieck group $\KK_0(\HSm)$. Combining \eqref{eq:add:5-5} with the obvious combinatorial identity $\sum_{j=1}^{l}(-1)^{\d \sigma_j - \d \sigma} =1$, we obtain a very simple formula
\begin{equation}
\sum_{j=1}^{l} \chi_h((1-\LL)^{\dim \sigma_{j}-1} \cdot [\tl{T_{\sigma_{j}}^{\circ}}]) =\chi_h((1-\LL)^{n- \d \gamma -1}\cdot [ \tl{T_{\sigma_{j_0}}^{\circ}} ]).
\end{equation}
Hence by Lemma \ref{MMFI}, for the face at infinity $\gamma$ of $\Gamma_{\infty}(f)$ such that $\sharp S_{\gamma}=n$ the equality
\begin{equation}\label{eq:5-10}
\sum_{j=1}^{l} \chi_h((1-\LL)^{\dim \sigma_{j}-1} \cdot [\tl{T_{\sigma_{j}}^{\circ}}]) =\chi_h((1-\LL)^{m_{\gamma}}\cdot [Z_{\Delta_{\gamma}}^*])
\end{equation}
holds. In the same way, we can show similar equalities also for the faces at infinity $\gamma$ of $\Gamma_{\infty}(f)$ such that $\sharp S_{\gamma}< n$. More precisely, for such $\gamma$ let $(\RR^{S_{\gamma}})^{\perp} \simeq \RR^{n -\sharp S_{\gamma}}$ be the orthogonal complement of $\RR^{S_{\gamma}} \subset \RR^n$ in $(\RR^n)^*$. Then some $\sigma_j$ in the ($n-\d \gamma )$-dimensional cone $\sigma$ associated with $\gamma$ may not satisfy the condition $(\RR^{S_{\gamma}})^{\perp} \prec \sigma_j$. We can prove a formula similar to \eqref{eq:5-10} by dividing the set of the cones $\sigma_j$ into $\{ \sigma_j \ | \ (\RR^{S_{\gamma}})^{\perp} \prec \sigma_j \}$ and $\{ \sigma_j \ | \ (\RR^{S_{\gamma}})^{\perp} \not\prec \sigma_j \}$. We omit the detail. This completes the proof of the assertion (i). The assertions (ii) and (iii) can be deduced from (i) and Theorem \ref{thm:7-6-2}. \qed
\end{proof}

\begin{remark}
Since we used a homotopy in proving \eqref{eq:add:5-5}, we can prove the equality of Theorem \ref{thm:7-12} (i) only in the Grothendieck group $\KK_0(\HSm)$ of Hodge structures. See also Remark \ref{EQH}. 
\end{remark}

Note that by using the results in Section \ref{sec:2} we can always calculate $e^{p,q}(\chi_h((1-\LL)^{m_{\gamma}}\cdot [Z_{\Delta_{\gamma}}^*]))_{\lambda}$ explicitly. Here we shall give some closed formulas for the numbers of the Jordan blocks with large sizes in $\Phi_{n-1}^{\infty}$. First let us consider the numbers of the Jordan blocks for the eigenvalues $\lambda \in \CC \setminus \{1\}$. Let $q_1,\ldots,q_l$ (resp. $\gamma_1,\ldots, \gamma_{l^{\prime}}$) be the $0$-dimensional (resp. $1$-dimensional) faces of $\Gamma_{\infty}(f)$ such that $q_i\in \Int (\RR_+^n)$ (resp. the relative interior $\relint(\gamma_i)$ of $\gamma_i$ is contained in $\Int(\RR_+^n)$). Obviously these faces are at infinity. For each $q_i$ (resp. $\gamma_i$), denote by $d_i >0$ (resp. $e_i>0$) the lattice distance $\dist(q_i, 0)$ (resp. $\dist(\gamma_i,0)$) of it from the origin $0\in \RR^n$. For $1\leq i \leq l^{\prime}$, let $\Delta_i$ be the convex hull of $\{0\}\sqcup \gamma_i$ in $\RR^n$. Then for $\lambda \in \CC \setminus \{1\}$ and $1 \leq i \leq l^{\prime}$ such that $\lambda^{e_i}=1$ we set
\begin{equation}
n(\lambda)_i
= \sharp\{ v\in \ZZ^n \cap \relint(\Delta_i) \ |\ \height (v, \gamma_i)=k\} +\sharp \{ v\in \ZZ^n \cap \relint(\Delta_i) \ |\ \height (v, \gamma_i)=e_i-k\}, 
\end{equation}
where $k$ is the minimal positive integer satisfying $\lambda=\zeta_{e_i}^{k}$ and for $v\in \ZZ^n \cap \relint(\Delta_i)$ we denote by $\height (v, \gamma_i)$ the lattice height of $v$ from the base $\gamma_i$ of $\Delta_i$. 

\begin{theorem}\label{thm:7-11}
Let $f$ be as above and $\lambda \in \CC^* \setminus \{1\}$. Then we have
\begin{enumerate}
\item The number of the Jordan blocks for the eigenvalue $\lambda$ with the maximal possible size $n$ in $\Phi_{n-1}^{\infty} \colon H^{n-1}(f^{-1}(R);\CC) \simto H^{n-1}(f^{-1}(R);\CC)$ ($R \gg 0$) is equal to $\sharp \{q_i \ |\ \lambda^{d_i}=1\}$. 
\item The number of the Jordan blocks for the eigenvalue $\lambda$ with size $n-1$ in $\Phi_{n-1}^{\infty}$ is equal to $\sum_{i \colon \lambda^{e_i}=1} n(\lambda)_i$.
\end{enumerate}
\end{theorem}

\begin{proof}
\noindent (i) By Theorem \ref{thm:7-12} (ii), the number of the Jordan blocks for the eigenvalue $\lambda \in \CC^* \setminus \{1\}$ with the maximal possible size $n$ in $\Phi_{n-1}^{\infty}$ is
\begin{align}
(-1)^{n-1} e^{n-1,n-1}(\chi_h(\SS_f^{\infty}))_{\lambda}
&=(-1)^{n-1} \sum_{i=1}^l e^{n-1,n-1}(\chi_h((1-\LL)^{n-1}\cdot [ Z_{\Delta_{q_i}}^* ]))_{\lambda}\\
&= \sum_{i=1}^l e^{0,0}(\chi_h( [Z_{\Delta_{q_i}}^* ]))_{\lambda}.
\end{align}
Note that $Z_{\Delta_{q_i}}^*$ is a finite subset of $\CC^*$ consisting of $d_i$ points. Then (i) follows from 
\begin{equation}
\sum_{i=1}^l e^{0,0}(\chi_h([Z_{\Delta_{q_i}}^* ]))_{\lambda}=\sharp \{q_i \ |\ \lambda^{d_i}=1\}.
\end{equation}
The assertion (ii) can be proved similarly by expressing $e^{n-1,n-2}(\chi_h(\SS_f^{\infty}))_{\lambda} +e^{n-2,n-1}(\chi_h(\SS_f^{\infty}))_{\lambda}$ in terms of the $1$-dimensional faces at infinity $\gamma_i$ of $\Gamma_{\infty}(f)$. \qed
\end{proof}

\begin{example}
Let $f(x,y)\in \CC[x,y]$ be a convenient polynomial whose Newton polyhedron at infinity $\Gamma_{\infty}(f)$ has the following shape.

\vspace{1mm}
\begin{center}
\includegraphics[scale=.8]{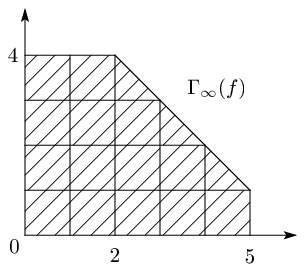}

Figure 13
\end{center}

\vspace{1mm}
Assume moreover that $f$ is non-degenerate at infinity. Then by Libgober-Sperber's theorem (Theorem \ref{thm:3-5}) the characteristic polynomial $P(\lambda)$ of $\Phi_1^{\infty} \colon H^1(f^{-1}(R);\CC) \simto H^1(f^{-1}(R);\CC)$ ($R \gg 0$) is given by
\begin{equation}
P(\lambda)=(\lambda-1)(\lambda^4-1)(\lambda^6-1)^3.
\end{equation}
In particular, the total multiplicity of the roots $-1$ in $P( \lambda )=0$ is $4$. For $\lambda \in \CC$, denote by $H^1(f^{-1}(R);\CC)_{\lambda}$ the generalized $\lambda$-eigenspace of the monodromy operator $\Phi_1^{\infty}$ at infinity. First, by the monodromy theorem the restriction of $\Phi_1^{\infty}$ to $H^1(f^{-1}(R);\CC)_1 \simeq \CC^5$ is semisimple. Moreover by Theorem \ref{thm:7-11} (i) the Jordan normal form of the restriction of $\Phi_1^{\infty}$ to $H^1(f^{-1}(R);\CC)_{-1} \simeq \CC^4$ is
\begin{equation}
\begin{pmatrix}
-1 &1&0&0\\
0&-1&0&0\\
0&0&-1&0\\
0&0&0&-1
\end{pmatrix}.
\end{equation}
In the same way, we can show that for $\lambda=\zeta_6, \sqrt{-1}, \zeta_3, \zeta_3^2, -\sqrt{-1}, \zeta_6^5$ the restriction of $\Phi_1^{\infty}$ to $H^1(f^{-1}(R);\CC)_{\lambda}$ is semisimple.
\end{example}

Next we consider the number of the Jordan blocks for the eigenvalue $1$ in $\Phi_{n-1}^{\infty}$. By Proposition \ref{prp:2-19}, we can rewrite Theorem \ref{thm:7-12} (iii) as follows. Denote by $\Pi_f$ the number of the lattice points on the $1$-skeleton of $\partial \Gamma_{\infty}(f) \cap \Int (\RR^n_+)$.

\begin{theorem}\label{thm:7-16}
In the situation as above, the number of the Jordan blocks for the eigenvalue $1$ with the maximal possible size $n-1$ in $\Phi_{n-1}^{\infty}$ is $\Pi_f$.
\end{theorem}

\begin{proof}
For a face at infinity $\gamma \prec \Gamma_{\infty}(f)$, denote by $\Pi(\gamma)$ the number of the lattice points on the $1$-skeleton of $\gamma$. Since for each face at infinity $\gamma \prec \Gamma_{\infty}(f)$ we have $\Pi(\Delta_{\gamma})_1-1=\Pi(\gamma)$ (for the definition of $\Pi(\Delta_{\gamma})_1$, see Section \ref{sec:2}), the assertion follows from Theorem \ref{thm:7-12} (iii) and Proposition \ref{prp:2-19}. \qed
\end{proof}

For a face at infinity $\gamma \prec \Gamma_{\infty}(f)$, denote by $l^*(\gamma)$ the number of the lattice points on the relative interior $\relint(\gamma)$ of $\gamma$. Then by Theorem \ref{thm:7-12} (iii) and Proposition \ref{prp:new}, we also obtain the following result.

\begin{theorem}\label{thm:7-17}
In the situation as above, the number of the Jordan blocks for the eigenvalue $1$ with size $n-2$ in $\Phi_{n-1}^{\infty}$ is equal to $2 \sum_{\gamma} l^*(\gamma)$, where $\gamma$ ranges through the faces at infinity of $\Gamma_{\infty}(f)$ such that $\d \gamma =2$ and $\relint(\gamma) \subset \Int (\RR^n_+)$. In particular, this number is even.
\end{theorem}

From now on, we assume that any face at infinity $\gamma \prec \Gamma_{\infty}(f)$ is prime in the sense of Definition \ref{dfn:2-16} (i) and rewrite Theorem \ref{thm:7-12} (ii) and (iii) more explicitly. First, recall that by Proposition \ref{prp:2-15} for $\lambda \in \CC^* \setminus \{1\}$ and a face at infinity $\gamma\prec \Gamma_{\infty}(f)$ we have $e^{p,q}(Z_{\Delta_{\gamma}}^*)_{\lambda}=0$ for any $p,q \geq 0$ such that $p+q >\d \Delta_{\gamma}-1=\dim \gamma$. So the non-negative integers $r \geq 0$ such that $\sum_{p+q=r}e^{p,q}(Z_{\Delta_{\gamma}}^*)_{\lambda}\neq 0$ are contained in the closed interval $[0,\d \gamma]\subset \RR$.

\begin{definition}
For a face at infinity $\gamma \prec \Gamma_{\infty}(f)$ and $k \geq 1$, we define a finite subset $J_{\gamma,k}\subset [0,\d \gamma] \cap \ZZ$ by
\begin{equation}
J_{\gamma,k}=\{0 \leq r\leq \d \gamma \ |\ n-2+k \equiv r \mod 2\}.
\end{equation}
For each $r\in J_{\gamma,k}$, set
\begin{equation}
d_{k,r}=\dfrac{n-2+k-r}{2}\in \ZZ_+.
\end{equation}
\end{definition}

Since for any face at infinity $\gamma \prec \Gamma_{\infty}(f)$ the polytope $\Delta_{\gamma}$ is pseudo-prime in the sense of Definition \ref{dfn:2-16} (ii), by Corollary \ref{cor:2-18} for $\lambda \in \CC^* \setminus \{1\}$ and an integer $r \geq 0$ such that $r\in [0,\d \gamma] $ we have
\begin{equation}
\sum_{p+q=r}e^{p,q}(\chi_h([Z_{\Delta_{\gamma}}^*]))_{\lambda}=(-1)^{\d \gamma +r+1} \sum_{\begin{subarray}{c} \Gamma\prec \Delta_{\gamma} \\ \d \Gamma=r+1\end{subarray}} \left\{ \sum_{\Gamma^{\prime} \prec \Gamma} (-1)^{\d \Gamma^{\prime}} \tl{\varphi}_{\lambda}(\Gamma^{\prime})\right\}.
\end{equation}
For simplicity, we denote this last integer by $e(\gamma,\lambda)_r$. Then by Theorem \ref{thm:7-12} (ii) we obtain the following result.

\begin{theorem}\label{thm:7-15}
In the situation as above, let $\lambda \in \CC^* \setminus\{1\}$ and $k\geq 1$. Then the number of the Jordan blocks for the eigenvalue $\lambda$ with sizes $\geq k$ in $\Phi_{n-1}^{\infty} \colon H^{n-1}(f^{-1}(R);\CC) \simto H^{n-1}(f^{-1}(R);\CC)$ ($R\gg 0$) is equal to
\begin{equation}
(-1)^{n-1}\sum_{\gamma} \left\{ \sum_{r \in J_{\gamma, k}} (-1)^{d_{k,r}} \binom{m_{\gamma}}{d_{k,r}} \cdot e(\gamma,\lambda)_r + \sum_{r \in J_{\gamma, k+1}} (-1)^{d_{k+1,r}} \binom{m_{\gamma}}{d_{k+1,r}} \cdot e(\gamma,\lambda)_r\right\},
\end{equation}
where in the sum $\sum_{\gamma}$ the face $\gamma$ of $\Gamma_{\infty}(f)$ ranges through those at infinity (we used also the convention $\binom{a}{b}=0$ ($0 \leq a <b$) for binomial coefficients).
\end{theorem}

By combining the proof of \cite[Theorem 5.6]{D-K} and Proposition \ref{prp:2-17-2} with Theorem \ref{thm:7-12} (iii), if any face at infinity $\gamma \prec \Gamma_{\infty}(f)$ is prime we can also explicitly describe the number of the Jordan blocks for the eigenvalue $1$ in $\Phi_{n-1}^{\infty}$.

Finally to end this section, we prove a global analogue of the Steenbrink conjecture proved by Varchenko-Khovanskii \cite{K-V} and Saito \cite{Saito-3}. We return to the general case.

\begin{definition}{\rm (Sabbah \cite{Sabbah} and Steenbrink-Zucker \cite{S-Z})} As a Puiseux series, we define the spectrum at infinity $\sp_f^{\infty}(t)$ of $f$ by
\begin{equation}
\sp_f^{\infty}(t)
= \sum_{\beta \in (0,1] \cap \QQ} \left[ \sum_{i=0}^{n-1} (-1)^{n-1}\left\{ \sum_{q \geq 0} e^{i,q}([H_f^{\infty}])_{\exp(2\pi \sqrt{-1}\beta)}\right\} t^{i+\beta}\right]+(-1)^{n}t^n.
\end{equation}
\end{definition}

When $f$ is tame at infinity, by Theorem \ref{cor:7-6-2} we can easily prove that the support of $\sp_f^{\infty}(t)$ is contained in the open interval $(0,n)$ and has the symmetry
\begin{equation}
\sp_f^{\infty}(t)=t^n \sp_f^{\infty}\( \frac{1}{t}\)
\end{equation}
with center at $\frac{n}{2}$. From now on, we assume that $f$ is convenient and non-degenerate at infinity. In order to describe $\sp_f^{\infty}(t)$ by $\Gamma_{\infty}(f)$, for each face at infinity $\gamma$ of $\Gamma_{\infty}(f)$ let $s_{\gamma}=\sharp S_{\gamma}\in \ZZ_{\geq 1}$ be the dimension of the minimal coordinate plane containing $\gamma$ and set $\Cone(\gamma)=\RR_+\gamma$. Next, let $h_f \colon \RR_+^n \longrightarrow \RR$ be the continuous function on $\RR_+^n$ which is linear on each cone $\Cone(\gamma)$ and satisfies the condition $h_f|_{\partial \Gamma_{\infty}(f) \cap \Int(\RR_+^n)} \equiv 1$. For a face at infinity $\gamma$ of $\Gamma_{\infty}(f)$, let $L_{\gamma}$ be the semigroup $\Cone(\gamma) \cap \ZZ_+^n$ and define its Poincar{\'e} series $P_{\gamma}(t)$ by
\begin{equation}
P_{\gamma}(t)=\sum_{\beta \in \QQ_+} \sharp \{ v \in L_{\gamma} \ |\ h_f(v) =\beta\} t^{\beta}.
\end{equation}

\begin{theorem}\label{thm:7-19}
Assume that $f$ is convenient and non-degenerate at infinity. Then we have
\begin{equation}
\sp_f^{\infty}(t)=\sum_{\gamma} (-1)^{n-1-\d \gamma} (1-t)^{s_{\gamma}}P_{\gamma}(t) +(-1)^n,
\end{equation}
where in the above sum $\gamma$ ranges through the faces at infinity of $\Gamma_{\infty}(f)$.
\end{theorem}

\vspace*{-3mm}
\begin{proof}
For $\beta \in (0,1] \cap \QQ$ and a face at infinity $\gamma$ of $\Gamma_{\infty}(f)$, set
\begin{equation}
P_{\gamma, \beta}(t)=\begin{cases}\sum_{i=0}^{\infty} \sharp \{ v \in L_{\gamma} \ |\ h_f(v) =i+ \beta\} t^{i+ \beta}  & (0< \beta <1), \\
\sum_{i=0}^{\infty} \sharp \{ v \in L_{\gamma} \ |\ h_f(v) =i \} t^{i} & (\beta = 1)
\end{cases}
\end{equation}
so that we have
\begin{equation}
\sum_{\beta \in (0,1] \cap \QQ} P_{\gamma, \beta}(t)=P_{\gamma}(t).
\end{equation}
Then for $\beta\in \QQ$ such that $0<\beta<1$ and a face at infinity $\gamma$ of $\Gamma_{\infty}(f)$, by Theorem \ref{thm:2-14} and \eqref{E:sym} we have
\begin{eqnarray}
\lefteqn{\sum_{i=0}^{n-1}(-1)^{n-1}\left\{ \sum_{q \geq 0} e^{i,q}(\chi_h((1-\LL)^{m_{\gamma}}\cdot [Z_{\Delta_{\gamma}}^*] ))_{\exp(2\pi \sqrt{-1}\beta)}\right\} t^{i+\beta}}\nonumber \\
&=& (-1)^{n-1-\d \gamma} t^{\beta}(1-t)^{m_{\gamma}}\sum_{i \geq 0} \varphi_{\exp(2\pi \sqrt{-1}\beta), \d \gamma +1-i}(\Delta_{\gamma}) t^i\\
&=& (-1)^{n-1-\d \gamma} t^{\beta}(1-t)^{s_{\gamma}-1-\d \gamma} \dfrac{1}{t}\sum_{i \geq 0} \psi_{\exp(-2\pi \sqrt{-1}\beta), i+1}(\Delta_{\gamma})t^{i+1}\\
&=& (-1)^{n-1-\d \gamma} t^{\beta} (1-t)^{s_{\gamma}+1} \sum_{k \geq 1}l(k\Delta_{\gamma})_{\exp(-2\pi \sqrt{-1}\beta)}t^{k-1}\\
&=& (-1)^{n-1-\d \gamma} t^{\beta}(1-t)^{s_{\gamma}}\nonumber \\
& & \ \times (1-t) \left\{ l(\Delta_{\gamma})_{\exp(-2\pi \sqrt{-1}\beta)}+l(2\Delta_{\gamma})_{\exp(-2\pi \sqrt{-1}\beta)} t+\cdots \right\} \\
&=& (-1)^{n-1-\d \gamma} (1-t)^{s_{\gamma}} P_{\gamma, \beta}(t).\end{eqnarray}
Therefore, the assertion for the non-integral part of $\sp_f^{\infty}(t)$ follows immediately from Theorem \ref{thm:7-12} (i). Moreover, the integral part
\begin{equation}
\sum_{i=0}^{n-1}(-1)^{n-1} \left\{ \sum_{q \geq 0} e^{i,q}(\chi_h(\SS_f^{\infty}))_1 \right\} t^{i+1}+(-1)^{n}t^n
\end{equation}
of $\sp_f^{\infty}$ is calculated as follows. For a face at infinity $\gamma$ of $\Gamma_{\infty}(f)$, we have
\begin{eqnarray}
\lefteqn{ \sum_{i=0}^{n-1} (-1)^{n-1} \left\{ \sum_{q\geq 0}e^{i,q}( \chi_h((1-\LL)^{m_{\gamma}}\cdot [Z_{\Delta_{\gamma}}^*]))_1\right\} t^{i+1}} \nonumber \\
&=& (-1)^{n-1-\d \gamma} (1-t)^{m_{\gamma}} \sum_{i \geq 0} \left\{ (-1)^i \binom{\d \gamma+1}{i+1} +\varphi_{1, \d \gamma +1-i}(\Delta_{\gamma})\right\} t^{i+1}\\
&=& (-1)^{n-1-\d \gamma} (1-t)^{s_{\gamma}-1-\d \gamma} \left\{ -(1-t)^{\d \gamma +1}+1 +\sum_{i \geq 0}\psi_{1,i+1}(\Delta_{\gamma})t^{i+1}\right\}\\
&=& (-1)^{n-1-\d \gamma} (1-t)^{s_{\gamma}-1-\d \gamma}\nonumber \\
& & \ \times \bigg[ -(1-t)^{\d \gamma+1} +(1-t)^{\d \gamma +2} \left\{ l(0)_1+l(\Delta_{\gamma})_1t+l(2\Delta_{\gamma})_1t^2+\cdots \right\} \bigg]\\
&=& (-1)^{n-\d \gamma} (1-t)^{s_{\gamma}} +(-1)^{n-1-\d \gamma} (1-t)^{s_{\gamma}}P_{\gamma, 1}(t).
\end{eqnarray}
Summing up these terms over the faces at infinity $\gamma$ of $\Gamma_{\infty}(f)$, we obtain
\begin{equation}
\sum_{i=0}^{n-1} (-1)^{n-1} \left\{ \sum_{q \geq 0} e^{i,q}(\chi_h(\SS_f^{\infty}))_1 \right\}t^{i+1}=\sum_{\gamma} (-1)^{n-1-\dim \gamma}(1-t)^{s_{\gamma}}P_{\gamma, 1}(t)+(-1)^{n+1}t^n+(-1)^n.
\end{equation}
This completes the proof. \qed
\end{proof}


\appendix
\section{Appendix by Claude Sabbah}

In this appendix, we prove Theorems \ref{cor:7-6-2} and \ref{thm:7-6-2} of the main article.

\subsection{Symmetry of Hodge numbers}\label{sec:mhm}

Let $U$ be a smooth affine complex variety  and let $\cO_U^\rH$ denote the mixed Hodge module also denoted by $\QQ_U^\rH$ in \cite{Saito-2} (we change the notation because we will mainly work with filtered $\cD$-modules). We set $n=\dim U$ and $m=n-1$. Let $\DD$ be the duality functor of algebraic mixed Hodge modules. We have $\DD\cO_U^\rH\simeq\cO_U^\rH(n)$.

Let $f:U\to\Afu$ be a regular function. We denote by $f_*,f_!$ be the push-forward and proper push-forward functors (mainly used at the level of filtered $\cD$-modules), where $f_!=\DD f_*\DD$ (cf. \cite[(4.3.5)]{Saito-2}).

Let $t$ be the coordinate on $\Afu$. We will also use the nearby cycle functor $\psi_{1/t}$, that we decompose as $\psi_{1/t}=\psi_{1/f,1}\oplus\psi_{1/f,\neq1}$ with respect to the eigenvalues of the monodromy. We have the following commutation relations in $\textup{MHM}(\Afu)$ (cf. \cite[Prop.\kern2pt 2.6]{Saito-2}):
\[
\psi_{1/t}\DD=(\DD\psi_{1/t})(1).
\]

According to the previous relations, we have
\[
\psi_{1/t}(\cH^0f_!\cO_U^\rH)\simeq\psi_{1/t}(\DD\cH^0f_*\DD\cO_U^\rH)\simeq\DD\big(\psi_{1/t}(\cH^0f_*\cO_U^\rH)\big)(-m).
\]

We denote by $h_!^{p,q}$ the Hodge numbers of the left-hand term, and by $h_*^{p,q}$ those of $\psi_{1/t}(\cH^0f_*\cO_U^\rH)$. We then get
\[\tag*{$(!*)$}
\forall p,q\in\ZZ,\qquad h_!^{p,q}=h_*^{m-p,m-q},
\]
or equivalently, for each eigenvalue $\alpha\in\exp(2\pi i\QQ)$,
\[\tag*{$(!*)_\alpha$}
\forall p,q\in\ZZ,\qquad h_{!,\alpha}^{p,q}=h_{*,\alpha^{-1}}^{m-p,m-q}=h_{*,\alpha}^{m-q,m-p},
\]
according to the behaviour of eigenvalues by duality and complex conjugation, if we note that $\alpha^{-1}=\overline{\alpha}$ for $\alpha\in\exp(2\pi i\QQ)$. From now on, we assume that $f$ is cohomologically tame, in the sense of \cite{Sabbah-2}. If $U=\bbA^n$ and $n\geqslant2$, $\psi_{1/t,\neq1}(\cH^kf_!\cO_U^\rH)=\nobreak0$ for $k\neq0$ (see \cite[Rem.~10.3]{Sabbah-2} for $\cH^kf_*$ and use duality). Therefore, using the notation $e^{p,q}$ of (2.8) and (2.9) in the main part of the article, we have $e^{p,q}_{\neq1}=h_{!,\neq1}^{p,q}$, and we wish to show the symmetry $h_{!,\alpha}^{p,q}=h_{!,\alpha}^{m-p,m-q}$ for $\alpha\neq1$. If $\alpha=1$, the point is to show the symmetry $h_{!,1}^{p,q}=h_{!,1}^{m-1-p,m-1-q}$, since $f_!\cO_U^\rH$ has cohomology in degrees $0$ and $m$ at most and $\cH^m$ has rank one. By $(!*)_\alpha$, these symmetries are equivalent to $h_{*,\alpha}^{p,q}=h_{*,\alpha}^{m-p,m-q}$ for $\alpha\neq1$, and $h_{*,1}^{p,q}=h_{*,1}^{m+1-p,m+1-q}$. Both are a direct consequence of the following proposition, since $\wt\rN$ is a morphism of type $(-1,-1)$.

\begin{proposition}\label{prop:main}
The weight filtration on $\psi_{1/t,1}(\cH^0f_*\cO_U^\rH)$ (resp.\ on $\psi_{1/t,\neq1}(\cH^0f_*\cO_U^\rH)$) is equal to the monodromy filtration of the nilpotent part of the monodromy, centered at $m+1$ (resp.\ $m$).
\end{proposition}

Notice also that Theorem \ref{thm:7-6-2} of the main article is a consequence of this statement. One can obtain the proposition as a consequence of Theorem 13.1 in \cite{Sabbah-2}, but we will propose another proof, which avoids the main results of \cite{Sabbah-2} related to Fourier transform, Brieskorn lattices and spectrum at infinity.

Let us first treat the case $\alpha\neq1$. Let $F:X\to\CC$ be a compactification of~$f$ with no vanishing cycle for $\QQ_U$ on $X\moins U$ (tameness), and let $\IC_X(\QQ_U)$ be the intersection complex of~$X$. It corresponds to a pure Hodge module $j_{!*}\cO_U^\rH$, according to M.~Saito, and $\cH^0F_*(j_{!*}\cO_U^\rH)$ is pure. Moreover, we have two morphisms in $\MHM(\Afu)$
\[
\cH^0f_!\cO_U^\rH\to\cH^0F_*(j_{!*}\cO_U^\rH)\to\cH^0f_*\cO_U^\rH
\]
and for each morphism, the kernel and cokernel (in $\MHM(\Afu)$) are constant mixed Hodge modules.

It follows that the computation of $\psi_{1/t,\neq1}\cH^0f_*\cO_U^\rH$ or $\psi_{1/t,\neq1}\cH^0f_!\cO_U^\rH$ (where~$t$ is the coordinate on $\Afu$) coincides with the computation of $\psi_{1/t,\neq1}\cH^0F_*(j_{!*}\cO_U^\rH)$. There, we can apply the properties of pure Hodge modules and get that the weight filtration is the monodromy filtration shifted by $m$, according to \cite{Saito-1}. The case $\alpha=1$ will occupy the next sections.

\subsection{A preliminary result}
Let $H$ be a finite dimensional vector space equipped with a nilpotent endomorphism $\wt\rN$. We denote by $\rM(\wt\rN,H)_\bbullet$ the monodromy filtration of~$\wt\rN$ on $H$ (centered at $0$), so that $\wt\rN(\rM(\wt\rN,H)_k)\subset \rM(\wt\rN,H)_{k-2}$ for  any $k\in\ZZ$ and, for any $\ell\in\NN^*$, $\wt\rN^\ell$ induces an isomorphism $\gr_\ell^{\rM(\wt\rN,H)}H\isom\gr_{-\ell}^{\rM(\wt\rN,H)}H$.

The space $H/\Im \wt\rN$ is naturally decomposed into primitive subspaces $\rP_0(H,\wt\rN)\oplus\cdots\oplus \rP_\ell(H,\wt\rN)\oplus\cdots$, and the filtration induced by $\rM(\wt\rN,H)_\bbullet$ on $H/\Im \wt\rN$ is the filtration by the degree of the primitive part. The following is straightforward, by using the Jordan normal form for instance.

\begin{lemma}\label{lem:W=M}
Let $L_\bbullet H$ be an increasing exhaustive filtration of $H$ such that $L_{-2}=\nobreak0$ and $L_{-1}H=\Im \wt\rN$. Then the following properties are equivalent:
\begin{enumeratea}\renewcommand{\labelenumi}{{\rm (\alph{enumi})}}
\item\label{lem:elem1}
$\rM(\wt\rN,H)_\bbullet$ is equal to the monodromy filtration of $\wt\rN$ relative to $L_\bbullet H$,
\item\label{lem:elem2}
for $k\geqslant0$, $L_kH=\rM(\wt\rN,H)_k+\Im \wt\rN$.
\end{enumeratea}
\end{lemma}

\subsection{Vanishing of hypercohomology}
Let $M$ be a regular holonomic $\cD$-module on the affine line $\Afu$ with coordinate $t$. The following operation defines a new regular holonomic $\cD$-module $\wt M$ such that the de~Rham hypercohomology $\bH^*(\Afu,\DR(\wt M))$ is zero. Note that, because we work with regular holonomic $\cD$-modules, there is no difference between the algebraic and the analytic de~Rham hypercohomologies. Working with $\Clt$-modules, this amounts to asking that $\partial_t:\wt M\to\wt M$ is bijective. This can easily be realized by the following operation:
\[
\wt M=\CC[\partial_t,\partial_t^{-1}]\otimes_{\CC[\partial_t]}M,
\]
but this operation is not easily extended to mixed Hodge modules, which is our main purpose.

We will now consider the $\Clt$-module $\CC[\partial_t,\partial_t^{-1}]$ as a mixed Hodge module, and we will denote it $\CC[\partial_t,\partial_t^{-1}]^\rH$. It is constructed as follows. Firstly, as a $\Clt$-module, we have a natural exact sequence
\[
0\to\CC[\partial_t]\to\CC[\partial_t,\partial_t^{-1}]\to\CC[t]\to0
\]
by presenting $\CC[\partial_t,\partial_t^{-1}]$ as $\Clt/(t\partial_t+1)$. Denoting by $j:(\Afu)^*=\Afu\moins\{0\}\hto\Afu$ the open inclusion and by $i:\{0\}\hto\Afu$ the complementary closed inclusion, it corresponds to the mixed Hodge module $j_!\cO_{(\Afu)^*}^\rH$. The previous exact sequence is the weight exact sequence:
\[
W_0(j_!\cO_{(\Afu)^*}^\rH)=i_*\QQ_{0}^\rH,\quad \gr_1^W(j_!\cO_{(\Afu)^*}^\rH)=\cO_{\Afu}^\rH.
\]
(Recall that, in the theory of mixed Hodge modules, $\cO_{\Afu}^\rH$ has weight $\dim\Afu=1$).

If $M$ is a regular holonomic $\cD_{\Afu}$-module, we thus set
\[
\wt M=\cH^0s_*(M\boxtimes j_!\cO_{(\Afu)^*})
\]
where $s:\Afu\times\Afu\to\Afu$ is the sum function $(x,y)\mto x+y$ and the direct image is taken in the sense of $\cD$-modules. Similarly, if $M$ is a mixed Hodge module, we can regard the previous definition within the frame of mixed Hodge modules and define $\wt M$ as a mixed Hodge module. We have a natural morphism $M\to\wt M$, whose kernel and cokernel are constant mixed Hodge modules.

Let us assume that $M$ is a pure Hodge module on $\Afu$, of weight $w$. Then its image in $\wt M$ is also a pure of weight $w$, and we still denote it by $M$. In other words, we will assume that $M$ has no constant submodule. Then we have an exact sequence in $\MHM(\Afu)$:
\begin{equation}\label{eq:*}
0\to M\to\wt M\to M''\to0
\end{equation}
and $M''$ is constant and has weights $\geqslant w+1$.

\subsection{Nearby cycles}
The exact sequence \eqref{eq:*} induces an exact sequence of mixed Hodge structures after taking nearby cycles at infinity:
\[
0\to\psi_{1/t,1}M\to\psi_{1/t,1}\wt M\to \psi_{1/t,1}M''\to0.
\]
The weight filtration on $\psi_{1/t,1}M$ is the monodromy filtration of the nilpotent part~$\rN$ of the monodromy at infinity, centered at $w-1$, that we write $\rM(\rN,\psi_{1/t,1}M)[w-\nobreak1]_\bbullet$. The weight filtration $W_\bbullet$ of $\psi_{1/t,1}\wt M$ is the monodromy filtration of $\wt\rN$ on $\psi_{1/t,1}\wt M$ relative to the filtration $L_\bbullet$ induced by $W_{\bbullet+1}\wt M$. Lastly, $M''$ is constant, so $\rN''=0$ on $\psi_{1/t,1}M''$ and the weight filtration $W_\bbullet \psi_{1/t,1}M''$ is equal to $\psi_{1/t,1}W_{\bbullet+1}M''$.

\begin{proposition}\label{prop:W=M}
Under these assumptions, the weight filtration $W_\bbullet$ on $\psi_{1/t,1}\wt M$ is equal to the (absolute) monodromy filtration $\rM(\wt\rN)[w]_\bbullet$ of $\wt\rN$ centered at~$w$, and~$L_\bbullet$ is given by Lemma \ref{lem:W=M}, up to a shift by $w$.
\end{proposition}

\begin{proof}
We will show that the filtration $L_\bbullet\psi_{1/t,1}\wt M$ defined above satisfies the assumption of Lemma \ref{lem:W=M} (up to a shift by $w$) and that $\rM(\wt\rN)[w]_\bbullet$ is the weight filtration of $\psi_{1/t,1}\wt M$. Therefore, the property \ref{lem:W=M}\eqref{lem:elem1} will be fulfilled, and thus the filtration $L_\bbullet$ satisfies \ref{lem:W=M}\eqref{lem:elem2}.

Let us first give some properties of the filtration $L_\bbullet$. In the exact sequence \eqref{eq:*}, the weight filtration of $\wt M$ satisfies $W_{w-1}\wt M=0$, $W_w\wt M=M$ and $M''$ has weights $\geqslant w+1$. Each $W_{k+1}M''$ ($k\geqslant w$) is a constant Hodge module, which is completely determined by $W_k\bH^{-1}(\Afu,\DR(M''))$ (where the $\cD$-module convention is used for the de~Rham complex, that is, $\DR(M'')$ has terms in degrees $-1$ and $0$). Since $\wt M$ has no global hypercohomology, we have an isomorphism of mixed Hodge structures
\[
\bH^{-1}(\Afu,\DR M'')\isom \bH^0(\Afu,\DR M).
\]
The mixed Hodge structure on $\bH^0(\Afu,\DR M)$ is described as follows. It has weights $\geqslant w$. Let us denote by $\cM_{\min}$ the minimal extension of $M$ by the inclusion $j:\Afu\hto\nobreak\PP^1$ and by $\cM$ the maximal extension $j_*M$. Then $\cM_{\min}$ is a pure Hodge module of weight $w$ on $\PP^1$, and
\[
\gr_w^W\bH^0(\Afu,\DR M)=W_w\bH^0(\Afu,\DR M)=\bH^0(\PP^1,\DR\cM_{\min})\subset \bH^0(\Afu,\DR M).
\]
The quotient Hodge structure $\bH^0(\Afu,\DR M)/W_w\bH^0(\Afu,\DR M)$ is identified with $\bH^0(\PP^1,\DR (\cM/\cM_{\min}))$. Note that $\cM/\cM_{\min}$ is supported at infinity, and is identified with the direct image by the inclusion $\infty\hto\PP^1$ of $\phi_{1/t,1}(\cM/\cM_{\min})$. Moreover, $\phi_{1/t,1}\cM_{\min}$ is identified with $\Im\rN:\psi_{1/t,1}\cM_{\min}\to\psi_{1/t,1}\cM_{\min}(-1)$  (\cf\cite[Lemme 5.1.4]{Saito-1}) and, since $\psi_{1/t,1}\cM_{\min}\to\psi_{1/t,1}\cM$ and $\var:\phi_{1/t,1}\cM\to\psi_{1/t,1}\cM(-1)$ are isomorphisms compatible with $\rN$, we get an identification of $\phi_{1/t,1}(\cM/\cM_{\min})$ with $\coker\rN:\psi_{1/t,1}M\to\psi_{1/t,1}M(-1)$. The graded pieces are thus given by
\[
\gr_{w+k+1}^W\bH^0(\Afu,\DR M)\isom \rP_k(\rN,\psi_{1/t,1}M)(-1),\quad\forall k\geqslant0.
\]
(Recall that $\psi_{1/t,1}M$ is a mixed Hodge module having weight filtration given by $W_\bbullet\psi_{1/t,1}M=\rM(\rN)[w-1]_\bbullet$; then $\psi_{1/t,1}M(-1)$ is a mixed Hodge module with $W_\bbullet\psi_{1/t,1}M(-1)=\rM(\rN)[w+1]_\bbullet$.)

This computation gives the weight filtration on $\bH^{-1}(\Afu,\DR M'')$, and thus on~$M''$ since this is a constant Hodge module:
\begin{align*}
W_{w+1}M''&=\bH^0(\PP^1,\DR\cM_{\min})\otimes_\CC\CC[t]^\rH,\\
\gr_{w+k+2}^WM''&\simeq \rP_k(\rN,\psi_{1/t,1}M)(-1)\otimes_\CC\CC[t]^\rH,\quad\forall k\geqslant0.
\end{align*}
As a consequence, we get
\begin{align}
\gr^L_w\psi_{1/t,1}\wt M=\gr^W_w\psi_{1/t,1}M''&\simeq\bH^0(\PP^1,\DR\cM_{\min}),\notag\\
\gr_{w+k+1}^L\psi_{1/t,1}\wt M=\gr_{w+k+1}^W\psi_{1/t,1}M''&\simeq \rP_k(\rN,\psi_{1/t,1}M)(-1),\quad\forall k\geqslant0.\label{eq:Pk}
\end{align}
On the other hand,
\[
L_{w-1}\psi_{1/t,1}\wt M=\gr^L_{w-1}\psi_{1/t,1}\wt M=\psi_{1/t,1}M.
\]

\subsubsection*{Proof that $L_{w-1}\psi_{1/t,1}\wt M=\Im\wt\rN$}

Since $\rN''=0$, we have $\Im\wt\rN\subset\psi_{1/t,1}M=L_{w-1}\psi_{1/t,1}\wt M$. We will prove equality by an argument of Fourier transform. Recall that the Fourier transform $\Fou M$ of $M$ is a $\Cltau$-module, through the correspondence $\Cltau\isom\Clt$, $\tau\mto\partial_t$, $\partial_\tau\mto-t$.

\begin{lemma}
There is a functorial isomorphism
\begin{equation}\label{eq:**}
(\phi_{\tau,1}\Fou M,\Fou \rN)\isom(\psi_{1/t,1}M,\rN)
\end{equation}
for any regular holonomic $\Clt$-module $M$.
\end{lemma}

\begin{proof}
This is ``well-known''. The proof of \cite[Prop.\,4.1(ivb)]{Bibi05b} can be adapted to $\cD$-modules to show that a similar assertion holds on the product space $\PP^1_t\times\Afu_\tau$ for the pull-back $p^*M$ of $M$ twisted by the exponential $\cD$-module $\cE^{-t\tau}$ (kernel of the Laplace transform). Applying direct image by the projection $q:\PP^1_t\times\Afu_\tau\to\Afu_\tau$ and the compatibility of the functor $\phi_{\tau,1}$ with direct images (\cf \eg \cite{M-S86b}), we obtain \eqref{eq:**}, since $\Fou M$ can also be computed as $\cH^0p_*(p^*M\otimes\cE^{-t\tau})$. \qed
\end{proof}

Notice that $\Fou\wt M$ is the localization with respect to $\tau$ of $\Fou M$. Then the natural map $\psi_{1/t,1}M\to\psi_{1/t,1}\wt M$ is identified with the natural morphism (variation) $\phi_{\tau,1}\Fou M\to\psi_{\tau,1}\Fou M$ via the commutative diagram:
\[
\xymatrix{
\psi_{1/t,1}M\ar[d]&\ar[l]_-\sim\phi_{\tau,1}\Fou M\ar[r]^-{\var_\tau}\ar[d]&\psi_{\tau,1}\Fou M\ar[d]^\wr\\
\psi_{1/t,1}\wt M&\ar[l]_-\sim\phi_{\tau,1}\Fou\wt M\ar[r]_-{\var_\tau}^-\sim&\psi_{\tau,1}\Fou \wt M
}
\]

The point is now that $M$ is a semi-simple $\Clt$-module, as it underlies a pure Hodge module. Moreover, the natural morphism from $\Fou M$ to its localization $\Fou\wt M$ is injective. Therefore, $\Fou M$ is a semi-simple $\Cltau$-module and has no submodule supported on $\tau=0$. Hence, the dual $\Cltau$-module satisfies the same properties, and therefore is included in its localization at $\tau=0$. As a consequence, $\Fou M$ is a minimal extension at $\tau=0$ (\ie has no sub or quotient module supported at $\tau=0$), which implies that $\phi_{\tau,1}\Fou M\simeq\Im(\Fou \rN:\psi_{\tau,1}\Fou M\to\psi_{\tau,1}\Fou M)$ (\cf \cite[Lemme 5.1.4]{Saito-1}), and using the previous diagram, this is equivalent to $\phi_{\tau,1}\Fou M\simeq\Im(\Fou \wt\rN:\phi_{\tau,1}\Fou\wt M\to\phi_{\tau,1}\Fou\wt M)$. Taking the inverse isomorphism \eqref{eq:**} gives the assertion.\qed

\subsubsection*{Purity of $\gr_\ell^{\rM(\wt\rN)}\psi_{1/t,1}\wt M$ for $\ell\neq0$}

According to \cite[Lemme 5.1.12]{Saito-1}, $\var_\tau$ is strictly compatible with the monodromy filtration after a shift by $-1$. Using the previous commutative diagram, we conclude that the same property holds for the inclusion $\psi_{1/t,1}M\to\psi_{1/t,1}\wt M$. On the other hand, $\rM_{-1}(\wt\rN)\subset\Im\wt\rN=\psi_{1/t,1}M$. Therefore, the previous inclusion induces isomorphisms
\begin{equation}\label{eq:gr}
\gr_{\ell+1}^{\rM(\rN)}\psi_{1/t,1}M\isom\gr_\ell^{\rM(\wt\rN)}\psi_{1/t,1}\wt M
\end{equation}
for each $\ell\leqslant -1$. Remark now that such morphisms underly morphisms of mixed Hodge structures, since $\rM_\bbullet$ is a filtration by mixed Hodge structures. By strictness, the corresponding morphisms of mixed Hodge structures are isomorphisms. Since the left-hand term is pure of weight $w+\ell$, the right-hand term is so. Lastly, since $\wt\rN^\ell:\gr_\ell^{\rM(\wt\rN)}\psi_{1/t,1}\wt M\isom\gr_{-\ell}^{\rM(\wt\rN)}\psi_{1/t,1}\wt M(-\ell)$ is an isomorphism of mixed Hodge structures, we conclude that $\gr_\ell^{\rM(\wt\rN)}\psi_{1/t,1}\wt M$ is pure of weight $w+\ell$ for $\ell\geqslant1$.\qed

\subsubsection*{Dimension of $\gr_{w+\ell}^W\psi_{1/t,1}\wt M$}

We now consider the weight filtration $W_\bbullet\psi_{1/t,1}\wt M$ of the mixed Hodge structure $\psi_{1/t,1}\wt M$. We claim that
\begin{equation}\label{eq:claim}
\forall\ell,\quad\dim\gr_{w+\ell}^W\psi_{1/t,1}\wt M=\dim\gr_\ell^{\rM(\wt\rN)}\psi_{1/t,1}\wt M.
\end{equation}
Notice that, since both filtrations are exhaustive, it is enough to prove the claim for $\ell\neq0$. Assume first that $\ell\leqslant-1$. On the one hand, we have by \eqref{eq:gr}
\[
\dim\gr_\ell^{\rM(\wt\rN)}\psi_{1/t,1}\wt M=\dim\gr_{\ell+1}^{\rM(\rN)}\psi_{1/t,1}M=\dim\gr_{w+\ell}^W\psi_{1/t,1}M.
\]
On the other hand, since $W_\bbullet\psi_{1/t,1}M''=L_\bbullet\psi_{1/t,1}M''$ and $L_{w+\ell}\psi_{1/t,1}M''=0$ for $\ell\leqslant-1$, the natural morphism
\begin{equation}\label{eq:W}
W_{w+\ell}\psi_{1/t,1}M\to W_{w+\ell}\psi_{1/t,1}\wt M
\end{equation}
is an isomorphism, hence the assertion for $\ell\leqslant-1$. Assume now that $\ell\geqslant1$. We have
\begin{align*}
\dim\gr_{w+\ell}^W\psi_{1/t,1}\wt M&=\dim\gr_{w+\ell}^W\psi_{1/t,1}M+\dim\gr_{w+\ell}^W\psi_{1/t,1}M''\\
&=\dim\gr_{\ell+1}^{\rM(\rN)}\psi_{1/t,1}M+\dim\rP\gr_{\ell-1}^{\rM(\rN)}\psi_{1/t,1}M\\
&=\dim\gr_{\ell-1}^{\rM(\rN)}\psi_{1/t,1}M.
\end{align*}
On the other hand,
\begin{align*}
\dim\gr_\ell^{\rM(\wt\rN)}\psi_{1/t,1}\wt M&=\dim\gr_{-\ell}^{\rM(\wt\rN)}\psi_{1/t,1}\wt M\\
&=\dim\gr_{-\ell+1}^{\rM(\rN)}\psi_{1/t,1}M=\dim\gr_{\ell-1}^{\rM(\rN)}\psi_{1/t,1}M,
\end{align*}
so \eqref{eq:claim} is proved.\qed

\subsubsection*{End of the proof of Proposition \ref{prop:W=M}}
The purity of $\gr_\ell^{\rM(\wt\rN)}\psi_{1/t,1}\wt M$ shows that $W_\bbullet$ induces the trivial filtration with one jump from $w+\ell-1$ to $w+\ell$ on $\gr_\ell^{\rM(\wt\rN)}\psi_{1/t,1}\wt M$, for $\ell\neq0$. In particular, for $\ell\neq0$,
\begin{equation}\label{eq:2b}W_{w+\ell-1}\cap \rM(\wt\rN)_\ell\subset \rM(\wt\rN)_{\ell-1}.
\end{equation}

Let $\ell_o\gg0$ be such that $W_{w+\ell_o}=\rM(\wt\rN)_{\ell_o}=\psi_{1/t,1}\wt M$. Then \eqref{eq:2b} shows that $W_{w+\ell_o-1}\subset\rM(\wt\rN)_{\ell_o-1}$, and \eqref{eq:claim} for $\ell=\ell_o$ implies equality. A similar argument can be applied by decreasing induction up to $\ell=1$, giving $W_{w+\ell}=\rM(\wt\rN)_{\ell}$ for any $\ell\geqslant0$. Assume now that $\ell\leqslant-1$. Then \eqref{eq:W} shows that $W_{w+\ell}=\rM(\wt\rN_{|\Im\wt\rN})_{\ell+1}$. It is easy to check that this is nothing but $\rM(\wt\rN)_\ell$. \qed
\end{proof}

\subsection{End of the proof of Proposition \ref{prop:main}}
Recall that we assume that $U=\CC^n$. Let us first show that, if we set $M=\cH^0F_*(j_{!*}\cO_U^\rH)$, which is a pure Hodge module of weight $n$, according to M.~Saito \cite{Saito-1}, we have $\cH^0f_*\cO_U^\rH\simeq\wt M$ as a mixed Hodge module. Indeed, by functoriality of the $\wt{\phantom{M}}$ operation, we have commutative diagram in $\MHM(\Afu_t)$:
\[
\xymatrix@R=.7cm{
M\ar[r]^-a\ar[d]&\cH^0f_*\cO_U^\rH\ar[d]^b\\
\wt M\ar[r]^-{\wt a}&\wt{\cH^0f_*\cO_U^\rH}
}
\]
Since the kernel and the cokernel of $a$ are constant, $\wt a$ is an isomorphism, and since the operator $\partial_t$ is invertible on $\cH^0f_*\cO_U$, $b$ is an isomorphism. As a consequence, Proposition \ref{prop:W=M} applies to $\cH^0f_*\cO_U^\rH$.\qed




\begin{thebibliography}{99}

\bibitem{A'Campo}
A'Campo, N. ``La fonction z{\^ e}ta d'une monodromie", \textit{Comment. Math. Helv.}, 50 (1975): 233-248.

\bibitem{Broughton}
Broughton, S. A. ``Milnor numbers and the topology of polynomial hypersurfaces", \textit{Invent. Math.}, 92 (1988): 217-241.

\bibitem{Danilov-2}
Danilov, V. I. ``The geometry of toric varieties", \textit{Russ. Math. Surveys}, 33 (1978): 97-154.


\bibitem{D-K}
Danilov, V. I. and Khovanskii, A. G. ``Newton polyhedra and an algorithm for computing Hodge-Deligne numbers", \textit{Math. Ussr Izvestiya}, 29 (1987): 279-298.

\bibitem{D-L-1}
Denef, J. and Loeser, F. ``Motivic Igusa zeta functions", \textit{J. Alg. Geom.}, 7 (1998): 505-537.

\bibitem{D-L-2}
Denef, J. and Loeser, F. ``Geometry on arc spaces of algebraic varieties", \textit{Progr. Math.}, 201 (2001): 327-348.

\bibitem{Dimca2}
Dimca, A. ``Monodromy at infinity for polynomials in two variables", \textit{J. Alg. Geom.}, 7 (1998): 771-779.

\bibitem{Dimca}
Dimca, A. \textit{Sheaves in topology}, Universitext, Springer-Verlag, Berlin, 2004.

\bibitem{D-N}
Dimca, A. and N{\'e}methi, A. ``On the monodromy of complex polynomials", \textit{Duke Math. J.}, 108 (2001): 199-209.

\bibitem{D-S-1}
Dimca, A. and Saito, M. ``Monodromy at infinity and the weights of cohomology", \textit{Compositio Math.}, 138 (2003): 55-71.

\bibitem{D-S-new}
Dimca, A. and Saito, M. ``Weight filtration of the limit mixed Hodge structure at infinity for tame polynomials", arXiv:1110.4840v1.



\bibitem{E-T}
Esterov, A. and Takeuchi, K. ``Motivic Milnor fibers over complete intersection varieties and their virtual Betti numbers", arXiv:1009.0230, to appear in \textit{Int. Math. Res. Not.}

\bibitem{Fulton}
Fulton, W. \textit{Introduction to toric varieties}, Princeton University Press, 1993.

\bibitem{L-N-2}
Garc{\'i}a L{\'o}pez, R. and N{\'e}methi, A. ``Hodge numbers attached to a polynomial map", \textit{Ann. Inst. Fourier}, 49 (1999): 1547-1579.

\bibitem{G-L-M}
Guibert, G., Loeser, F. and Merle, M. ``Iterated vanishing cycles, convolution, and a motivic analogue of a conjecture of Steenbrink", \textit{Duke Math. J.}, 132 (2006): 409-457.

\bibitem{H-T-T}
Hotta, R., Takeuchi, K. and Tanisaki, T. \textit{D-modules, perverse sheaves, and representation theory}, Birkh{\"a}user Boston, 2008.

\bibitem{K-S}
Kashiwara, M. and Schapira, P. \textit{Sheaves on manifolds}, Springer-Verlag, 1990.

\bibitem{Khovanskii}
Khovanskii, A. G. ``Newton polyhedra and toroidal varieties", \textit{Funct. Anal. Appl.}, 11 (1978): 289-296.

\bibitem{Kushnirenko}
Kouchnirenko, A. G. ``Poly\'edres de Newton et nombres de Milnor", \textit{Invent. Math.}, 32 (1976): 1-31.


\bibitem{Le}
L{\^ e}, D. T. ``Some remarks on relative monodromy", \textit{Real and complex singularities (Proc. Ninth Nordic Summer School/NAVF Sympos. Math., Oslo, 1976)} (1977): 397-403.

\bibitem{L-S}
Libgober, A. and Sperber, S. ``On the zeta function of monodromy of a polynomial map", \textit{Compositio Math.}, 95 (1995): 287-307.

\bibitem{Looijenga}
Looijenga, E. ``Motivic measures", \textit{Ast{\'e}risque} 276 (2002): 267-297.

\bibitem{Macdonald}
Macdonald, I. G. ``Polynomials associated with finite cell-complexes", \textit{J. London Math. Soc.}, 4 (1971): 181-192.

\bibitem{M-T-new2}
Matsui, Y. and Takeuchi, K. ``Milnor fibers over singular toric varieties and nearby cycle sheaves", \textit{Tohoku Math. J.}, 63 (2011): 113-136.

\bibitem{M-T-new3}
Matsui, Y. and Takeuchi,  K. ``Monodromy zeta functions at infinity, Newton polyhedra and constructible sheaves", \textit{Mathematische Zeitschrift}, 268 (2011): 409-439.

\bibitem{M-T-new1}
Matsui, Y. and Takeuchi, K. ``A geometric degree formula for $A$-discriminants and Euler obstructions of toric varieties", \textit{Adv. in Math.}, 226 (2011): 2040-2064.

\bibitem{M-T-1}
Matsui, Y. and Takeuchi, K. ``Motivic Milnor fibers and Jordan 
normal forms of Milnor monodromies", submitted.


\bibitem{M-T-2}
Matsui, Y. and Takeuchi, K. ``On the sizes of the Jordan blocks 
of monodromies at infinity", submitted.


\bibitem{M-S86b}
Mebkhout, Z. and Sabbah, C. ``{\S\kern .15em III.4 {$\mathscr{D}$}-modules et cycles \'evanescents}, Le formalisme des six op\'erations de Grothendieck pour les {$\mathscr{D}$}-modules coh\'erents", \textit{Travaux en cours} 35, Hermann, Paris, (1989): 201-239.

\bibitem{Milnor}
Milnor, J. \textit{Singular points of complex hypersurfaces}, Princeton University Press, 1968.

\bibitem{N-N}
Neumann, W. D. and Norbury, P. ``Vanishing cycles and monodromy of complex polynomials", \textit{Duke Math. J.}, 101 (2000): 487-497.

\bibitem{Oda}
Oda, T. \textit{Convex bodies and algebraic geometry. An introduction to the theory of toric varieties}, Springer-Verlag, 1988.

\bibitem{Oka}
Oka, M. \textit{Non-degenerate complete intersection singularity}, Hermann, Paris (1997).



\bibitem{Raibaut}
Raibaut, M. ``Fibre de Milnor motivique {\`a} l'infini", \textit{C. R. Acad. Sci. Paris S{\'er}. I Math.}, 348 (2010): 419-422.

\bibitem{Sabbah}
Sabbah, C. ``Monodromy at infinity and Fourier transform", \textit{Publ. Res. Inst. Math. Sci.}, 33 (1997): 643-685.


\bibitem{Bibi05b}
Sabbah, C. ``Monodromy at infinity and Fourier~transform~II", \textit{Publ. Res. Inst. Math. Sci.}, 42 (2006): 803-835.

\bibitem{Sabbah-2}
Sabbah, C. ``Hypergeometric periods for a tame polynomial", \textit{Port. Math.}, 63 (2006): 173-226.

\bibitem{Sabbah-3}
Sabbah, C. ``Fourier-Laplace transform of a variation of polarized complex Hodge structure", \textit{J. Reine Angew. Math.}, 621 (2008): 123-158.

\bibitem{Saito-1}
Saito, M. ``Modules de Hodge polarisables", \textit{Publ. Res. Inst. Math. Sci.}, 24 (1988): 849-995.

\bibitem{Saito-3}
Saito, M. ``Exponents and Newton polyhedra of isolated hypersurface singularities", \textit{Math. Ann.}, 281 (1988): 411-417.

\bibitem{Saito-2}
Saito, M. ``Mixed Hodge modules", \textit{Publ. Res. Inst. Math. Sci.}, 26 (1990): 221-333.

\bibitem{S-T-1}
Siersma, D. and Tib{\u a}r, M. ``Singularities at infinity and their vanishing cycles", \textit{Duke Math. J.}, 80 (1995): 771-783.

\bibitem{Steenbrink}
Steenbrink, J. H. M. ``Mixed Hodge structures on the vanishing cohomology", \textit{Real and Complex Singularities}, Sijthoff and Noordhoff, Alphen aan den Rijn, (1977): 525-563.

\bibitem{S-Z}
Steenbrink, J. H. M. and Zucker, S. ``Variation of mixed Hodge structure I", \textit{Invent. Math.}, 80 (1985): 489-542.

\bibitem{Takeuchi}
Takeuchi, K. ``Perverse sheaves and Milnor fibers over singular varieties", \textit{Adv. Stud. Pure Math.}, 46 (2007): 211-222.

\bibitem{Takeuchi-2}
Takeuchi, K. ``Monodromy at infinity of $A$-hypergeometric functions and toric compactifications", \textit{Math. Ann.}, 348 (2010): 815-831.


\bibitem{Varchenko}
Varchenko, A. N. ``Zeta-function of monodromy and Newton's diagram", \textit{Invent. Math.}, 37 (1976): 253-262.

\bibitem{K-V}
Varchenko, A. N. and Khovanskii, A. G. ``Asymptotic behavior of integrals over vanishing cycles and the Newton polyhedron", \textit{Dokl. Akad. Nauk SSSR}, 283 (1985): 521-525.

\bibitem{Voisin}
Voisin, C. \textit{Hodge theory and complex algebraic geometry, I}, Cambridge University Press, 2007.

\end{thebibliography}
\end{document}